\documentclass{jhrs}

\newcommand{\f}{\mathrm{f}}

\newcommand{\ZZ}{\mathbb{Z}}

\def\Groups{\mathsf{Groups}}
\def\Groupoids{\mathsf{Groupoids}}
\def\A{\mathsf{A}}
\def\B{\mathsf{B}}
\def\C{\mathsf{C}}
\def\LL{\mathsf{L}}
\def\Gla{\mathsf{GL}}

\usepackage{epsfig,enumerate,color}
\usepackage{amsmath}
\usepackage{amssymb}
\usepackage{amsfonts}
\usepackage{amscd}
\usepackage[all]{xy}
\usepackage{euscript}
\usepackage{latexsym}
\usepackage{multicol}
\usepackage[numbers,sort&compress]{natbib}
\usepackage{geometry}

\parskip\medskipamount

\def\leq{\leqslant}
\def\geq{\geqslant}
\newtheorem{example}{Example}[section]

\newtheorem{examples}[example]{Examples}
\newtheorem{Def}[example]{Definition}

\newtheorem{question}[example]{Problem/Question}
\newtheorem{questions}[example]{Problems/Questions}
\newtheorem{prop}[example]{Proposition}
\newtheorem{thm}[example]{Theorem}
\newtheorem{remarks}[example]{Remarks}

\newtheorem{rem}[example]{Remark}

\newtheorem{cor}[example]{Corollary}
\newtheorem{lem}[example]{Lemma}

\def\epsilon{\varepsilon}
\def\Act{\operatorname{Act}}
\def\St{\operatorname{St}}
\def\Gl{\operatorname{GL}}
\def\E{\operatorname{E}}
\def\Ob{\operatorname{Ob}}
\def\Arr{\operatorname{Arr}}
\def\Ker{\operatorname{Ker}}
\def\Aut{\operatorname{Aut}}
\def\endbox{\rule{0em}{1ex}\hfill $\Box$}
\def\Equiv{\operatorname{Equiv}}

\def\<{\langle}
\def\>{\rangle}

\begin{document}

\title{Global Actions, Groupoid Atlases and Related Topics\\
(\emph{dedicated to the memory of Saunders MacLane}(1909-2005)}

\author{A. Bak}
\email{bak@mathematik.uni-bielefeld.de}       
\address{Department of Mathematics,\\
University of Bielefeld,\\
P. O. Box 100131,\\
33501 Bielefeld,\\
GERMANY}

\author{R. Brown}            
\email{r.brown@bangor.ac.uk }       
\address{Department of Mathematics,\\
         University of Wales, Bangor,\\
         Dean Street,\\
         Bangor,\\
         Gwynedd LL57 1UT\\
         UK. }

\author{G. Minian}            
\email{gminian@dm.uba.ar}       
\address{Departamento de Matem\'{a}tica,\\
Facultad de Ciencias Exactas y Naturales,\\
Universidad de Buenos Aires,\\ Argentina}

\author{T. Porter}            
\email{t.porter@bangor.ac.uk}       
\address{Department of Mathematics,\\
         University of Wales, Bangor,\\
         Dean Street,\\
         Bangor,\\                 
         Gwynedd LL57 1UT\\
         UK. } 


\begin{abstract}
A. Bak  developed  a combinatorial approach to higher $K$-theory,
in which control is kept of the elementary operations involved,
through paths and `paths of paths' in what he called a {\it global
action}. The homotopy theory of these was developed by G. Minian.
R. Brown and T. Porter developed applications to identities among
relations for groups, and also the extension to {\it
groupoid atlases}. This paper is intended as an  introduction to
this circle of ideas, and so to give a basis for exploration and
development of this area.
\end{abstract}

\received{February 9, 2006}   
\revised{June 2, 2006}    
\published{June 12, 2006}  
\submitted{Hvedri Inassaridze}  

\volumeyear{2006} 
\volumenumber{1}  
\issuenumber{1}   

\startpage{101}     

\maketitle

\tableofcontents

\section{Introduction}
The motivation for the  introduction of global actions by A. Bak
\cite{bak1,bak2} was to provide an algebraic setting in which to
bring higher algebraic $K$-theory nearer to the intuitions of the
original work of J.H.C. Whitehead on $K_1(R)$. In this work,
elementary matrices, and sequences of their actions on the general
linear group, play a key r\^ole.

The geometric origin for $K_1$ came from Whitehead's plan to seek
a  generalisation to all dimensions of the Tietze equivalence
theorem in combinatorial group theory: this theorem states that
two finite presentations of isomorphic groups may be transformed
from one to the other by a finite sequence of elementary moves,
called Tietze transformations. These algebraic moves were
translated by Whitehead into the geometric moves of {\it
elementary collapses and expansions} of finite simplicial
complexes: this gave the notion of {\it simple homotopy
equivalence}. The astonishing conclusion of his work was that
there was an obstruction to a homotopy equivalence $f: X \to Y$ of
finite simplicial complexes being a simple homotopy equivalence,
and that this lay in a quotient group of a group he defined,
namely $K_1(R)$ where $R$ is the integral group ring of
$\pi_1(X)$.

Since then, there was a search for higher order groups $K_i(R)$,
which were finally defined by Quillen as homotopy groups of a
space $F(R)$, the homotopy fibre of a map $B\Gl(R) \to B\Gl(R)^+$.
This was a great result, but this excursion into topology has
meant that the original combinatorial intuitions get somewhat
lost.

Bak, in \cite{bak1, bak2}, for instance, found that to deal with stability questions and higher
$K$-groups it was useful to consider a family of elementary
matrices, and that the theory in general could be seen as a family
of group actions,  indexed by a set with a relation $\leq$, often
a partial order, and with certain `patching conditions'. This
became his `global action', which was viewed as a kind of
`algebraic manifold', analogous to
the notion of topological manifold, but in an algebraic setting. A global action was, by analogy with the atlases of the theory of manifolds, an \emph{atlas of actions}. Thus `local' meant at one
group action, and `global' meant understanding the interaction of
them all. The elaboration of the definition was intended to cope
with paths, to define an analogue of $K_1$, and paths of paths, to
deal with higher order questions.

In discussions at Bangor and Bielefeld, it was seen that: (a)
there were interesting applications of global actions to
identities among relations for groups with a specified family of subgroups,
and (b) there were advantages in using the well known transition
from group actions to groupoids given by the action groupoid of an
action, to rephrase the definition of global action so that it became part of a slightly wider concept, \emph{atlases of groupoids}, or as we will nearly always say,
\emph{groupoid atlases}. The individual groupoids would play the  analogous role to the \emph{charts} in the definition of a manifold.  This generalisation would allow a wider scope for the theory,
since groupoids can generalise not only groups and group actions
but also equivalence relations. Further, the notion of a
presentation of a groupoid allows for potential extensions of
these ideas to resolutions, and so to higher dimensions.

There are even further possibilities.

The intention is to loosen the restriction to having algebraic
control of a system to be specified in terms of {\it globally defined}
algebraic structures, but to have also the possibility of a {\it
locally defined} algebraic structure. (One could imagine that this idea might be important in some areas of
physics, where immediate information, such as symmetry,  is given
locally.)

It is interesting to compare this idea with other tools for
local-to-global problems involving `patching', of which the most
notable and very well explored is the notion of sheaf. An
advantage of global actions and groupoid atlases over other
previously defined tools is the possibility of incorporating in
the theory the notion of {\it elementary process}, and paths and
homotopies of paths. This allows concrete notions of higher
homotopies to be studied and used.

The notes on which this paper is based were  produced by Tim
Porter, on the basis of lectures by Bak, together with discussions in Bangor
and Bielefeld. A version has been available since then as a Bangor
Preprint \cite{bakbmp99}.
\footnote{Acknowledgements:  The research discussions summarised here
have been supported by the British Council/ARC programme (`Global
actions and homotopy theory' Project Number 859) and supplemented
by an INTAS grant (Algebraic K-theory, groups and categories,
INTAS 93-436). Brown was also supported by a Leverhulme Emeritus
Fellowship, 2002-2004.}

\section{Global actions}

The motivating idea is of a family of interacting and overlapping
local $G$-sets for varying $G$. The prime example is the
underlying set $\Gl _n (R)$ operated on by the family of subgroups
$\Gl _n (R)_{\alpha}$ which are generated by elementary matrices
of a certain form $\alpha$. We will give the details of this
example shortly, but first we will set up some basic terminology and notation concerning group actions ($G$-sets), before we give the definition of a global action.

 A (left) group action consists of a group $G$ and a set $X$ on which $G$ acts on the left; we will write $G\curvearrowright X$. This means that there is a function from the set $G\times X$ to $X$ written as  $(g,x)$ goes to $g.x$, such that $g_1.(g_2.x) = (g_1g_2).x$ and $1_G.x =x$ for all $g_1,g_2\in G$ and $x\in X$.
 
 It is often convenient to omit the dot so we may write $gx$ instead of $g.x$.
 
 A morphism of group actions, $(\varphi, \psi) : G \curvearrowright X \rightarrow
H\curvearrowright Y$, consists of a homomorphism of groups $\varphi :G \rightarrow H$ and a function $\psi : X\rightarrow Y$ such that $\psi (g.x) = \varphi(g).\psi(x)$.

The promised `global version' of this is:
\begin{Def}
A {\it global action } $\A$ consists of a set $X_\A$ together
with:
\begin{enumerate}[(i)]
\item an indexing set $\Phi_\A$, called the {\it coordinate
system} of $\A$; \item a reflexive relation, written $\leq$, on
$\Phi_\A$;\item a family $\{(G_\A)_{\alpha}\curvearrowright
(X_\A)_\alpha \mid \alpha \in \Phi_\A\} $ of group actions on
subsets  $ (X_\A)_{\alpha} \subseteq X_\A $; the $(G_\A)_{\alpha}$
are called the {\it local groups} of the global action;  \item for
each pair  $\alpha \leq \beta$ in $\Phi_\A$,  a group morphism
\[
 ({G}_\A)_{\alpha \leq \beta}: ({G}_\A)_{\alpha}\to (G_\A)_{\beta}
 .
\]
  \item[] This data is required to satisfy: \item if $\alpha \leq
\beta$ in $\Phi_\A$, then $  ({G}_\A)_{\alpha \leq \beta}$ leaves
$(X_\A)_{\alpha}
  \cap (X_\A)_{\beta}$ invariant;
  \item if $\sigma \in (G_\A)_\alpha$ and $x \in (X_\A)_{\alpha}
  \cap (X_\A)_{\beta}$, then \begin{equation}
   \sigma x = (({G}_\A)_{\alpha \leq
  \beta}(\sigma))x. \tag*{$\Box$}\end{equation}
\end{enumerate}
\end{Def}
The diagram $G_\A : \Phi_\A \rightarrow  \Groups$ is called the
{\it global group} of $\A$. The set $X_\A$ is the {\it enveloping
set} or {\it   underlying set} of $\A$. The notation $|X_\A|$ or
$|\A|$ for $X_\A$ is sometimes used for emphasis or to avoid
confusion since
$$X_\A : \Phi_\A \rightarrow  \mathcal{P} (X_\A)$$
is also a useful notation, where $ \mathcal{P} (X_\A)$ is the
powerset of $X_\A$.

\begin{remarks}
\hspace{0.5em}
\begin{enumerate}[a)]

\item For technical reasons it is not assumed that the collection
$(X_\A)_{\alpha} \subseteq X_\A$ necessarily  covers $X_\A$.  This
holds in all the basic  examples we will examine but is not a
requirement.

\item The relation $\leq$ is not assumed to be transitive
on $\Phi_\A$, so really $G_\A$ is not a functor. However, the
difference is minor as, if $F(\Phi_\A)$ denotes the free category
on the graph of $(\Phi_\A, \leq)$, then $G_\A$ extends to a
functor $G_\A : F(\Phi_\A) \rightarrow \Groups$.  We will
usually refer, as here, to $G_\A$ as a diagram of groups and
will sometimes use `natural transformation' to mean a generalised
natural transformation defined on the generating graphs, which
yields an actual natural transformation on the corresponding
extensions. It will sometimes be useful to consider groups as single object groupoids, in which case the above yields a diagram of groupoids\footnote{By a groupoid we mean a small
category in which every arrow is an isomorphism.}.\end{enumerate}

\end{remarks}

The simplest global actions come with just a single domain: a
global action $\A$ is said to be {\it single domain} if for each
$\alpha \in \Phi_\A, (X_\A)_{\alpha} = |\A|$.

\begin{example}
Let $G$ be a group, $\mathcal{H} = \{H_i \mid i \in \Phi\}$ a
family of subgroups of $G$. For the moment $\Phi$ is just a set
(that is : $\alpha \leq \beta$ in $\Phi$ if and only if $\alpha =
\beta$). Define $\A=\A(G, \mathcal{H})$ to be the global action
with
\begin{align*}
X = |X_\A| & = |G|,\text{ the underlying set of } G\\
\Phi_\A & = \Phi\\
(X_\A)_{\alpha} & = X_\A \text{ for all } \alpha \in \Phi
\\H_i &\curvearrowright X  \text{ by left multiplication}\end{align*}
(so the local orbits of the $H_i$-action are the left cosets of
$H_i$).

Later on in section \ref{subdivAGH}, we will need to refine this construction, \label{secondAGH}taking $\Phi_\A
$ to be the family of finite non-empty subsets of $\Phi$ ordered
by opposite inclusion and with $(G_\A)_\alpha = \bigcap_{i \in
\alpha}H_i$ if $\alpha \in \Phi_\A$.\endbox\end{example}

We will later look in some detail at certain specific such single
domain global actions. The following prime motivating example is
similar to these, but the indexing set/coordinate system is
slightly more complex.

\begin{example}{\bf The General Linear Global Action  $\Gla _n (R)$.}\label{GLn}
Let $R$ be an associative ring with identity and $n$ a positive
integer.

Let $\Delta = \{ (i, j) ~|~ i \neq j, 1 \leq i, j \leq n\}$ be the
set of non-diagonal positions in an $n \times n$ array. Call a
subset $\alpha \subseteq \Delta$ \emph{closed} if
$$(i, j) \in \alpha \text{ and } (j, k) \in \alpha \text{ implies } (i,
k) \in \alpha$$

Note that if $(i, j) \in \alpha$ and $\alpha$ is closed then $(j,
i) \notin \alpha$.

Let $\Phi = \{ \alpha \subseteq \Delta \mid \alpha \text{ is closed}\}$. We
put a reflexive relation $\leq$ on $\Phi$ by $\alpha \leq \beta$ if
$\alpha \subseteq \beta$.

Now suppose $(i, j) \in \Delta$ and $r \in R$. The {\it elementary
matrix} $\epsilon_{ij} (r)$ is the matrix obtained from the
identity $n \times n$ matrix by putting the element $r$ in
position $(i, j)$,
\[
\text{i.e.} \qquad \epsilon_{ij} (r)_{k, l} = \begin{cases}
1 & \text{if } k = l\\
r & \text{if } (k, l) = (i, j)\\
0 & \text{otherwise .}
\end{cases}
\]
Let $\Gl _n (R)_{\alpha}$, for $\alpha \in \Phi$, denote the
subgroup of $\Gl _n (R)$ generated by $$\{ \epsilon_{ij} (r) \mid
(i, j) \in \alpha, r \in R\}.$$ It is easy to see that $(a_{kl})
\in \Gl _n (R)_{\alpha}$ if and only if
\[
a_{k, l} =  \begin{cases}
1 & \text{if } k = l\\
\text{arbitrary} & \text{if } (i, j) \in \alpha\\
0 & \text{if } (i, j) \in \Delta \backslash \alpha.
\end{cases}
\]
For $\alpha \leq \beta$, there is an inclusion of $\Gl
_n(R)_{\alpha}$ into $\Gl _n(R)_{\beta}$. This will give the
homomorphism
$$\Gl _n(R)_{\alpha \leq \beta} : \Gl _n (R)_{\alpha} \rightarrow \Gl _n
(R)_{\beta}.$$ Let $\Gl _n(R)_{\alpha}$ act by left multiplication
on $\Gl _n(R)$.

This completes the description of the single domain global action
$\Gla_n(R)$. Later we will see how to define the fundamental group and more generally the higher homotopy
groups of a global action. The $(i - 1)^{\text{th}}$-homotopy
group of $\Gla_n(R)$ is the algebraic $K$-theory group $K_i
(n, R)$ and the usual algebraic $K$-group, $K_i(R)$ is the direct
limit of $K_i (n, R)s$ by the obvious maps induced from the
inclusions $\Gl_n (R) \rightarrow  \Gl_{n + 1}
(R)$. \endbox
\end{example}

The way that a global action extends local information to become
global information can be observed from the simplest cases of the
$\A (G, \mathcal{H})$.

If $\mathcal{H}$ has just a single group $H$ in it, then the global
action is just the collection of orbits, i.e. right cosets. There is no
interaction between them.

If $\mathcal{H}$ consists of distinct subgroups $ \{H_1, H_2\}$,
then any $H_1$-orbit intersects with any $H_2$-orbit, so now
orbits do interact. How they interact can be very  influential on
the homotopy properties of the situation. 

\begin{example}\label{S3}
As a simple example consider
the symmetric group $S_3 \equiv \langle a, b ~|~ a^3 = b^2 =
(ab)^2 = 1 \rangle$, with $a$ denoting the 3-cycle (1~2~3) and $b$
the transposition (1 ~2). Take $H_1 = \langle a\rangle = \{ 1,
(1~2~3), (1~3~2) \}$ yielding two orbits for its left action on
$S_3, H_1 \text{ and } H_1 b$. Similarly take $H_2 = \langle
b\rangle$ giving  local orbits $H_2, H_2 a, H_2a^2$. Any
$H_1$-orbit intersects with any $H_2$-orbit, but of course they do
not overlap themselves. This gives an intersection diagram:
\[
\xymatrix{
&H_1 \ar@{-}[ld] \ar@{-}[rd] \ar@{-}[rrrd]&  & H_1 b \ar@{-}[llld]
\ar@{-}[ld]\ar@{-}[rd] &\\
H_2 & & H_2 a &  & H_2 a^2
}
\]
This graph makes it clear that, even in such a simple case, it is
possible to find loops and circuits within the global action,
following an element through a local orbit and within an intersection
crossing to the next orbit, eventually getting back to the starting
position.

The element $1 \in H_2$ multiplied on the left by $b \in H_2$ ends up
in $H_2 \cap H_1 b$, multiplied on the left by $a \in H_1$ yields $ab
\in H_1b \cap H_2 a^2$ and so on. The circuit
\[
\xymatrix{
H_2 \ar@{>}[r]&  H_1 b \ar@{>}[r] &H_2 a^2 \ar@{>}[r] & H_1 \ar@{>}[r]& H_2\\
I \ar@{>}[r]^{b\times}& b \ar@{>}[r]^{a\times}& ab \ar@{>}[r]^{b\times}& bab \ar@{>}[r]^{a\times}& abab = 1
}
\]
relates the  structure of the single domain global action with the
combinatorial information encoded in the presentation. This will be
examined in more detail later.\end{example}

\section{Morphisms}

Morphisms between  global actions come in various
strengths depending  on what part of the data is preserved. Preservation
of just the local orbit information corresponds to a ``morphism'',
compatibility  with the whole of the data then yields a ``regular
morphism''.

First we introduce a subsidiary notion which will be important at
several points in the later development.

\begin{Def}
Let $\A$ be a global action. Let $x \in (X_\A)_{\alpha}$ be some
point in a local set of $\A$.

 A \emph{local frame} at $x$ in $\alpha$ or $\alpha$\emph{-frame} at $x$ is a
sequence $x = x_0, \cdots, x_p$ of points in the local orbit of
the $(G_\A)_{\alpha}$-action on $(X_\A)_{\alpha}$ determined by
$x$. Thus for each $i$, $1 \leq i \leq p$, there is some $g_i \in
(G_\A)_{\alpha}$ with $g_i x = x_i$.  \endbox
\end{Def}
Note that in extreme cases, such as a trivial action, all the
$x_i$ may be equal, but if the action is faithful, each
$\alpha$-frame at $x$  consists essentially of $x$ and a sequence
$g_1, \cdots, g_p$ of elements of $(G_\A)_{\alpha}$. For some of
the homotopy theoretic side of the development this may be of use
as $g_1, g_2g_1^{- 1}, \cdots$ yields a $(p - 1)$-simplex in the
nerve of the group $(G_\A)_{\alpha}$.

\begin{Def}
If $\A$ and $\B $ are global actions, a {\it morphism} $f : \A
\rightarrow \B $ of global actions is a function $f: |\A|
\rightarrow |\B |$ on their
underlying sets, which preserves local frames. More precisely:\\
if $x_0, \cdots, x_p$ is an $\alpha$-frame at $x_0$ for some
$\alpha \in \Phi_\A$ then $f(x_0), \cdots, f(x_p)$ is a
$\beta$-frame at $f(x_0)$ for some $\beta \in \Phi_\B $.  \endbox
\end{Def}
Note that not all $\alpha$-frames may lead to the same $\beta$, so
this notion is {\it not} saying that the whole of the local orbit
of the $(G_\A)_{\alpha}$-action corresponding to $x_0$ must end up
within a single local orbit, merely that given $x_0, \cdots, x_p$,
there is some $\beta$ such that $f(x_0), \cdots, f(x_p)$ form a
$\beta$-frame. This is of course only significant when there are
infinitely many co-ordinates, as larger frames may lead to
different ``larger'' $\beta$s.

Intuitively a \emph{path} in a global action $\A$ is a sequence of
points $a_0, \cdots, a_n$ in $|\A|$ so that each $a_i, a_{i + 1},
i = 0, \cdots, n-1$ is a $\alpha$-frame for some (varying) $\alpha
\in \Phi_\A$. This idea can be captured using a morphism from a
global action model of a line, and this is done in the initial
papers on global actions, \cite{bak1,bak2}. Here we postpone this
until the third section as there is a certain technical advantage
in considering the line with a groupoid atlas structure and that
will be introduced there.

\begin{Def}
A {\it regular morphism} $\eta : \A \rightarrow \B $ of global
actions is a triple $(\eta_{\Phi}, \eta_G, \eta_X)$ satisfying the
following
\begin{align*}
-  ~\eta_{\Phi} : &  ~\Phi_\A \rightarrow \Phi_\B  \text{ is a
relation preserving
function} :\\
& \text{if } \alpha \leq \alpha', \text{ then } \eta_{\Phi} (\alpha) \leq
\eta_{\Phi} (\alpha').\\
- ~\eta_G : &~ G_\A \rightarrow (G_\A)_{\eta\Phi(~)} \text{ is a
natural transformation of group diagrams over
$\eta_{\Phi}$,}\end{align*} i.e. for each $\alpha \in \Phi_\A$,
$$\eta_G (\alpha) : (G_\A)_{\alpha} \rightarrow (G_\B )_{\eta_{\Phi} (\alpha)}$$
is a group homomorphism such that if $\alpha \leq \alpha^\prime$
in $\Phi_\A$, the diagram
\[
\xymatrix{ (G_\A)_{\alpha} \ar@{>}[d] \ar@{>}[r]^{\eta_G (\alpha)}
&
(G_\B )_{\eta_{\Phi} (\alpha)} \ar@{>}[d]\\
(G_\A)_{\alpha'} \ar@{>}[r]_{\eta_G (\alpha')} & (G_\B
)_{\eta_{\Phi} (\alpha')} }
\]
where the vertical maps are the structure maps of the respective
diagrams;\\
\hspace*{1cm}- ~$\eta_X : |\A|  \rightarrow |\B |$ is a function
such that $\eta_X ((X_\A)_{\alpha}) \subseteq (X_\B )_{\eta_{\Phi}
(\alpha)}$
for all $\alpha \in \Phi_\A$;\\
\hspace*{1cm}- for each $\alpha \in \Phi_\A$, the pair
$$ (\eta_G, \eta_X) : (G_\A)_{\alpha} \curvearrowright (X_\A)_{\alpha}
\rightarrow (G_\B )_{\eta_{\Phi} (\alpha)} \curvearrowright (X_\B
)_{\eta_{\Phi} (\alpha)}$$ is a morphism of group actions.
\endbox
\end{Def}
\begin{rem}
If $\eta$ is a regular morphism, it is clear that $\eta_X$
preserves local frames and so is a morphism in the weaker sense.

Composition of both types of morphism is defined in the obvious way
and so one obtains categories of global actions and morphisms and
of global actions and regular morphisms.

It is perhaps necessary to underline the meaning of a morphism of
group actions: If $G \curvearrowright X$ and $H \curvearrowright Y$
are group actions of $G$ on $X$ and $H$ on $Y$, respectively, a
morphism from $G \curvearrowright X$ to $H \curvearrowright Y$ is a
pair $(\varphi : G \rightarrow H, \psi : X \rightarrow Y)$ with
$\varphi$ a homomorphism and $\psi$ a function such that for $g \in G,
x \in X$,
$$\varphi (g). \psi (x) = \psi(g . x).$$
We need to note that, if $x$ and $x'$ are in the same orbit of $G
\curvearrowright X$ then $\psi(x)$ and $\psi(x')$ are in the same
orbit of $H \curvearrowright Y$.  \endbox
\end{rem}
\section{Actions as groupoids,  and groupoid atlases.}
\subsection{Actions as groupoids}
If $G \curvearrowright X$ is a group action then we can construct
an {\it action groupoid}  from it. 

$\Act(G, X)$ or $G \ltimes X$ will denote the category with $X$ as
its set of objects and $G \times X$ as its set of arrows. Given an
arrow $(g, x)$, its source is $x$ and its target $g.x$. We write
$s(g, x) = x, t(g, x) = g.x$ and represent this diagrammatically
by
\[
\xymatrix{
x \ar@{>}[r]^{(g, x)} &  g.x
}.
\]
The composite of $(g, x)$ and $(g', x')$ is defined only if the target
of $(g, x)$ is the source of $(g', x')$ so $x'= g.x$, then
\[
\xymatrix{ x \ar@{>}[r]^{(g, x)} & g.x \ar@{>}[r]^{(g', gx)} &g'
gx }\] gives a composite $(g'g, x)$. The identity at $x$ is $(1,
x)$. The inverse of $(g, x)$ is $(g^{- 1}, gx)$ so $G \ltimes X$
is in fact a groupoid.

\begin{example}
Let $X = \{0, 1\}, G = C_2$ with the obvious action on $X$
interchanging 0 and 1. If we write $C_2 = \{1, c\}$ we have $\Ob
  (G \ltimes X) = X = \{0, 1\}$,
$$
\Arr(G \ltimes X) =  \{ (1, 0) : 0 \rightarrow 0, (1, 1)
  : 1 \rightarrow 1,   \hspace{0,2cm}(c, 0) : 0 \rightarrow 1, (c, 1) : 1 \rightarrow 0 \}$$
Thus diagrammatically the groupoid is just
\[
\xymatrix{ G \ltimes X : = &  \underset{0}{.}\ar@(ul,dl)[]
\ar@/^/[rr]^{(c, 0)} &&\underset{1}{\cdot} \ar@/^/[ll]^{(c, 1)}
\ar@(ur,dr)[] }\] i.e. it is the groupoid often written as
$\mathcal{I}$, the {\it unit interval groupoid}.  \endbox
\end{example}

Back to the general situation:

Suppose $(\varphi, \psi) : G \curvearrowright X \rightarrow
H\curvearrowright Y$ is a morphism of group actions, then we can
define a morphism of groupoids by
\begin{align*}
\varphi \ltimes \psi : G \ltimes X & \rightarrow  H \ltimes Y\\
(\varphi \ltimes \psi) (x) & = \psi(x) \quad\hspace{1.3cm} \text{on objects}\\
(\varphi \ltimes \psi) (g, x) & =  (\varphi (g), \psi (x)) \quad
\text{ on arrows}.
\end{align*}

We check:
\begin{align*}
s(\varphi (g), \psi (x)) & = \psi(x) = \psi (s(g, x)),\\
t(\varphi (g), \psi (x)) & = \varphi (g).\psi (x) = \psi (g,
x)\\
& = \psi t (g, x)
\end{align*}
so $\varphi \ltimes \psi$ preserves source and target. It also
preserves identities and composition as is easily checked.
\subsection{Groupoid atlases}
The ``language'' of group actions thus translates well into the
language of groupoids. The notion of an orbit of a group action
becomes a connected component of a groupoid, so what is the analogue
of a global action? The translation is not difficult, but the obvious
term ``global groupoid'' does not seem to give the right intuition
about the concept. We noted that a global action was similar to the notion of an atlas in the theory of manifolds, so is an \emph{atlas of actions}, so instead we will use the term `groupoid atlas' or, synonymously, `atlas of groupoids'.

First a bit of notation: if $G$ is a groupoid with object set $X$ and
$X' \subset X$ is a subset of $X$ then $G\downharpoonright_{X'}$ will denote the
groupoid with object set $X'$ having
$$G\downharpoonright_{X'} (x, y) = G (x, y),$$
if $x, y \in X'$ and with the same composition and identities as
$G$, when this makes sense. This groupoid
$G\downharpoonright_{X'}$ is the {\it full
  sub-groupoid of } $G$ {\it determined by the objects in } $X'$ or
more simply, the {\it restriction} of $G$ to $X'$.
\begin{Def}
A {\it groupoid atlas } $\A$ consists of a set $X_\A$ together
with:
\begin{enumerate}[(i)]
\item an indexing set $\Phi_\A$, called the {\it coordinate
system} of $\A$; 
\item a reflexive relation, written $\leq$, on
$\Phi_\A$;
\item a family $\mathcal{G}_\A = \{(G_\A)_{\alpha} \mid \alpha \in
\Phi_\A\} $ of groupoids with object sets $ (X_\A)_{\alpha}$; the
$(G_\A)_{\alpha}$ are called the {\it local groupoids} of the
groupoid atlas; \item if $\alpha \leq \beta$ in $\Phi_\A$, a
groupoid morphism
\[
 ({G}_\A)_{\alpha} \downharpoonright_{(X_\A)_{\alpha}
\cap (X_\A)_{\beta}}  \longrightarrow
({G}_\A)_{\beta}\downharpoonright_{(X_\A)_{\alpha} \cap
(X_\A)_{\beta}}
\]
which is the identity map on objects. The notation we will use for
this morphism will usually be $\varphi_{\beta}^{\alpha}$ but the
more detailed $({ G}_\A)_{\alpha \leq
  \beta}$ may be used where more precision is needed. As before we
write $|\A| \text{ for } X_\A$, the underlying set of $\A$.
\item[] This data is required to satisfy: \item if $\alpha \leq
\beta$ in $\Phi_\A$, then $(X_\A)_{\alpha}
  \cap (X_\A)_{\beta}$ is a union of components of $({G}_\A)_{\alpha}$,\\
~~i.e. if $x \in (X_\A)_{\alpha} \cap (X_\A)_{\beta}$ and $g \in
({G}_\A)_{\alpha}$ is such that $s (g) = x$ then $t(g) \in
(X_\A)_{\alpha} \cap (X_\A)_{\beta}$.  \endbox
\end{enumerate}
\end{Def}

A morphism of groupoid atlases comes in several strengths as with
the special case of global actions.

A \emph{local frame} in a groupoid atlas, $\A$, is a sequence $(x_0,
\cdots, x_p)$ of objcts in a single connected component of some
$(\mathcal{G}_\A)_{\alpha}$, i.e. there is some $\alpha \in
\Phi_\A$, $x_0, \cdots, x_p \in (X_\A)_{\alpha}$ and arrows $g_i :
x_0 \rightarrow x_i, i = 1, \cdots, p$.

A function $f : |\A| \rightarrow |\B |$ supports a \emph{weak morphism}
structure if it preserves local frames. Similar comments apply to
those made above about morphisms of global actions.

The stronger form of morphism of groupoid atlases will just be called
a (strong) morphism.

A \emph{strong morphism} $\eta : \A \rightarrow \B $ of groupoid
atlases is a
triple $(\eta_X,\eta_{\Phi}, \eta_g)$ satisfying the following
\begin{enumerate}[(i)]
\item $\eta_X : X_A\rightarrow X_B$ is a function between the underlying
sets;
\item $\eta_{\Phi} : \Phi_\A \rightarrow \Phi_\B $ is a relation
preserving
function;
\item $\eta_G :\mathcal{G}_\A \rightarrow (\mathcal{G}_\B
)_{\eta_{\Phi}}$ is a (generalised) natural transformation of
diagrams of groupoids {\it over} the function $\eta_{\Phi}$ 
on the objects.
\end{enumerate}

To illustrate the difference between global actions and groupoid
atlases, we consider some simple examples.

\begin{example}
Let $X = \{0, 1, 2\}, G = C_3 = \{1, a, a^2\}$ (and, of course,
$a^3 = 1$), the cyclic group of order 3, acting by $a. 0 = 1, a. 1
= 2$, on $X$. This gives us $C_3 \ltimes X$ with 9 arrows. We set
$\B  = C_3 \ltimes X$ as groupoid or $C_3 \curvearrowright X$ as
$C_3$-set. We also have the example $\A = C_2 \ltimes \{0, 1\} =
\mathcal{I}$ considered earlier.

Both $\A$ and $\B $ will be considered initially as global actions
having $\Phi_\A$ and $\Phi_\B $  a single element.

Any function $f : \{0, 1\}\rightarrow \{0, 1, 2\}$ supports the
structure of a morphism of global actions since the only
non-trivial frame in $\A$ is based on the set $\{x_0,x_1\}$, where $x_0 = 0$, and $x = 1$
and this must get mapped to a frame in $\B $, since any non-empty
subset of $X$ is a frame in $\B $. On the other hand, a regular
morphism $\eta: \A \rightarrow \B $ must contain the information
on a group homomorphism
$$\eta_G : C_2 \rightarrow C_3$$
which must, of course, be trivial. Hence the only regular morphism
$\eta$ must map all of $\A$ to a single point in $\B $. There are
thus 9 morphisms from $\A$ to $\B $, but only 3 regular morphisms.
The regular morphisms are very rigid.  \endbox
\end{example}
\begin{rem}
It is not always the case that there are fewer regular morphisms
than (general) morphisms.  If $\A$ is a global action with one
point and a group acting on that point and $\B $ is similar with
group $H$, there is only one general morphism from $\A$ to $\B $,
but the set of regular morphisms is `the same as' the set of group
homomorphisms from $G$ to $H$.  \endbox
\end{rem}

\begin{example}
Now we continue the previous example by considering $\A$ and $\B $
as groupoid atlases. The element $(c, 0): 0 \rightarrow 1$ in the
single groupoid determining $\A$, must be sent to some arrow in
$\B $. The inverse of $(c, 0)$ is $(c, 1)$, so as soon as a
morphism, $\eta_G$ is specified on $(c, 0)$, it is determined on
$(c, 1)$ since $\eta_G (c, 1) = (\eta_G (c, 0))^{- 1}$. Thus if we
pick an arrow in $\B $, say,
$$(a^2, 0) : 0 \rightarrow 2,$$
we can define a morphism
$$\eta_G : \A \rightarrow \B $$
by specifying $\eta_G (c, 0) = (a^2, 0)$, so  $ \eta_G (0) = 0$,
$\eta_G (1) = 2$ etc. In other words the fact that $\A$ uses an
action by $C_2$ and $\B $ by $C_3$ does not inhibit the existence
of morphisms from $\A$ to $\B $. Any morphism of global actions
from $\A$ to $\B $ {\it in this case} will support the structure
of a morphism of the corresponding groupoid atlases, yet the extra
structure of a ``regularity condition'' is supported in this
latter setting. Of course the relationship between morphisms of
global actions and morphisms of the corresponding groupoid atlases
can be expected to be more subtle in general.  \endbox
\end{example}

\begin{question}
\emph{If $\A$ and $\B $ are global actions and $f: \A \rightarrow \B $
is a morphism, does $f$ support the structure of a morphism of the
corresponding groupoid atlases?}

In general the answer is `no' since if $\A$ is a global action
with $\Phi_\A = \{a,b \mid a\leq b\}$ with both $X_a$ and $X_b$
single points, and $\B $ is similar but with $\Phi_\A$ discrete,
then the general morphism which corresponds to the identity does
not support the structure of a (strong) morphism of the
corresponding groupoid atlases because of the need for a relation
function $\eta : \Phi_\A \rightarrow \Phi_\B $. Refining the
question, suppose we have a general morphism of global actions
together with a relation preserving function between the
coordinate systems, which is compatible with the morphism. In that
case the question is related to the following question
about groupoids:\\
\emph{if we have two groupoids $\A$ and $\B $ and a function $f$ from
the objects of $\A$ to the objects of $\B $ which sends connected
components of $\A$ to connected components of $\B $, what
obstructions are there for there to exist a functor $F$ from $\A$
to $\B $ such that $F$ restricted to the objects is the given $f$?}
\end{question}

Clearly any global action determines a corresponding groupoid atlas as
we have used above. As there are morphisms of action groupoids that do
not come from regular morphisms of actions, the groupoid morphisms give a new
notion of morphism of global actions, whose usefulness for the
motivating examples will need investigating. Are there ``useful''
groupoid atlases other than those coming from global actions? The
answer is most definitely: yes.
\subsection{$\Equiv(\A)$}

Equivalence relations are examples of groupoids.
\begin{example}
Let $X$ be a set. Any equivalence relation $R$
on $X$ determines a groupoid with object set $X$. We will denote this
groupoid by $R$ as well. It is specified by
\[
R(x, y) = \begin{cases}
\{(x, y)\} &\text{ if } x Ry\\
\emptyset &\text{ if }  x \text{ is not related to } y.
\end{cases}
\]
Now suppose $R_1, \cdots, R_n$ are a family of equivalence
relations on $X$. Then define $\A$ to have coordinate system
$$\Phi_\A = \{1, \cdots, n\} \quad \text{ with discrete } \leq$$
and $(G_\A)_i = R_i$. This gives a groupoid atlas that does  not
in general arise from a global action.  \endbox
\end{example}
\begin{example}
Let $G$ be a group, $X$ a $G$-set and $R$ an equivalence relation
on $X$. Let $\Phi = \{1, 2\}$, with $\leq$ still to be specified.
Take $G_1 = G \ltimes X$, $G_2 = R$ and $X_1 = X_2 = X$. Assume we
have a groupoid atlas structure with this as partial data. If
$\leq$ is discrete, there is no interaction between the two
structures and no compatibility requirement. If $1 \leq 2$, each
$G$-orbit is contained in an equivalence class with $\varphi_2^1
(x, g) = (x, gx)$, i.e. the $G$-orbit structure is finer than the
partition into equivalence classes. If $2 \leq 1$, the partition
is finer than the orbit structure (the connected components of the
groupoid $G_1$) and if $x Ry$ then there is some $g_{x, y} \in G$
such that $g_{x, y} x = y$.  \endbox
\end{example}
This last case is closely related to a useful construction on
global actions.

\begin{example}
Let $\A = (X_\A, G_\A, \Phi_\A)$ be a global action. Let $\alpha
\in \Phi_\A$ and $(G_\A)_{\alpha} \curvearrowright
(X_\A)_{\alpha}$ be the corresponding action. Set $R_{\alpha}$ to
be the equivalence relation determined by the
$(G_\A)_{\alpha}$-action. Thus $x R_{\alpha} x'$ if and only if
there is some $g \in (G_\A)_{\alpha}$ with $g x = x'$. Of course
the partition of $(X_\A)_{\alpha}$ into $R_{\alpha}$-equivalence
classes is exactly that given by the $(G_\A)_{\alpha}$-orbits (or
the $(\mathcal{G}_\A)_{\alpha}$-components where
$(\mathcal{G}_\A)_{\alpha}$ is the corresponding groupoid).

If $\alpha \leq \beta$ then the compatibility conditions are
satisfied between $R_{\alpha}$ and $R_{\beta}$ making $(X_\A,
R_\A, \Phi_\A)$ with $R_\A = \{ R_{\alpha} \mid \alpha \in \Phi_\A\}$
into a groupoid atlas which will be denoted $\Equiv(\A)$.

The functions $(\mathcal{G}_\A)_{\alpha} \rightarrow R_{\alpha}$
mapping the groupoid of the $(G_\A)_{\alpha}$-action to the
corresponding equivalence relation yield a natural transformation
of groupoid diagrams and hence a strong morphism
$$\A \rightarrow \Equiv(\A)$$
with obvious universal properties. Of course the same construction
works if $\A$ is an arbitrary groupoid atlas, that is, one not
necessarily arising from a global action. The result gives a left
adjoint to the inclusion of the full subcategory of atlases of
equivalence relations into that of groupoid atlases. The
usefulness of this construction is another reason for extending
our view beyond global actions to include groupoid atlases. The
notion of morphism of global actions, $f :\A \rightarrow \B $,
translates to the notion of strong morphism, $f : \Equiv(\A)
\rightarrow \Equiv(\B )$ of the corresponding groupoid atlases, at
least for examples with finite orbits.  \endbox
\end{example}

\subsection{The Line}

We have seen that the simple action with $G = C_2$, $X = \{0, 1\}$
gives the groupoid $\mathcal{I}$ (also sometimes written [1] as it
is the groupoid version of the 1-simplex). We want an analogue of
a line so as to describe paths and loops. The line, $\LL$, is
obtained by placing infinitely many copies of $\mathcal{I}$ end to
end. It is a global action, but, as the morphisms that give paths
in a global action $\A$ will need to be non-regular morphisms in
general, it is often expedient to think of it as a groupoid atlas.

The set, $|\LL|$, of points of $\LL$ is $  \ZZ$, the set of
integers; $\Phi_\LL = \ZZ \cup \{\square\}$, where $\square$ is an
element satisfying $\square < n$ for all $n \in \ZZ$, and
otherwise the relation $\leq$ is equality. (Thus $\square \leq
\square$, for all $n\in \ZZ$, $\square < n$ and $n \leq n$, but
that gives all related pairs.) If $n \in \Phi_\LL$, $(X_\LL)_n =
\{n, n + 1\}$, whilst $(X_\LL)_{\square} = |\LL|$ itself.

The groupoid $(\mathcal{G}_\LL)_n$ is a copy of $\mathcal{I}$,
whilst $(\mathcal{G}_\LL)_{\square}$ is discrete with trivial
vertex groups.

The underlying structure of $\LL$ rests firmly on the locally
finite simplicial complex structure of the ordinary real line.
There the (abstract) simplicial complex structure is given by:
\begin{align*}
\text{Vertices } & = \ZZ, \text{ the set of integers};\\
\text{Set of 1-simplices } & = \{ \{n, n + 1\} \mid n \in \ZZ \}, \text{ the set of adjacent pairs in } \ZZ.
\end{align*}
We will see shortly that there is a close link between simplicial
complexes and this context of global actions/ groupoid atlases.

\section{Curves, paths and connected components}

Suppose $\A$ is a global action or more generally a groupoid
atlas. A curve in $\A$ is simply a (weak) morphism
$$f: \LL \rightarrow \A$$
where $\LL$ is the line groupoid atlas introduced above.

This implies that $f: |\LL| \rightarrow |\A|$ is a function for
which local frames are preserved. In $\LL$ the local frames are
simply the adjacent pairs $\{n, n + 1\}$ and the singleton sets
$\{n\}$. Thus the condition that $f: \LL \rightarrow \A$ be a path
is that the sequence of points
$$\cdots, f(n), f(n + 1), \cdots$$
is such that for each $n$, there is a $\beta \in \Phi_\A$ and
$g_{\beta} : f (n) \rightarrow f(n + 1)$ in
$(\mathcal{G}_\A)_{\beta}$. (If you prefer global action notation
$g_{\beta} \in (G_\A)_{\beta}$ and $g_{\beta} f(n) = f(n + 1)$.)

Note that $f$ does not specify $\beta$ and $g_{\beta}$, merely
requiring their existence. This observation leads to a notion of a
{\it strong curve} in $\A$ which is a morphism of groupoid atlases
$$f: \LL \rightarrow \A$$
so for each $n$ one gets a $\beta = \eta_{\Phi}(n)\in \Phi_\A$ and $\eta_G : \mathcal{G}_\LL \rightarrow
(\mathcal{G}_\A)_{\eta_{\Phi}}$ is a natural transformation of
groupoid diagrams. This condition only amounts to specifying
$\eta_G (n, \LL) = g :  f(n) \rightarrow f(n + 1)$, but this time
the data is part of the specification of the curve. We can thus
write a strong curve as $(\cdots, f(n), g_n, f(n + 1), \cdots )$,
that is a sequence of points of $|\A|$ together with locally
defined arrows
$$g_n : f(n) \rightarrow f(n + 1)$$
in the chosen local groupoid $(\mathcal{G}_\A)_{\beta}$. Changing
the $\beta$ or the $g_n$ changes the morphism. We will later see
the r\^{o}le of strong curves, strong paths, etc.

A ({\it free}) {\it path} in $\A$ will be a curve that stabilises
to a constant value on both its left and right ends. More
precisely it is a curve $f : \LL \rightarrow \A$ such that there
are integers $N^- \leq N^+$ with the property that
$$\text{for all } n \leq N^-,~ f(n) = f(N^-);$$
$$\text{for all } n \geq N^+,~ f(n) = f(N^+).$$

We will call $(N^-, N^+)$ a \emph{stabilisation pair} for $f$.

A ``based path'' can be defined if $\A$ has a distinguished base
point. This occurs naturally in such cases as $\A = \Gla_n (R)$ or $\A = \A(G, \mathcal{H})$ for $\mathcal{H}$ a
family of subgroups of a group $G$, but is also defined abstractly
by adding the specification of the chosen base point explicitly to
the data. This situation is well known from topology where a
notation such as $(\A, a_0)$ would be used. We will adopt similar
conventions.

If $(\A, a_0)$ is based groupoid atlas, a {\it based path} in
$(\A, a_0)$ is a free path that stabilises to $a_0$ on the left,
i.e., in the notation above, $f(N^-) = a_0$.

A {\it loop} in $(\A, a_0)$ is a based path that stabilises to
$a_0$ on both the left and the right so $f(N^-) = f(N^+) = a_0$.

The analogue in this setting of concepts such as ``connected
component'' should now be clear. We say that points $p$ and $q$ of
$\A$, a global action or groupoid atlas, are {\it free path
equivalent} if there is a free path in $\A$ which stabilises to
$p$ on the left and to $q$ on the right.

Clearly free path equivalence is reflexive. It is also symmetric
since if $g_n : f(n) \rightarrow f(n + 1)$ in a local patch then
$g^{- 1}_n : f (n + 1) \rightarrow f(n)$. Once a free path from
$p$ to $q$ has reached $q$ (i.e. has stabilised at $q$) then it
can be concatenated with a path from $q$ to $r$, say, hence free
path equivalence is also transitive. The equivalence classes for
free path equivalence will be called {\it connected components},
with $\pi_0 (\A)$ denoting the {\it   set of connected components}
of $\A$. If $\A$ has just one  connected component then it is said
to be {\it connected}.

\begin{examples}
\hspace{0.5em}

\begin{enumerate}[\hspace{0.5em} 1.]
\item The prime and motivating example is the set of connected
components $\pi_0{\Gla_n (R)}$ of the general linear global
action.

Suppose $x, y \in\Gla_n (R)$. Suppose $f: \LL\rightarrow\Gla_n
(R)$ is a free path from $x$ to $y$, so there are $N^- \leq N^+$
as above with
\begin{align*}
\text{ if } n \leq N^-, & ~f(n) = f(N^-) = x,\\
\text{ if } n \geq N^+, & ~f(n) = f(N^+) = y.
\end{align*}
 For each $i  \in [N^-, N^+]$, there is some local arrow
$$g_i : f(n) \rightarrow f(i + 1)$$
and since $\Gla_n (R)$ is a global action, this means there is
some $\alpha_i \in \Phi$ and $\varepsilon_i \in\Gla_n
(R)_{\alpha_i}$ such that $\varepsilon_i f(i) = f(i + 1)$. (The
specification of $f$ gives the {\it existence} of such an
$\varepsilon_i$ but does not actually specify which of possibly
many $\varepsilon_i s$ to take, so we choose one. In fact of course, $f(i)$, and $f(i+1)$ are invertible matrices, so there is only one $\varepsilon_i$ possible.) We thus have
$$\varepsilon_{N^+} \varepsilon_{N^+ - 1} \cdots \varepsilon_{N^-} x =
  y.$$
If $\E_n (R)$ is the subgroup of elementary matrices of $\Gl _n
(R)$, this is the subgroup generated by all the $\Gl
_n(R)_{\alpha}$ for $ \alpha \in \Phi$ and so if $x$ and $y$ are
free path equivalent $$y \in \E_n (R) x,$$ i.e., $x$ and $y$ are in
the same right coset of $\E_n (R)$.

Conversely if $y \in \E_n (R) x$, there is an element $\varepsilon \in
\E_N (R)$ such that $y = \varepsilon x$, but $\varepsilon$ can be
written (in possibly many ways) as a product of elementary matrices
$$\varepsilon = \varepsilon_N \cdots \varepsilon_1$$
with $\varepsilon_i = \Gl _n (R)_{\alpha_i}$, say. Then defining
$$f : \LL \rightarrow \Gla_n (R)$$
by
\[ f(n) = \begin{cases}
x & n \leq 0\\
\varepsilon_n \cdots \varepsilon_1 x & 1 \leq n \leq N\\
y & n \geq N
\end{cases}
\]
gives a free path from $x$ to $y$ in $\Gla_n (R)$.

Thus $\pi_0 (\Gla_n (R)) = \Gl_n (R) / \E_n (R)$, the set of right
cosets of $\Gl_n (R)$ modulo elementary matrices. This is, of
course, the algebraic $K$-group $K_1(n, R)$ if $R$ is a
commutative ring.

We can naturally ask the question: `is $K_2 (n, R) \cong \pi_1 (\Gla_n
(R))$?' even if we have not yet defined the righthand side of this.

We note the use of the strong rather than the weak version of paths would not
change the resulting $\pi_0$.

\item
 Suppose $\A = \A (G, \mathcal{H})$. Can one calculate $\pi_0 (\A)$? A
similar argument to that in 1 above shows that if $x, y \in |\A| =
|G|$, then they are free path equivalent if and only if there are
indices $\alpha_i \in\Phi$ and elements $h_{\alpha_i} \in
H_{\alpha_i}$,  such that
$$h_{\alpha_n} \cdots h_{\alpha_0} x = y$$
for some $n$. Thus writing $\langle \mathcal{H} \rangle = \langle H_i \mid i \in \Phi \rangle$ for the subgroup of $G$ generated by the family
$\mathcal{H} = \{H_i \mid i\in \Phi\}$, we clearly have
$$\pi_0 (\A(G, \mathcal{H})) = G/ \langle \mathcal{H} \rangle.$$
Again the question arises as to $\pi_1 (\A(G, \mathcal{H}))$: what
is it and what does it tell us?
\end{enumerate}
\end{examples}
\section{ Fundamental groups and fundamental \\groupoids.}
Ideas for the construction of $\pi_1 (\A, a_0)$ for a pointed
global action or groupoid atlas, $\A$, seem clear enough. There
are three possible approaches:
\begin{enumerate}[(i)]
\item take some notion of homotopy of paths and define $\pi_1
  (\A, a_0)$ to be the set of homotopy classes of loops at $a_0$;
  hopefully this would be a group for the natural notion of
  composition via concatenation of paths. Alternatively define a
  fundamental groupoid $\Pi_1 \A$ using a similar plan and then take
  $\pi_1 (\A, a_0)$ to be the vertex group of $\Pi_1 \A$ at $a_0$.
\item define a global action or groupoid atlas structure on the
  set $\Omega (\A, a_0)$ of loops at $a_0$, then take $\pi_1 (\A, a_0) =
  \pi_0 (\Omega (\A, a_0))$.
\item define covering morphisms of global actions or groupoid
  atlases, then use a universal or simply connected covering to find a
  classifying group for connected coverings. This should be $\pi_1 (\A,
  a_0)$.
\end{enumerate}

We will look at the first two of these in this section, handling the
third one later.
\subsection{Products}

We start by checking that products exist in the various categories we are looking at.  In this section we will only need them in very special cases, but they will be needed later on in full strength.

Let $\A$ and $\B $ be groupoid atlases, $\A = (X_\A,
\mathcal{G}_\A, \Phi_\A), \B  = (X_\B , G_\B , \Phi_\B )$, and consider the structure that we will denote by $\A \times \B $  and which is
 given by
$$|X_{\A \times \B }| = |X_\A| \times |X_\B |;$$
$$\Phi_{\A \times \B } = \Phi_\A \times \Phi_\B $$
with $(\alpha, \beta) \leq (\alpha', \beta')$ if and only if  $\alpha
\leq \alpha'$ and $\beta \leq \beta'$;
$$(\mathcal{G}_{\A \times \B })_{(\alpha, \beta)} =(\mathcal{G}_\A)_{\alpha} \times (\mathcal{G}_\B )_{\beta},$$
the product groupoid, for $(\alpha, \beta) \in \Phi_{\A \times \B
}$ with the obvious product homomorphisms as coordinate changes. (We thus have $\mathcal{G}_{\A\times \B}$ is the product groupoid diagram
$$\xymatrix{\Phi_{\A\times \B} = \Phi_\A \times \Phi_\B \ar[r]^{\mathcal{G}_\A\times \mathcal{G}_\B\quad}& \Groupoids \times \Groupoids \ar[r]& \Groupoids.})$$

\begin{lem}
$\A\times \B$ with this structure is a groupoid atlas. 
\end{lem} 
\begin{proof}The proof is by routine checking so is omitted.\end{proof}

There are obvious `projections' $p_\A : \A \times \B \to \A$, $p_\B : \A \times \B \to \B$, at least at the level of underlying sets, so we need to check that they enrich nicely to give the various strengths of morphism.  We start with weak morphisms.

\begin{lem}
$p_\A$ and $p_\B$ support the structure of weak morphisms.
\end{lem}
\begin{proof}
A local frame in $\A \times \B$ will be a sequence $((x_0,y_0), \ldots, (x_p,y_p))$ of objects in a connected component of some $(\mathcal{G}_{\A\times \B})_{(\alpha,\beta)}$, so there are arrows $(g_i,g^\prime_i) : (x_0,y_0)\to (x_i,y_i)$ in $(\mathcal{G}_\A)_\alpha \times (\mathcal{G}_\B)_\beta$.  Clearly this means that $(x_0,\ldots, x_p)$ is an $\alpha$-frame and $(y_0, \ldots, y_p)$ is a $\beta$-frame, so $p_\A$ and $p_\B$ are weak morphisms.
\end{proof}

This lemma  is an immediate corollary of the next one, but has the advantage that its very simple proof shows directly the structure of local frames in $\A \times \B$ and how they relate to those in $\A$ and $\B$.  This will aid our intuition when looking at the links with simplicial complexes later on.
\begin{lem}
a) $p_\A$ and $p_\B$ are strong morphisms of groupoid atlases.\\
b) If $\A$ and $\B$ are global actions, considered as groupoid atlases, then so is $\A \times \B$
\end{lem}
\begin{proof}
(  b) will be left as an `exercise', as it is easy to check.)

a)  We have already specified $p_\A$ at set level.  On the coordinate system $p_{\A,\Phi} : \Phi_{\A \times \B} \rightarrow \Phi_\A$ is again just the projection and this is clearly order preserving.  Finally $p_{\A,\mathcal{G}} : (\mathcal{G}_{\A \times \B})_{(\alpha,\beta)} \rightarrow (\mathcal{G}_{\A})_{p_{\A,\Phi}(\alpha,\beta)}$ is the projection from $(\mathcal{G}_\A)_{\alpha} \times (\mathcal{G}_\B )_{\beta}$ to its first factor. This coincides with $(p_\A)$ on objects and satisfies all the naturality conditions.
\end{proof}

We next return to weak morphisms to see if $(\A \times \B, p_\A, p_\B)$ has the universal property for products (in the relevant category of groupoid atlases and weak morphisms).

Suppose that $f :\C\to \A$ and $g: C\to \B$ are weak morphisms, then we clearly get a mapping $(f,g) : \C \to \A\times \B$ given by $(f,g)(c) = (f(c),g(c))$. This preserves local frames as is easily seen and is unique with the property that $p_\A(f,g) = f$ and $p_\B(f,g) = g$. We thus nearly have:
\begin{prop}
The categories of groupoid atlases and of global actions, both with weak morphisms, have all finite products.
\end{prop}
\begin{proof}
The only thing left to note is that these categories have a terminal object, namely the singleton trivial global action.
\end{proof}

As we would expect we have a similar result for strong morphisms.
\begin{prop}
The categories of groupoid atlases and of global actions, with strong morphisms, has all finite products.
\end{prop}
\begin{proof}
The only thing left to prove is the universal property of the `product' $\A\times \B$ in ths strong case, so suppose $f :\C\to \A$ and $g: C\to \B$ are now morphisms, then we have $f = (f_X, f_\Phi, f_\mathcal{G})$, etc.  As the construction of $\A\times \B$ uses products of sets, reflexive relations and groupoids, we get the corresponding existence and uniqueness of $(f,g) = ((f_X,g_X),(f_\Phi,g_\Phi),(\f_\mathcal{G},g_\mathcal{G}))$ more or less for free.  The terminal object is as before the singleton trivial global action.
\end{proof}
\begin{remarks}\hspace{0.5em}
\begin{enumerate}[(i)]
\item It is now easy to generalise from finite products to arbitrary products $\prod \A_i = (\prod X_{\A_i},\prod \Phi_{\A_i}, \prod \mathcal{G}_{\A_i})$.  We leave the detailed check to the reader as we will not need this result for the moment.
\item If we are looking at products of global actions, the corresponding projection morphisms are regular.\end{enumerate}
\end{remarks}
\subsection{Homotopies of paths}

We only need one example of a product for the moment, namely $\LL
\times \LL$, which is a groupoid atlas model of $\mathbb{R}^2$. Just as a path
has to start and end somewhere so does a homotopy of paths.

Given a global action $\A$, points $a, b \in \A$ and paths $f_0,
f_1 : \LL \rightarrow \A$ joining  $a$ and $b$ (and hence
stabilising to these values to the left and right respectively), a
(fixed end point) {\it   homotopy} between $f_0$ and $f_1$ is a
morphism
$$ h : \LL \times \LL \rightarrow \A$$
such that:\\
there exist $N^-, N^+ \in \ZZ$, $N^- \leq N^+$ such that
\begin{enumerate}
\item[-] for all $n \leq N^-$, and all $m \in |\LL|$, $h(m, n) =
f_0 (m);$ \item[-] for all $n \geq N^+$, and all $m \in |\LL|$,
$h(m, n) = f_1 (m);$ \item[-] for all $m \leq N^-$, and all $n \in
\ZZ$, $ h(m, n) = a;$ \item[-] for all $m \geq N^+$, and all $n
\in \ZZ$,  $ h(m, n) = b.$
\end{enumerate}

\begin{remarks}
\hspace{0.5em}

\begin{enumerate}[(i)]
\item The idea is that if we consider $\LL  \times \LL$ as being
  based on the integer lattice of the plane, the morphism $h$ must
  stabilise along all horizontal and vertical lines outside the square
  with corners $(N^+, N^+)$, $(N^-, N^+)$, $(N^-, N^-)$  and $ (N^+,
  N^-)$. Although the paths $f_0, f_1$ coming with given ``lengths'',
  i.e. a given number of steps from $a$ to $b$, we allow a homotopy to
  increase, or decrease, the number of those steps an arbitrary (finite) amount.
\item Given any $f: \LL \rightarrow \A$, a path from $a$ to $b$,
we   can re-index $f$ to get
$$f' : \LL \rightarrow \A$$
with $f'(n) = a$ if $n \leq 0$, and a new $N^{+,}$ so that $f'(n) =
b$ if $n \geq N^{+,}$, simply by taking the old stabilisation pair
$(N^-, N^+)$ for $f$ and defining
$$f'(n) = f(n - N^-), \quad n \in \ZZ = |\LL|.$$
The resulting $f'$ is homotopic to $f$. Although this is fairly clear
intuitively, it is useful as an exercise as it brings home the
complexity of the processes involved, but also their inherent
simplicity.
\end{enumerate}
\end{remarks}

\begin{example}
Let $f_0 : \LL \rightarrow \A$ be a path from $a$ to $b$ with
$(N_0^-, N^+_0)$ being a stabilisation pair for $f$. Define $f_1 :
\LL \rightarrow \A$ by
$$f_1 (n) = f_0 (n - 1),\quad \quad n \in |\LL|,$$
i.e., $f_1$ is $f_0$ shifted one ``notch'' to the right on $\LL$. Then\\
\begin{tabular}{rl}
&(i) a suitable stabilisation pair for $f_1$ is $(N_0^- + 1, N^+_0 +
1)$\\
& (ii) $f_1$ is a path from $a$ to $b$\\
and &(iii) $f_1$ is homotopic to $f_0$.
\end{tabular}\\
The only claim that is not obvious is (iii). To construct a suitable homotopy $h$, we
construct many intermediate steps. For simplicity we will start
defining $h$ on the upper half-plane (we can always extend it to a
suitable square afterwards by a vertically constant extension):
$$h (n, 0) = f_0 (n), \qquad n \in |\LL|$$
and we make a choice of a local arrow $g_n: f_0 (n) \rightarrow f_0(n +
1)$, for each $n$ (of course, for $n \leq N^-_0$ or $n \geq N^+_0$, $g_n$
will be an identity of the local groupoid patch),
\begin{eqnarray*}
h (n, 1)  = & f_0 (n) &\text{ for } n \leq N^+_0 - 1\\
h (N^+_0, 1)  = & f_0 (N^+_0 - 1) &
\end{eqnarray*}
with the identity on $f_0(N^+_0 - 1)$ as corresponding local arrow from $h(N^+_0 - 1,1)$ to $h(N^+_0, 1)$.
$$h(N^+_0 + 1, 1) = f_0 (N^+_0) \quad \text{ and stabilise
  horizontally.}$$
Thus so far we have inserted an identity one place from the end and
shifted the end stage one to the right. We give next a local arrow
from $h(n, 0)$ to $h (n, 1)$ for each $n$. For most this will be the
identity arrow but for the local arrow from  $h (N_0^+, 0)$ to
$h(N_0^+, 1)$ we take $g_{N_0^+}^{- 1}$. The same idea is used for
$h(n, 2)$ but with the identity inserted one step back to the left. At
each successive stage of the homotopy, the ``ripple'' that is the
identity moves an extra step to the left. (In the diagram we write
$N$ for $N_0^+$.)
\[
\xymatrix{
h(-, 2) : \cdots \ar@{>}[r] & f_0 (N - 2) \ar@{>}[r]^{id}
& f_0 (N - 2) \ar@{>}[r]^{g_{N - 2}}  & f_0 (N - 1) \ar@{>}[r]^{g_{N - 1}}
& f_0 (N) \ar@{>}[r] & S\\
h(-, 1) : \cdots \ar@{>}[r] & f_0 (N - 2) \ar@{>}[r]^{g_{N - 2}}
\ar@{>}[u]_{id} & f_0
(N - 1) \ar@{>}[r]^{id} \ar@{>}[u]_{\overset{- 1}{g_{N - 2}}}& f_0 (N
- 1) \ar@{>}[r]^{g_{N - 1}} \ar@{>}[u]_{id} &
f_0 (N)\ar@{>}[r] \ar@{>}[u]_{id}& S\\
h(-, 0) : \cdots \ar@{>}[r] & f_0 (N - 2)\ar@{>}[r]^{g_{N - 2}}
\ar@{>}[u]_{id} & f_0
(N - 1) \ar@{>}[r]^{g_{N - 1}} \ar@{>}[u]_{id} & f_0 (N)
\ar@{>}[r]^{id} \ar@{>}[u]_{\overset{- 1}{g_{N - 2}}}&
f_0 (N) \ar@{>}[r] \ar@{>}[u] & S
}
\]
(the symbol $S$ indicates the sequence stabilises to the last
specified value.) \end{example}

Thus within a homotopy class we can ``ripple homotopy'' a path to have
specified $N^-$ (or for that matter $N^+$).

Now suppose $f: \LL \rightarrow \A$ is a path from $a$ to $b$ and
$g : \LL \rightarrow \A$ one from $b$ to $c$. We can assume that
the stabilisation pair for $g$ is to the right of that for $f$,
i.e., if $(N^- (f), N^+ (f))$ and $(N^- (g), N^+ (g))$ are
suitable stabilisation pairs, $N^+(f) \leq N^- (g)$. Then we can
form a concatenated path: $f \ast g$ by first going along $f$
until it stabilises at $b$ then along $g$. Of course $f \ast g$
will depend on the choice of stabilisation pairs, but using
``ripple homotopies'' we can change positions of $f$ and $g$ at
will and these homotopies will be reflected by homotopies of the
corresponding $f \ast g$. We may, for instance, start $g$
immediately after $f$ stabilises to $b$. This means that the
composition is well defined on homotopy classes of paths.

Likewise using vertical composition, i.e., exchanging the roles of
horizontal and vertical on homotopies it is elementary to prove that
(fixed end point) homotopy is an equivalence relation on
paths. Reflexivity of the homotopy relation is proved by taking the
inverses of all vertical local arrows in a homotopy. To reverse paths,
$f \leadsto  f^{(r)}$, and to prove the ``reverse'' is an inverse
modulo homotopy is also simple using the move:
\[
\xymatrix{
\cdot\ar@{>}[r]^{g_k} & \cdot \ar@{>}[r]^{g_{k}^{-1}} & \cdot\\
~ \ar@{>}[u]^{id} \ar@{>}[r]^{id} & \cdot
\ar@{>}[u]_{g_k} \ar@{>}[r]^{id} & ~
\ar@{>}[u]_{id}
}
\]
followed by a ripple homotopy to move the identities to the end.
As concatenation does not require reindexation (unlike paths in
spaces where $f : [0 , 1] \rightarrow X$ uses a {\it unit length}
interval) proof of associativity is easy: one concatenates
immediately on stabilisation to get a unique chosen composite and
then associativity is assured.

This set of properties allows one to define the {\it fundamental
   groupoid} $\Pi_1 \A$ of a global action or groupoid atlas $\A$ in
  the obvious way. The objects of $\Pi_1 \A$ are the points of $|\A|$
  whilst if $a, b$ are points of $|\A|$, $\Pi_1 \A (a, b)$ will be the
  set of (fixed end point) homotopy classes of paths from $a$ to $b$
  within $\A$. Composition is by concatenation as above and ``inversion
  is by reversion'': if $w = [f]$, $w^{- 1} = [f^{(r)}]$, where $f^{(r)}$
  is the ``reverse'' of $f$.

There is a strong variant of this construction. All the homotopies
etc. used above manipulate a strong path that represents the
chosen path, i.e., we chose the local arrows $g_n: f (n)
\rightarrow f(n + 1)$ and worked with them. There is a clear
notion of (fixed end point) strong homotopy of paths and strong
homotopy classes of strong paths compose in the same way giving a
{\it strong fundamental groupoid} $\Pi_1^{\text{Str}} \A$.  
\subsection{Objects of curves, paths and loops}

We aim to define for at least a large class of groupoid atlases
(including most if not all global actions of significance for
$K$-theory), a ``loop space'' analogous to that defined for
topological spaces. We expect to be able to concatenate loops
within that structure giving some embryonic analogue of the
$H$-space structure on a loop space. More precisely given a nice
enough groupoid atlas $ \A = (X_\A, \mathcal{G}_\A, \Phi_\A)$, we
want a new groupoid atlas $\Omega \A$ and a concatenation operator
$$\Omega \A \times \Omega \A \rightarrow \Omega \A,$$
which will induce at least a monoid structure on $\pi_0 (\Omega
\A)$ and hopefully for a large class of examples will give a group
isomorphic to $\pi_1 (\A)$.

Loops are best thought of via paths and thus via curves. We thus
start by searching for a suitable structure on the set of all
curves $\text{Mor}(\LL, \A)$ in $\A$. (Remember a curve is merely
a (weak) morphism from $\LL$ to $\A$.) If $f : \LL \rightarrow \A$
is a curve of $\A, f$ will pass through a sequence of local sets.
The intuition is that a local set containing $f$ in $\text{Mor}
(\LL, \A)$ will consist of curves passing through the same local
sets $(X_\A)_{\alpha}$ in the same sequence.

To simplify notation we will write $(\LL, \A)$ for the set of
morphisms from $\LL$ to $\A$.
\begin{Def}

Given a curve $f : \LL \rightarrow \A$ in $\A$, a function $\beta
: |\LL|
\rightarrow \Phi_\A$ {\it frames} $f$ if $\beta$ is a function such that \\
(i)
for $m \in |\LL| = \ZZ,$ $ f(m) \in (X_\A)_{\beta (m)}$; \\
 (ii)  for $m \in |\LL|$, there is a $b$
  in $\Phi_\A$ with  $b \geq \beta (m)$, $b \geq \beta (m+1)$ and a $g :
 f(m) \rightarrow f(m + 1)$ in $(\mathcal{G}_\A)_b$.
\end{Def}
Thus $\beta$ picks out the local sets $(X_\A)_{\beta (m)}$ which
are to receive $f(m)$. The condition (ii) ensures that these
choices are compatible with the requirement that $f$ be a curve.
Note that there may be curves that have no framing, especially if
the pseudo-order on  $\Phi_\A$ has few related pairs, e.g. is
discrete. For instance in a single domain global action of the
form $\A(G, \mathcal{H})$ we have used a discrete order on
$\Phi_\A$. Hence the condition $b \geq \beta (m)$, $b \geq \beta
(m+1)$ must imply that $\beta(m) = \beta(m+1)$, yet in our
examples we have seen non-trivial paths going through several
local orbits.  Thus in such a case the coordinate system we will
define shortly does not cover the set of paths in $\A$ and for
this reason we will consider in detail, later on,   the  second
global action structure on such $\A$ that was mentioned at the
start of example \ref{secondAGH}.

 \begin{lem}
Let $f: \LL \rightarrow \A$ be a curve and $\beta : |\LL|
\rightarrow \Phi_\A$ a framing for $f$.

Suppose $\sigma : |\LL| \rightarrow |\prod (\mathcal{G}_\A)_{\beta
  (\cdot)}|$ is a function defining a sequence $(\sigma_m)$ of arrows
  in the local groupoids of $\A$ such that $\sigma_m \in (\mathcal{G}_\A)_{\beta
  (m)}$ and, in fact, the source sequence of $\sigma$ is the underlying
  function of $f$, i.e.,
$$s(\sigma_m) = f(m).$$
Then the sequence $(f^\prime(m))$, where  $f^\prime(m) = t (\sigma_m)$, supports a weak morphism
structure, $f' : \LL \rightarrow \A$, and $\beta$ frames $f'$ as
well.
\end{lem}
\begin{proof}
We have $b \geq \beta(m)$ and $b \geq \beta(m + 1)$ with $g : f(m)
\rightarrow f(m + 1)$ in $(\mathcal{G}_\A)_b$. We have $\sigma_m :
f(m) \rightarrow f'(m)$ and hence its inverse is in  $(\mathcal{G}_\A)_{\beta (m)}$, so
using the structural morphism
$$\varphi_b^{\beta(m)} : (\mathcal{G}_\A)_{\beta (m)} \rightarrow
(\mathcal{G}_\A)_b,$$ we get $\varphi_b^{\beta (m)} (\sigma_m) :
f(m) \rightarrow f'(m) \text{ in } (\mathcal{G}_\A)_b$. Similarly
$\varphi_b^{\beta (m)} (\sigma_{m+1}) : f(m + 1) \rightarrow f'(m
+ 1)$.

It is now clear that using the composite
$$f'(m) \rightarrow f(m) \rightarrow f(m + 1) \rightarrow f'(m + 1)$$
of these three arrows in $(\mathcal{G}_\A)_b$ shows that $\beta$
frames $f'$. As the only frames in $\LL$ have size two, i.e., are
adjacent pairs $\{m, m + 1\}$, this also shows that $f'$ is a
curve in $\A. \qquad \Box$

It should be clear what our next step will be.

Take $\Phi_{(\LL, \A)} = \{ \beta : |\LL| \rightarrow \Phi_\A ~|~
\beta \text{ frames some curve } f\}.$

The set of objects $(X_{(\LL, \A)})_{\beta}$ of
$(\mathcal{G}_{(\LL,
  \A)})_{\beta}$ will be
$$(X_{(\LL, \A)})_{\beta} = \{ f : \LL \rightarrow \A ~|~ f \in (\LL, \A),
\beta \text{ frames } f\}$$ then take $(\mathcal{G}_{(\LL,
\A)})_{\beta}$ to be the set of sequences $(\sigma_m)$ with
$\sigma_m$ an arrow in $(\mathcal{G}_\A)_{\beta (m)}$, with the
property that $(s(\sigma_m))$ supports the structure of a curve in
$(X_{(\LL, \A)})_{\beta}$ , i.e., framed by $\beta$.  The lemma
above ensures that in this case $(t(\sigma_m))$ is also in
$(X_{(\LL, \A)})_{\beta}$, and that $(\mathcal{G}_{(\LL,
\A)})_{\beta}$ is a groupoid.

We have yet to specify the relation on $\Phi_{(\LL, \A)}$, but a
component-wise definition is the obvious one to try:
$$\quad \beta \leq \beta' \text{ if and only if } \beta (m) \leq \beta' (m)
\text{ for all } m \in |\LL|.$$ As each $\beta(m) \leq \beta'(m)$
results in an induced morphism of groupoids
$$(\mathcal{G}_\A)_{\beta (x)} \downharpoonright \rightarrow (\mathcal{G}_\A)_{\beta '(x)}\downharpoonright$$
over the intersection $(X_\A)_{\beta (m)} \cap (X_\A)_{\beta'
(x)}$, there is an induced morphism of groupoids
$$(\mathcal{G}_{(\LL, \A)})_{\beta}\downharpoonright \rightarrow (\mathcal{G}_{(\LL,
  \A)})_{\beta'}\downharpoonright$$
over $(X_{(\LL, \A)})_{\beta} \cap (X_{(\LL, \A)})_{\beta '}$.

If $f \in (\LL, \A)$ is framed by both $\beta$ and $\beta '$, then
for any $\sigma : f\rightarrow f'$ in $(\mathcal{G}_{(\LL,
\A)})_{\beta}$, we have seen in the above lemma that $f'$ is framed by both
$\beta$ and $\beta' $. We thus have all the elements of the
verification of the following:
\begin{prop}
With the above notation, $\A^\LL = ((\LL, \A), \mathcal{G}_{(\LL,
\A)}, \Phi_{(\LL, \A)})$ is a groupoid atlas. If $\A$ is a global
action, then so is $\A^\LL$.\endbox\end{prop}

The only part not covered by the previous discussion is that
relating to global actions, however taking $(G_{(\LL,
\A)})_{\beta}$ to be the product of the $(G_\A)_{\beta (m)}$, one
gets an action of this on $(X_\A)_{\beta}$ giving exactly the
groupoid $(\mathcal{G}_{(\LL,\A)})_{\beta}$ of the proposition.

We thus have a groupoid atlas $\A^\LL$ of curves in $\A$, and if
$\A$ is a global action, $\A^\LL$ is one as well. Note however
that $\A^\LL$ will not usually be a single domain global action
even when $\A$ is one.  This is one of the reasons why the general
case is necessary.

Our earlier comments also show that to assume that the $(X_\A)_\alpha$ cover $X_\A$ in our original definition would have been unduly retrictive here.  Of course, our ability to study properties of elements of $X_\A$ which lie outside the coordinate patches is restricted.

To obtain a path space $P(\A)$, we merely restrict to those curves
that are paths with an adjustment made to the local groupoids to
allow for the fact that a path can be linked by an arrow to a
general curve.

More explicitly we have
\[
\begin{array}{ll}
|P(\A)| & = \text{ the set of paths } f : \LL \rightarrow \A\\
\Phi_{P(\A)} & = \{ \beta : |\LL| \rightarrow \Phi_\A ~|~ \beta
\text{
  stably frames some path } f\}\\
(X_{P(\A)})_{\beta}& = \{ f \in |P(\A)| ~|~  \beta \text{ stably
frames
  } f\}\\
(\mathcal{G}_{P(\A)})_{\beta} & = \text{ the groupoid of stable
  sequences } (\sigma_m) \text{ with } \sigma_m \text{ an arrow in }\\
& \quad (\mathcal{G}_{P (\A)})_{\beta (m)} \text{ and such that }
(s(\sigma_m)) \text{
  supports the structure of}\\
& \quad \text{a path in } (X_{P(\A)})_{\beta}, \text{
  i.e. is stably framed by } \beta.
\end{array}
\]
The references to ``stable'' in this are the needed restriction to
ensure ``paths'' not ``curves'' are involved. Recall $f : \LL
\rightarrow \A$ is a path if it is a curve and there are integers
$N^-_f \leq N^+_f$ with $f(n) = f(N^-_f)$ for $n \leq N^-_f$ and
$f(n) = f(N^+_f)$ for $n \geq N^+_f$. (The pair $(N^-_f, N^+_f)$
was earlier called a stabilisation pair.) When $\beta : |\LL|
\rightarrow \Phi_\A$ frames a path, $f$, it would clearly be
possible to have $\beta$ varying beyond the end of the ``active
interval'', $N^-_f \leq n \leq N^+_f$, but is this necessary? We
will use the term ``stable frame'' of $f$ if for some
$(N^-_{\beta}, N^+_{\beta}), \beta (n)$ is constant for smaller
and for larger $n$. We do not specify how the stabilisation pair
for $f$ is related to one for $\beta$ if at all. Similar comments
apply to stable sequences, $(\sigma_m)$. There is a stable version
of the above lemma. The proof should be clear.

\begin{lem}
If $f$ is a path, $\beta : \LL \rightarrow \Phi_\A$ a stable
framing of $f$ and $\sigma = (\sigma_m)$ a stable sequence of
arrows with $s(\sigma_n)$ the underlying function of $f$ then $f'
= t(\sigma)$ supports the structure of a path and $\beta$ stably
frames $f'$ as well.\endbox \end{lem}

Finally we will want a based path space $\Gamma(\A,a_0)$ and a ``loop space'' $\Omega \A$. We note first that if $f :
\LL \rightarrow \A$ is in $P(\A)$ then there are integers $N^-_f
\leq N^+_f$ with $f(n) = f(N^-_f)$ for $n \leq N^-_f$ and $f(n) =
f(N^+_f)$ for $n \geq N^+_f$. Define two functions
$$e^0, e^1: P(\A) \rightarrow \A$$
by $e^0 (f) = f(N^-_f), e^1 (f) = f(N^+_f)$. These are clearly
independent of the choice of stabilisation pair for $f$ used in their
definition. Clearly we expect these functions to support (weak) morphisms on
the corresponding groupoid atlases or global actions.

Suppose we have a local frame in the groupoid atlas $P\A$. Then we
have some $\beta \in \Phi_{P\A}$, and  paths $f_0, f_1, \cdots,
f_p$ in some connected component of $(\mathcal{G}_{P(\A)})_{\beta}$. Thus
$\beta$ is a stable framing of $f_0$ and we have that there exist
$\sigma^{(i)} : f_0 \rightarrow f_i$, $ \sigma^{(i)} =
\left(\sigma^{(i)}_m\right)$. By the two lemmas, $\beta$ is a
stable framing of each $f_i$ as well - with the same ``b'' for all
of them.

Each $f_i$, each $\sigma^{(i)}$ and $\beta$ itself have stabilisation
pairs so we can find $N^-$ smaller than all the left hand ends of
these and $N^+$ larger than the right hand ends. Since
$\sigma^{(i)}_m$ will be constant for $n \leq N^-$ and also for $n
\geq N^+$, we have
$$f_0 (n), f_1(n), \cdots, f_p (n)$$
is a frame in $\beta (n)$, i.e.,  $e^0 (f_0), \cdots, e^0 (f_p)$
is a frame in $\A$ and we have shown $e^0$ is a morphism. Of
course a completely similar argument shows $e^1$ is a morphism.

We now assume that a base point $a_0$ is given in $\A$. We can
take $\Gamma (\A, a_0)$ to be the global action or groupoid atlas
defined by $e^{- 1}_0 (a_0)$. As we have not yet described how to
do such a construction, we consider a more general situation.

Suppose $f : \A \rightarrow \B $ is a weak morphism of global actions
or, more generally, groupoid atlases and let $b \in \B $, then the
set $f^{- 1} (b) \subset \A$ supports the following structure
\[
\begin{array}{ll}
|f^{- 1} (b)| & = \{ \alpha \in \A : f (a) = b\}\\
\Phi_{f^{- 1} (b)} & = \{ a \in \Phi_\A ~|~ (X_\A)_{\alpha} \cap
f^{- 1}
(b) \neq \emptyset \}\\
(X_{f^{_{- 1}} (b)})_\alpha & = (X_\A)_a \cap f^{- 1} (b),
\end{array}
\]
and the local action  / local groupoid is the restriction of that
in $\A$. That this last specification works needs a bit of care.
We first look at the global action case.

If $f(a) = b$, and $g \in (G_\A)_{\alpha}$ one usually would not
expect $g.a$ to be still ``over $b$'' so one has to take
$$(G_{f^{- 1}(b)})_{\alpha} = \{ g \in (G_\A)_{\alpha} : g.f^{_{-1}} (b) = gf^{- 1} (b) \},$$
then no problem arises. For the groupoid case the corresponding
$(\mathcal{G}_{f^{- 1} (b)})_{\alpha}$ is just the full
sub-groupoid of $(\mathcal{G}_\A)_{\alpha}$ determined by the
objects $(X_{f^{- 1}
  (b)})_{\alpha}$. The induced morphisms when $\alpha \leq \alpha'$
in $\Phi_{f^{- 1} (b)}$ are now easy to handle.

We thus have a global action/groupoid atlas $\Gamma (\A, a_0)$ of
based paths in $(\A, a_0)$ and a restricted ``end point morphism''
$$e_1 : \Gamma (\A, a_0) \rightarrow \A.$$
Of course, $\Omega (\A, a_0)$, the \emph{object of loops at}
$a_0$, will be $e^{- 1}_1 (a_0)$ considered as a subobject of
$\Gamma (\A, a_0)$ or $e^{- 1}_0 (a_0) \cap e^{- 1}_1 (a_0)$ when
thought of as being within $P(\A)$.

The discussion of concatenation of paths in the lead up to the
fundamental groupoid indicates that there is a concatenation of loops,
but that such an operation depends strongly on the choice of
stabilisation pairs for the individual loops. There is thus no single
composition map
$$\Omega(\A, a_0) \times \Omega (\A, a_0) \rightarrow \Omega(\A, a_0)$$
that is ``best possible'' or ``most natural''. Composition can be
defined if stabilisation pairs are chosen:

Given $f, g \in |\Omega (\A, a_0)|$ with chosen stabilisation
pairs $(N^-_f, N^+_f), (N^-_g, N^+_g)$ then define
$$f \ast g : \LL \rightarrow \Omega(\A, a_0)$$
\[
\text{by } \qquad (f \ast g) (n) =
\begin{cases}
f (n) \hspace{2.3cm} \text{if } n \leq N^+_f\\
g(n - N^+_f + N^-_g) \text{ if } n \geq N^+_f
\end{cases}
\]
The obvious stabilisation pair is $(N^-_f, N^+_f + N^+_g)$ and
with this choice we get an associative composition,
\underline{but} on loops with chosen stabilisation pairs. Of
course the composition is not really well defined on the loops
alone. As before it is well defined up to ``ripple homotopies''.
This structure is analogous to the `Moore loop space' construction
in a topological setting.

Our next task is to calculate $\pi_0 \Omega (\A, a_0)$ comparing
it with $\pi_1 (\A, a_0)$ as defined as the vertex group at $a_0$
of $\Pi_1 \A$, the fundamental groupoid of $\A$.

Suppose $f_0$ and $f_1$ are loops at $a_0$. When are they free
path equivalent as points of $\Omega (\A, a_0)$? Suppose $h : \LL
\rightarrow \Omega (\A, a_0)$ is a free path in $\Omega (\A, a_0)$
that stabilises to $f_0$ on the left and to $f_1$ on the right.
The ``obvious'' way to proceed is to try to use $h$ to construct a
homotopy between $f_0$ and $f_1$ as paths in $\A$. Picking a
stabilisation pair $(N^-_h, N^+_h)$ for $h$, we, of course, have
\begin{align*}
h(n) =  & ~ f_0\hspace{2.5cm} \text{ for } n \leq N^-_h,\\
h(n) =  & ~ f_1\hspace{2.5cm} \text{ for } n \geq N^+_h.
\end{align*}
Define \begin{align*}
H : |\LL| \times |\LL| & \rightarrow \A \text{ by}\\
H (m, n) &=  ~h (n) (m),
\end{align*}
and, noting that there are finitely many \underline{different}
$h(n)$ involved, pick for each of these $h(n)$ a suitable
stabilisation pair $(N^-_{h(n)}), (N^+_{h(n)})$ and set
\begin{align*}
N^- & = \min \left(\{ N^-_{h(n)} : n \in [ N^-_h, N^+_h]\}, N^-_h
\right)\\
N^+ & = \max \left(\{ N^+_{h(n)} : n \in [ N^-_h, N^+_h]\}, N^+_h \right)
\end{align*}
We claim that outside the square $[N^-, N^+] \times [N^-, N^+],$
$H(m,n) = a_0$, since, for instance, if $m < N^-$ and $n> N^+$,
then  $n > N^+_h$ so $h(n) = f_1$ and as $m < N^-_{h(n)}$, $h(n) (m)
= f_1 (m) = a_0$, as required. The other cases are similar.

Now assume that $f_0$ and $f_1$ are loops at $a_0$ in $\A$ which
yield the same element of the fundamental group $\pi_1 (\A, a_0)$.
Then there is a fixed end point homotopy
$$H : f_0 \simeq f_1$$
between them and reversing the above process, we obtain a free
path in $\Omega (\A, a_0)$. We have only to note that the
definitions of the concatenation operations in $\pi_1 (\A, a_0)$
and $\pi_0 (\Omega (\A, a_0))$ correspond exactly, to conclude
that
$$\pi_1 (\A, a_0) \cong \pi_0 (\Omega (\A, a_0)).$$

\section{Simplicial complexes from global actions.}
If $\A$ is a global action,  a path $f$ in $\A$ is a sequence of
points $ ~\cdots, f(n), \cdots~$ so that pairs of successive points are in
the orbit of some local action. A path can thus wander from one
local ``patch'' to the next by going via a point in their
intersection. It is only that each $f(n)$ is in  a local orbit
with $f (n - 1)$ and in a local orbit with $f(n + 1)$, not the
group elements used that matter. This suggests that $f$ yields a
path in a combinatorially defined simplicial complex constructed
by considering finite families of points within local orbits. This
is the case. We will describe this for a general groupoid atlas,
$\A$.

The simplicial complex $V(\A)$ of $\A$ has $|\A|$ as its set of
vertices. A subset
$$\sigma = \{x_0, \cdots, x_p\}$$
is a $p$-simplex of $V(\A)$ if there is some $\alpha \in\Phi_\A$
for which $\sigma$ is an $\alpha$-frame, i.e., $\sigma \subseteq
(X_\A)_{\alpha}$ and there are $g_i: x_0 \rightarrow x_i \in
(\mathcal{G}_\A)_{\alpha}$, so $\sigma$ is contained in a single
connected component of $(\mathcal{G}_\A)_{\alpha}$.

It is clear that $\Pi_1 V(\A) \cong \Pi_1 \A$ and if $a_0 \in \A$
is a base point then $$\pi_1 (V(\A), \{a_0\}) \cong \pi_1 (\A,
a_0).$$

The one disadvantage of $V(\A)$ is that it has as many vertices as
$\A$ has elements and this can obscure the essential combinatorial
structure involved. This construction of $V(\A)$ is an analogue of
the Vietoris construction used in Alexander-\v{C}ech cohomology
theory in algebraic topology. It is an instance of a general
construction that associates two simplicial complexes to a
relation from one set to another. This construction was studied by
C. H. Dowker, \cite{Dowker}. We outline his results.
\subsection{The two nerves of a relation: Dowker's construction}
Let $X, Y$ be sets and $R$ a relation between $X$ and $Y$, so $R
\subseteqq X \times Y$. We write $x Ry$ for $(x, y) \in R$.

For our case of interest $X$ is the set of points of $\A$ and $Y$
is the set of local components of the local groupoids if $\A$ is
groupoid atlas and is thus the set of local orbits in the global
action case. The relation is `$x Ry$ if and only if $x \in y$'.
Two other exemplary cases should be mentioned.
\begin{example}
Let $X$ be a set (usually a topological space) and $Y$ be a
collection of subsets of $X$ covering $X$, i.e. $\bigcup Y = X$.
The classical case is when $Y$ is an index set for an open cover
of $X$. The relation is the same as above i.e. $x Ry$ if and only
if $x \in y$ or more exactly $x$ is in the subset indexed by $y$.
\endbox
\end{example}
\begin{example}
If $K$ is a simplicial complex, its structure is specified by a
collection of non-empty finite subsets of its set of vertices
namely those sets of vertices declared to be simplices. This
collection of simplices is supposed to be downward closed, i.e.,
if $\sigma$ is a simplex and $\tau \subseteqq \sigma$ with $\tau \neq \emptyset$, then $\tau$
is a simplex. For our purposes here, set $X = V_K$ to be the set
of vertices of $K$ and $Y = S_K$, the set of simplices of $K$ with
$x Ry$ if $x$ is a vertex of the simplex $y$. \endbox
\end{example}
Returning to the general situation we define two simplicial
complexes associated to $R$, as follows:
\begin{enumerate}
\item[(i)] $K = K_R:$
\begin{enumerate}
\item[-]the set of vertices is the set $X$;
\item[-] a $p$-simplex of $K$ is a set $\{x_0, \cdots, x_p\}
  \subseteq X$ such that there is some $y \in Y$ with $x_i Ry$ for $i
  = 0, 1,~ \cdots, ~p$.
\end{enumerate}
\item[(ii)] $L = L_R :$
\begin{enumerate}
\item[-]the set of vertices is the set, $Y$;
\item[-]  $p$-simplex of $K$ is a set $\{y_0, \cdots, y_p\}
  \subseteq Y$ such that there is some $x \in X$ with $x Ry_j$ for $j
  = 0, 1,~ \cdots,~ p$.
\end{enumerate}
\end{enumerate}
Clearly the two constructions are in some sense dual to each
other. For our situation of global actions/groupoid atlases, $K_R$
is $V(\A)$. The corresponding $L_R$ does not yet seem to have been
considered in exactly this context. We will denote it $N(\A)$ so
if $\sigma \in N(\A)_p$
$$\sigma = \left\{ U_{\alpha_0}, \cdots, U_{\alpha_p}\right\}$$
with $U_{\alpha_i}$ a local orbit for $G_{\alpha_i} \curvearrowright
X_{\alpha_i}$ or a connected component of $\mathcal{G}_{\alpha_i}$
with the requirement that $\bigcap \sigma = \underset{i = 0}{\bigcap}^p
U_{\alpha_i} \neq \emptyset$.

In the case of $X$ a space with $Y$ an open cover, $K_R$ is the
\emph{Vietoris complex} of $X$ relative to $Y$ whilst $L_R$ is the nerve of
the open cover (often called the \emph{\u{C}ech complex} of $X$ relative to
the cover). We will consider the other example in detail later on.

\subsection{Barycentric  subdivisions}

Combinatorially, if $K$ is a simplicial complex with vertex set $V_K$,
then one associates to $K$ the partially ordered set of its
simplices. Explicitly we write $S_K$ for the set of simplices of $K$
and $(S_K, \subseteq)$ for the partially ordered set with $\subseteq$
being the obvious inclusion. The barycentric subdivision, $K'$, of $K$
has $S_K$ as its set of vertices and a finite set of vertices of $K'$
(i.e. simplices of $K$) is a simplex of $K'$ if it can be totally
ordered by inclusion. (Thus $K'$ is the simplicial complex given by
taking the nerve of the poset, $(S_K, \subseteq)$. We may sometimes
write $Sd(K)$ instead of $K'$.)

\begin{rem}
It is important to note that there is in general no natural simplicial
map from $K'$ to $K$. If however $V_K$ is ordered in such a way that
the vertices of any simplex in $K$ are totally ordered (for instance
by picking a total order on $V_K$), then one can easily specify a map
$$ \varphi : K' \rightarrow K$$
by:\\
if $\sigma' = \{ x_0, \cdots, x_p\}$ is a vertex of $K'$ (so $\sigma'
\in S_K$), let $\varphi \sigma'$ be the least vertex of $\sigma'$ in
the given fixed order.

This preserves simplices, but reverses order so
if $\sigma_1' \subset \sigma_2'$ then $\varphi (\sigma_1') \geq
\varphi(\sigma_2')$.
\end{rem}

If one changes the order, then the resulting map is
\emph{contiguous}:\\ Let $\varphi,~ \psi : K \rightarrow L$ be
two simplicial maps between simplicial complexes. They are \emph{contiguous}
if for any simplex $\sigma$ of $K$, $\varphi (\sigma) \cup \psi
(\sigma)$ forms a simplex in $L$.\\
Contiguity gives a constructive form of homotopy applicable to
simplicial maps.

If $\psi : K \rightarrow L$ is a simplicial map, then it induces
$\psi' : K' \rightarrow L'$ after subdivision. As there is no way
of knowing/picking compatible orders on $V_K$ and $V_L$ in
advance, we get that on constructing
$$\varphi_K : K' \rightarrow K$$
and
$$\varphi_L : L' \rightarrow L$$
that $\varphi_L \psi'$ and $\psi \varphi$ will be contiguous to each
other but rarely equal.

Returning to $K_R$ and $L_R$, we order the elements of $X$ and
$Y$. Then suppose $y'$ is a vertex of $L'_R$, so $y' = \{y_0,
\cdots, y_p\}$, a simplex of $L_R$ and there is an element $x \in
X$ with $x Ry_i, i = 0, 1, \cdots, p$. Set $\psi y' = x$ for one
such $x$.

If $\sigma = \{ y'_0, \cdots, y'_q\}$ is a $q$-simplex of $L'_R$,
assume $y'_0$ is its least vertex (in the inclusion ordering)
$$\varphi_L (y^\prime_0) \in y'_0 \subset y^\prime \text{ for each } y_i \in
\sigma$$
hence $\psi y'_i R \varphi_L (y'_0)$ and the elements $\psi y'_0,
\cdots, \psi y'_q$ form a simplex in $K_R$, so $\psi : L'_R
\rightarrow K_R$ is a simplicial map. It, of course, depends on the
ordering used and on the choice of $x$,  but any other choice $\bar x$
for $\psi y'$ gives a contiguous map.

Reversing the r\^{o}les of $X$ and $Y$ in the above we get a
simplicial map
$$\bar\psi : K'_R \rightarrow L_R.$$
Applying barycentric subdivisions again gives
$$\bar\psi' : K''_R \rightarrow L'_R$$
and composing with $\psi : L'_R \rightarrow K_R$ gives a map
$$\psi \bar\psi' : K''_R \rightarrow K_R.$$
Of course, there is also a map
$$\varphi_K \varphi'_K : K''_R \rightarrow K_R.$$

\medskip

 \begin{prop} (Dowker, \cite{Dowker} p.88). The two maps $\varphi_K
\varphi'_K$ and $\psi \bar\psi'$ are contiguous.\endbox
\end{prop}

\medskip

Before proving this, note that contiguity implies homotopy and that
$\varphi \varphi'$ is homotopic to the identity map on $K_R$ after
realisation, i.e., this shows that
\begin{cor}
$$|K_R| \simeq |L_R|.$$\endbox\end{cor}
The homotopy depends on the ordering of the vertices and so is not
natural. 

\begin{proof} of Proposition.

Let $\sigma''' = \{ x_0'', x_1'', \cdots, x_q''\}$ be a simplex of
$K_R''$ and as usual assume $x_0''$ is its least vertex, then for all
$i > 0$
$$x_0'' \subset x_i''.$$
We have that $\varphi'_K$ is clearly order reversing so $\varphi'_K
x_i'' \subseteq \varphi'_K x_0''$. Let $y = \bar\varphi \varphi'_K
x_0''$, then for each $x \in \varphi'_K x_0''$, $x Ry$. Since $\varphi_K
\varphi'_K x_i'' \in \varphi'_K x_i'' \subseteq \varphi'_K x_0''$, we
have $\varphi_K \varphi'_K x_i'' Ry$.

For each vertex $x'$ of $x_i'', \bar \psi x' \in \bar \psi'
x_i''$, hence as $\varphi'_K x_0'' \in x_0'' \subset x_i'', y =
\bar \psi \varphi'_K xx_0'' \in \bar \psi' x_i''$ for each $x_i
''$, so for each $x_i'', \psi \bar \psi' x_i'' Ry$, however we
therefore have $$ \qquad \varphi_k \varphi'_K (\sigma '') \cup
\psi \bar \psi (\sigma''') = \bigcup \varphi_k \varphi'_K (x_i'')
\cup \psi \bar \psi; x_i ''$$ forms a simplex in $K_R$, i.e.
$\varphi_K \varphi'_K$ and $\psi \bar \psi'$ are contiguous.
\end{proof}
\begin{example}
To illustrate both the Proposition and the remaining example of $K_R$ and
$L_R$, consider the simplicial complex, $K$:
\[
\xymatrix{
& 2 \ar@{-}[rd] \ar@{-}[dd] \ar@{-}[ld]&\\
1 \ar@{-}[rd]&  & 4 \ar@{-}[ld]\\
& 3 &
}
\]consisting of two 2-simplices joined along a common edge.  More precisely, take $X = V_K = \{ 1, 2, 3, 4\}$ with this as given order and $Y$ to
be the set $S_K$ of simplices of $K$, so $S_K$ consists of all non-empty
subsets of $V_K$ that do not contain $\{1, 4\}$.

There are 11 elements in $Y$.

The relation $R$ from $X$ to $Y$ in $x Ry$ if and only if $x$ is a
vertex of simplex $y$.

In $K_R, \{ x_0, \cdots, x_p\}$ is a simplex if there is a $y \in Y$
such that each $x_i \in y$, so with $K_R$ we retrieve exactly $K$
itself.

Before looking at $L_R$ consider a simpler example.

If we consider $\Delta [n]$, the $n$-simplex, with vertices $X = \{0,
1, \cdots, n\}$ and the non-empty subsets of $X$ as simplices then
$K_R$ will be $\Delta [n]$, but the vertices of $L_R$ will be the set
of simplices of $\Delta [n]$, the 1-simplices of $L_R$ will be pairs
of simplices with non-empty intersection. In particular for each
vertex, $i$, of $\Delta [n]$, there will be a $(2^n - 1)$-simplex in
$L_R$ namely that obtained by considering the power set of $X
\backslash \{i\}$ (yielding $2^n$ elements) and adding in the
singleton $\{i\}$ to each of these sets. For instance  for $n = 2$,
$ X = \{0, 1, 2\}$ and there is a 3-simplex $\left\{\{0\}, \{0,
  1\}, \{0, 2 \}, \{0, 1, 2\} \right\}$ in $L_R$. Thus $L_R$ has much
higher dimension than the original $K$.

Among the simplices of $L_R$, however, we have all of those that
are totally ordered in the inclusion ordering and these give a
sub-complex of $L_R$ that is isomorphic to $K'$, the barycentric
subdivision of $K$. This is true in general and in our example of
the two 2-simplices with a shared edge, the complex $L_R$ contains
the barycentric subdivision of $K_R$, but also has some higher
dimensional simplices such as
$$\sigma_{2}^{\max} = \left\{ \{2\}, \{1, 2\}, \{2, 3\}, \{2, 4\},
  \{1, 2, 3\}, \{2, 3, 4\} \right\}.$$
Of course the inclusion  map $K'_R \rightarrow L_R$ is part of that structure
used in the lemma. The map from $L'_R$ to $K_R$ is now relatively easy
to describe. The above 5-simplex $\sigma_{2}^{\max}$ is a simplex
because the element 2 is in all of the parts so $\psi
\sigma_{2}^{\max} = 2$. In general of course there will be a choice of
element, for instance,
$$\sigma_{\{2, 3\}}^{\max} = \left\{ \{2, 3\}, \{1, 2, 3\}, \{2, 3,
  4\} \right\},$$
and is a simplex of $L_R$ because its intersection is non-empty as it
contains both 2 and 3, thus there are two different maps one using
$\psi \sigma_{\{2, 3\}}^{\max} = 2$, the other using 3 as image
point. Of course they are contiguous. The complex $L_R$ seems to
include aspects of both the barycentric subdivision and the dual
complex. The explosion in dimension is, of course, typical here as, for
instance, in the case of an open cover of a topological space the nerve
of the cover yields a simplicial complex whose dimension indicates the
multiple overlaps in the cover and as the cover is varied reflects the
covering dimension of the space but typically the Vietoris complex is
of unbounded dimension.
\end{example}
Returning to global actions and groupoid atlases, combining earlier results, we have that:

\begin{prop}
If $\A$ is a global action or groupoid atlas, pointed at $a_0$,
then
$$\pi_1 (\A, a_0) \cong \pi_1 (N(\A), [a_0])$$
where $[a_0]$ is a connected component of some local groupoid
$(\mathcal{G}_\A)_{\alpha}$ with $a_0 \in
(X_\A)_{\alpha}$.\endbox
\end{prop}
\subsection{$V$ and $N$ on morphisms}
We would clearly expect these constructions, $V$ and $N$ to extend to give us functors from global actions /groupoid atlases to simplicial complexes.  Life is not quite that simple, but almost.  We have to check what they do to the various strengths of morphism.

\textbf{Weak morphisms} 

As weak morphisms are defined in terms of local frames, and simplices in $V(\A)$ are essentially just local frames, it is clear that
\begin{lem}
If $f : \A \rightarrow \B$ is a weak morphism, then $f$ induces a simplicial map $V(f): V(\A) \to V(\B)$, by $V(f)\langle x_0, \ldots, x_n\rangle = \langle f(x_0), \ldots, f(x_n)\rangle$, followed by elimination of repeats.
\end{lem}
The only point for comment is the last phrase: $\langle f(x_0), \ldots, f(x_n)\rangle$ may not actually be a simplex as it may involve repeated elements, but on eliminating these repeats we will get a simplex.  This minor technicality can be avoided using simplicial sets where degenerate simplices are part of the structure, but their use would entail other complications so we merely note that `technicality' here.  It causes no real problem.

The corresponding result for $N$ is much more complicated.  For $\{U_\alpha\}$ in $N(\A)$, we know there is some $x\in U_\alpha$ and hence $\{f(x)\}$ is a local frame in  some $(X_\B)_\beta$, so there is some $\beta \in \Phi_\B$ with $\{f(x)\}$ a $\beta$-frame and so an orbit or connected component $\{U_\beta^\prime\}$ containing it.  We could map $\{U_\alpha\}$ to  $\{U_\beta^\prime\}$, but there is no reason to suppose that $\{U_\beta^\prime\}$ will be the only such possibility, there may be many.  Of course, $\{U_\beta^\prime ~|~ f(x) \in U_\beta^\prime\}$, if finite, will define a simplex of $N(\B)$ and so we might, in that case, attempt to define the corresponding mapping from $N(\A)$ to $N^\prime(\B) = SdN(\B)$, that is the barycentric subdivision of $N(\B)$.  This could work with care, but then the functoriality gets complicated since if we have also $g : \B\to \C$ as well, the composite of our possible $N(f)$ with $N(g)$ is not possible as the former ends at $N^\prime(\B)$ whilst the latter starts at $N(\B)$.  We could apply subdivision to the second map and then compose but then the composite ends up at $N^{\prime\prime}(\C)$ not at $N^\prime(\C)$, which we would need for functoriality. This sort of situation is well understood in homotopy theory as it corresponds to the presence of homotopy coherence caused by the necessity of chosing an image amongst the possible vertices of a simplex.  The different choices are homotopic in a `coherent' way.  It however makes the nerve construction much more complicated to use than the Vietoris one if weak morphisms are being used.
 
One way around the difficulty is to take geometric realisations as  $|N(\B)| \cong |N^\prime(\B)| $, but this defeats the purpose of working with global actions in the first place, which was to avoid topological arguments as they tended to obscure the algebraic and combinatorial processes involved. Probably the safest way is to use $V(\A)$ when developing theoretic arguments involving weak morphisms, but using $N(\A)$ for calculations `up to homotopy' as $N(\A)$ is often much smaller than $V(\A)$.

\medskip

\textbf{Strong morphisms}

If we now turn to strong morphisms, as any strong morphism preserves local frames, it is also a weak morphism and so we have no difficulty in inducing a $V(f)$ from a given strong morphism $f : \A \rightarrow \B$.  We are thus left to see if the `strength' of $f$  allows us to avoid the difficulties we had above in defining a $N(f)$.

Suppose $\{U_\alpha\}$ is a vertex in $N(\A)$ as before, then $U_\alpha$ is a (non-empty) local orbit in $\mathcal{G}_\alpha$.  We have $f = (\eta_X,\eta_\Phi,\eta_\mathcal{G})$, $(\eta_\mathcal{G})_\alpha : \mathcal{G}_{\A,\alpha} \to \mathcal{G}_{\B,\eta_\Phi(\alpha)}$, and $(\eta_\mathcal{G})_\alpha(U_\alpha)$ is contained in a uniquely defined connected component of $\mathcal{G}_{\B,\eta_\Phi(\alpha)}$, which we will denote by $f_*(U_\alpha)$.

Suppose now we have $\sigma = \{U_{\alpha_0}, \ldots, U_{\alpha_n}\}$ in $N(\A)$, so there is some $x\in \bigcap \sigma$.  Consider the family $f_*(\sigma):=\{f_*(U_{\alpha_0}), \ldots, f_*(U_{\alpha_n})\}$, does it have non empty intersection? On objects $(\eta_\mathcal{|G})_{\alpha_i}$ is just $\eta_X$, so $\eta_Xx$ is an object of $f_*(U_{\alpha_i})$ for each $i$ and so $f_*(\sigma)$ is a simplex of $N(\B)$, with the usual proviso that repeats are removed.  The construction of $f_*(\sigma)$ from $\sigma$ has been made without choices and is well defined, moreover if we define $N(f)(\sigma)$ to be this $f_*(\sigma)$, $N(f)$ is a simplicial mapping and there is no problem with functoriality: $N(gf) = N(g)n(f)$.  We thus have that whilst $V$ behaves well with both types of morphsims, $N$ behaves well only with strong morphisms.
\subsection{`Subdivision' of $\A(G, \mathcal{H})$\label{subdivAGH}}

Earlier we saw that in the setting of a group $G$ and a family
$\mathcal{H} = \{H_i \mid i \in \Phi \}$ of subgroups, we could
construct at least two global actions.  In one of these
constructions we took $\Phi_\A = \Phi$ with the discrete order.
Although for this $\A = \A(G,\mathcal{H})$, the description of
$N(\A)$ was very simple (it  is a generalisation of the intersection diagram we used  in section 2).  We noted in the discussion of the
construction of $\Omega \A$ that this type of global action
suffers from the discreteness of its coordinate system as there
were few framings for curves.  In fact the only curves with
framings were those within a single local orbit.

This deficiency serves to highlight  the importance of the order
in $(\Phi_\A, \leq)$ and its influence on the homotopy properties
of $\A$.  If $\A$ does not have enough framings of curves, paths
or loops, then there will be a divergence between the properties
of $\Omega \A$ and those of the loops on $V(\A)$ or $N(\A)$.

Given this difficulty, how can we change the construction of
$\A(G,\mathcal{H})$ to gain more framings of curves?  In fact, if
$\Phi$ is not a singleton, this can be done in a variety of ways,
of graded strength.  We will look in detail at the strongest one.

The problem of framings was that, if $f :\LL \rightarrow \A$ was
a curve, a framing for $f$ was a mapping $\beta :|\LL|
\rightarrow \Phi_\A$, so that $f(m) \in (X_\A)_{\beta(m)}$ for
each $m$ and there was a $b$ in $\Phi_\A$ bigger than both
$\beta(m)$ and $\beta (m + 1)$, so that $f(m)$ and $f(m+1)$ were
linked by  some $g$ in $(\mathcal{G}_\A)_b$.

If $\Phi_\A$ is discrete, then $b \geq \beta(m)$ and $b \geq
\beta(m+1)$ implies equality of $\beta(m)$ and $\beta(m+1)$.  This
is not the intuition intended, but it is not the fault of
`framings', rather of the discreteness of $\Phi_\A$.  Intuitively
we expect $f(m)$ to be in two of the local orbits, so as to link
previous vertices to the new ones later in the sequence.  We thus
need these intersections there.  We could do this by replacing
$\mathcal{H}$ by $\{H_i \cap H_j \mid i,j \in \Phi\}$, but then
higher dimensional homotopy might perhaps suffer.  It is easier to
close $\mathcal{H}$ up under finite intersections as follows:

\medskip

Take $G$, $\mathcal{H}$ as before with $\mathcal{H} = \{H_i \mid i
\in \Phi\}$. Define a (new) global action $\A^\prime(G,
\mathcal{H})$ by
\begin{align*}
X = |X_{\A^\prime}| & = |G|,\text{ the underlying set of } G\\
\Phi_{\A^\prime} & = \textrm{ the set of non-empty subsets of }\Phi \\
&\quad\quad \textrm {ordered by}\supseteq,\textrm{ i.e. } \alpha \leq \beta
\textrm{ if }
\alpha\supseteq \beta;\\
(X_{\A^\prime})_{\alpha} & = X_{\A^\prime} \text{ for all } \alpha
\in
\Phi_{\A^\prime}\\
(G_{\A^\prime})_\alpha &= \bigcap_{i \in \alpha}H_i\text{
operating by left multiplication}
\end{align*}
with :\\
\hspace*{1cm}if $\alpha \leq \beta $, (so $\alpha\supseteq \beta$), then
$$(G_{\A^\prime})_{\alpha\leq \beta } :\bigcap_{i \in \alpha}H_i \rightarrow
\bigcap_{j \in \beta}H_j$$being inclusion.

We can think of $\A^\prime = \A^\prime(G, \mathcal{H})$ as a
`subdivision' of $\A(G, \mathcal{H})$, rather like the barycentric
subdivision above.  The effect
of this `subdivision' on the simplicial complexes $V(\A)$ and $N(\A)$ is interesting:\\
a) Vietoris: $V(\A^\prime) = V(\A)$, since any $\{x_0, \ldots,
x_n\} \in V(\A^\prime)_n$ is there because there is some orbit
$(\bigcap_{i \in
  \alpha}H_i)x_0$ containing it, but then $\{x_0, \ldots, x_n\}\subseteq
H_ix_0$ for any $i\in \alpha$, so there are no simplices in either of the two
complexes, not in the other.\\
b) Nerve:  the relationship is more complex.  We know from
Dowker's theorem (above) that $|N(\A)| \simeq |N(\A^\prime)|$ and,
denoting by $N(\A)^\prime$, the barycentric subdivision of
$N(\A)$, we can relate $N(\A)$ and $N(\A)^\prime$ to
$N(\A^\prime)$. The old local orbits of $\A$ are still there in
$\A^\prime$, so we have an inclusion $$N(\A)\rightarrow
N(\A^\prime)$$ corresponding to a (strong) morphism from $\A$ to
$\A^\prime$.

Now assume we have $\sigma \in N(\A^\prime)$.  The vertices of
$\sigma$ can be totally ordered by inclusion:
$$\sigma = \{\{H_ix_i \mid i \in \alpha_0\}, \ldots, \{H_i x_i\mid  i \in \alpha_n\}\}.$$
with $\alpha_0 \subseteq  \ldots\subseteq  \alpha_n$ and
$\bigcap\{H_i \mid i \in \alpha_n\}\neq \emptyset$.  We therefore
have an element, $a$ in this intersection and so $\sigma$ can also
be written $$\sigma = \{\{H_ia \mid i \in \alpha_0\}, \ldots, \{H_i
a\mid  i \in \alpha_n\}\}.$$We can assign a vertex of $N(\A^\prime)$
to each of the vertices of $\sigma$, by sending $\{H_ia \mid i \in
\alpha\}$ to its intersection $(\bigcap H_i)a$.  This will give a
simplicial map from $N(\A)^\prime$ to $N(\A^\prime)$, but will
often collapse simplices. (It is well behaved as a map of
simplicial \underline{sets}\footnote{If we put a total order on the vertices of a simplicial complex, $K$ it naturally gives us a simplicial set, $K^{\textrm{simp}}$ as we can control the notion of degenerate simplices, i.e. ones with repeated entries.  We will not be using this much in this paper but we note that it can be useful when handling maps that collapse simplices as we have seen earlier.  The corresponding simplicial set $V^\textrm{sing}(\A)$ is related to the bar resolution in the case that $\A = \A(G,\mathcal{H})$.} , but not that nice at the simplicial complex level.)  For instance, if $\mathcal{H} = \{H_1, H_2, H_3\}$
with $H_i \cap H_j = \{1\}$ if $i \neq j$, then the 2-simplex,
$$\{\{H_1\},\{H_1, H_2\},\{H_1, H_2, H_3\}\}$$gets mapped to $\{H_1,\{1\}\}$,
a 1-simplex, and there are a lot of other collapses as well.

\medskip

The global action / groupoid atlas $\A^\prime$ has the necessary
property with regard to framings:

Suppose $f :\LL \rightarrow \A^\prime$ is a curve, then for $m
\in \mathbb{Z}$, we have a $\beta- \in \Phi$ with a $g- :f(m)
\rightarrow f(m+1)$ in $(\cal{G}_\A)_{\beta -}$ and a $\beta+ \in
\Phi$ with a $g+ :f(m) \rightarrow f(m+1)$;  we take $\beta (m)$ to
be $\{ \beta -, \beta +\}$, or any family containing this.  Of
course, when we look at $f(m+1)$, we can take its $\beta- $ to be
the $\beta+ $ of $f(m)$, so $\beta+ \geq b(m)$, $\beta(m +1)$ and
the framing can be constructed. (Of course, the element of choice
here can be avoided by replacing $\beta+ $ by the family of all
$\alpha \in \Phi$ such that a suitable  $g+$ exists in
$(\cal{G}_\A)_\alpha$.)

If instead of $\LL$, we were mapping in higher dimensional
objects, we would need, not just pairwise families, but all, as we
have done. Effectively, in replacing $\A$ by $\A^\prime$, we have
introduced an object that is more `complete' with respect to local
frames than the original $\A$.  This `completeness' allows much better exponentiation properties: we would be
able to form $\A^\B $ for any $\B $ and the result will have the
`right' sort of behaviour.  The `completeness' property required
is called the `infimum condition', (cf. Bak, \cite{bak1,bak2}).

A groupoid atlas $\A$ satisfies the \emph{infimum condition} if
given any non-empty finite subset $U\subseteq |\A|$, the set
$$\{\alpha \in \Phi_\A :U \text{ is a local frame in } \alpha\}$$ is empty or
has an initial element in the order induced from $\Phi_\A$.

\medskip

\begin{example}
$ \A^\prime (G,\mathcal{H})$ is infimum (i.e. satisfies the
infimum condition).

If $U=\{x_0, \ldots, x_n\}$ is any finite set of elements of $G$,
let $$\alpha = \{H \in \Phi :x_1, \ldots, x_n \in Hx_0\}.$$ If
$\alpha$ is non-empty, then $\alpha \in\Phi_{\A^\prime}$ and $x_0,
\ldots, x_n$ is in $(\bigcup_\alpha H)x_0$. This $\alpha$ is thus
the initial element required, or is empty. \endbox
\end{example}
\medskip

\begin{question}  Given any groupoid atlas, $\A$, find a groupoid
atlas $\A^\prime$, which satisfies the infimum condition, comes
together with a strong morphism,$$\A\rightarrow \A^\prime$$ and, if
possible, is universal with these properties.\end{question}

\medskip

\textbf{Remarks on the problem.}

It seems that $ \A^\prime (G,\mathcal{H})$ will be the solution
for $ \A (G,\mathcal{H})$. Presumably other single domain global
actions will be `completed' in a similar way.

\medskip

 We plan to return to the infimum condition later in this sequence of papers.

\medskip

The passage from   $ \A (G,\mathcal{H})$ to $ \A^\prime
(G,\mathcal{H})$ really  corresponds to an operation on
$\mathcal{H}$, giving $\Phi$ itself more structure and closing
$\mathcal{H}$ up under intersections.  If we extend the notation $
\A (G,\mathcal{H})$ to include families $\mathcal{H}$ with
additional structure, then a convenient notation is $ \A
(G,\bar{\mathcal{H}})$ for what we have denoted $ \A^\prime
(G,\mathcal{H})$, thus emphasising that it is the `closure' of
$\mathcal{H}$, that is used.  This construction of $ \A
(G,\bar{\mathcal{H}})$, with two additional conditions, is closely
related to Volodin's definition of $K$-groups of rings with extra
structure, \cite{volodin}.  We will assume $\mathcal{H}$ is
already `closed' in this way:

Suppose as before that $G$ is a group.  Let  $\mathcal{H}$ be a
family of subgroups of $G$ indexed by $(\Phi_\A,\leq)$, where
\begin{enumerate}[(i)]
\item $H_\alpha = H_\beta$ if and only if $\alpha = \beta$; \item
 $\alpha \leq \beta$ if and only if $H_\alpha \subseteq H_\beta$;
\item  there is some $* \in \Phi_\A$ with $H_* = \{1\}$.
\end{enumerate}
We assume $\mathcal{H}$ is closed under arbitrary intersections,
and that if $H_\alpha$ and $H_\beta$ are contained in some
$H_{\gamma^\prime}$, $\gamma^\prime \in \Phi_\A$, then the
subgroup generated by $H_\alpha$ and $H_\beta$, denoted $\langle
H_\alpha, H_\beta\rangle$, is itself some $H_{\gamma}$ for
$\gamma\in \Phi_\A$.  (Thus the order structure of
$(\Phi_\A,\leq)$ is almost a lattice, but a top element need not
exist.)  In this case $ \A (G,\mathcal{H})$  is called a
\emph{Volodin model}.

\section{ Calculations of fundamental groups - some easy examples.}

Our examples will all be single domain global actions, i.e. the
local actions are all based on a single set, $(X_\A)_{\alpha} =
X_\A$ for all $\alpha \in \Phi_\A$. They will all be of the form
$\A = \A (G, \mathcal{H})$, where $G$ is a group and $\mathcal{H}
= \{ H_i \mid i \in \Phi\}$ is a family of subgroups (see section 1).

\begin{example}
(already considered in example \ref{S3})
\begin{align*}
G = S_3 & = \langle a, b ~|~ a^3 = b^2 = (ab)^2 = 1 \rangle, \text{ so
  } a = (1, 2, 3),~ b = (1, 2);\\
H_1 =\langle a\rangle & = \{1,~ (1, 2, 3),~ (1, 3, 2)\},\\
H_2 = \langle b\rangle & = \{1, ~(1, 2) \};\\
\mathcal{H} & = \{H_1, H_2\}
\end{align*}
\end{example}
The intersection diagram given in our earlier look at this example is in fact the nerve
$N(\A)$ having 5 vertices and 6 edges. The other complex $V(\A)$
is almost as simple. It has 6 vertices corresponding to the 6
elements of $S_3$, and each orbit yields a simplex\\ $H_1 = \{1,
a, a^2\}$ gives a 2-simplex (and 3 ~1 -simplices),\\ $H_1b = \{ b,
ab, a^2b\}$ also gives a 2-simplex;\\ $H_2 = \{1, b\}$ yields a
1-simplex, as do its cosets $H_2 a$ and $H_2 a^2$.

We thus have $V(\A)$ has two 2-simplices joined by 1-simplices at
the vertices, (see below).

As $N(\A)$ is a connected graph with 5 vertices and 6 edges, we know $\pi_1 N(\A)$ is free on 2
generators. (The number of generators is the number of edges
outside a maximal tree.) This same rank can be read of equally
easily from $V(\A)$ as that complex is homotopically equivalent to a bouquet
of 2 circles, (i.e. a figure eight).
\begin{figure}[ht]
\begin{center}
\includegraphics[scale = 1.0]{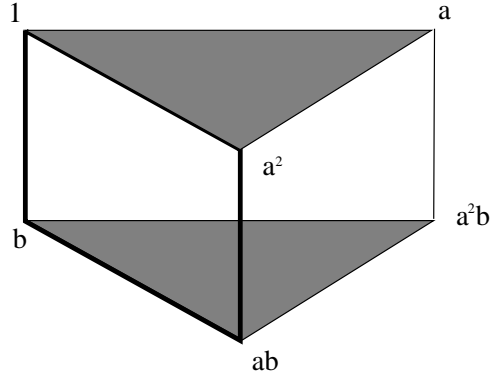}
\caption{$V(\A(S_3,\{\langle a\rangle,\langle b\rangle\}))$}
\end{center}
\end{figure}
The generators can be identified with words in the free product $H_1
\ast H_2$ (one choice being shown in example \ref{S3}) and relate to the
kernel of the natural homomorphism from $H_1 \ast H_2$ to $S_3$.\\

The heavy line in the figure corresponds to a loop at 1 given by
\[
\xymatrix{
1 \ar@{>}[r]^b & b \ar@{>}[r]^a &ab \ar@{>}[r]^b &a^2 \ar@{>}[r]^a & 1
}
\]
and the word is $abab\in C_2 \ast C_3$.

The reason that this happens in clear. Starting at 1, each part of
the loop corresponds to a left multiplication either by an element
of $H_1 \cong C_3$ or of $H_2 \cong C_2$. We thus get a word in
$H_1 \ast H_2 \cong C_2 \ast C_3$. As the loop also finishes at 1,
we must have that the corresponding word must evaluate to 1 when
projected down into $S_3$.

In more complex examples, the interpretation of $\pi_1 (V(\A), 1)$
will be the same, but sometimes when $G$ has more elements,
$N(\A)$ may be easier to analyse than $V(\A)$. The important idea
to retain is that the two complexes give the same information, so
either can be used or both together.

Some of the limitations of the information encoded by $\pi_1 (\A)$
are illustrated by the next two examples.

\begin{example}\label{K4}
$G = K_4$, the Klein 4 group, $\{ 1, a, b, c \} \cong C_2 \times
C_2$,
so $a^2 = b^2 = c^2 = 1$ and $ab = c$;\\
$\mathcal{H} = \{H_a, H_b, H_c\}$ where $H_a = \{1, a\}$, etc. Set
$\A_{K4} = \A(K_4, \mathcal{H})$.

The cosets are $H_a, H_a b, H_b, H_b a, H_c, H_c a$ each with two
elements so \\
\centerline{$V(\A_{K4}) \cong$ the 1-skeleton of $\Delta [3]$}
\[
\xymatrix@C=1pc{
& a \ar@{-}[dr] \ar@{-}[ld] \ar@{-}[rrrd]&&&\\
1  \ar@{-}[rrrr] \ar@{-}[rrrd]&& ~\ar@{-}[rd] && c \ar@{-}[ld]\\
&&& b &
}
\]
$N(\A_{K4})$ is ``prettier'':

Labelling the cosets from 1 to 6 in the order given above, we have 6
vertices, 12 1-simplices and 4 ~2-simplices. For instance $\{1, 3,
5\}$ has the identity in the intersection, $\{1, 4, 6\}$ gives $H_a
\cap H_b a \cap H_c a$, so contains $a$ and so on. The picture is of the
shell of an octahedron with 4 of the faces removed.\\

\begin{figure}[ht]
\begin{center}
\includegraphics[width = 4cm,height=5cm]{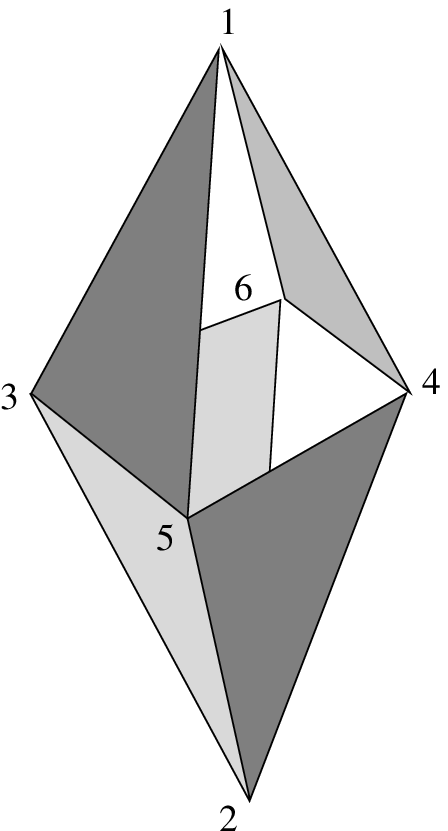}
\caption{$N(\A_{K4})$}
\end{center}
\end{figure}

From either diagram it is clear that $\pi_1 \A_{K4}$ is free of
rank 3.

Again explicit representations for elements are easy to give.

Using $V(\A)$ and the maximal tree given by the edges 1a, 1b and
1c, a typical generating loop would be
$$1 \rightarrow a \rightarrow b \rightarrow 1,$$
i.e., $(1, a, b, 1)$ as the sequence of points. There is a strong
representative for this, namely
\[
\xymatrix{
1 \ar@{>}[r]^a & a \ar@{>}[r]^c &b \ar@{>}[r]^b & 1
}
\]
and up to shifts, this is the only strong representative. \endbox
\end{example}

In general any based path at 1 in an $\A(G, \mathcal{H})$ will
yield a word in $\sqcup\mathcal{H}$, the free product of the
family $\mathcal{H}$. Whether or not that representative is unique
depends on whether or not there are complicated intersections and
``nestings'' between the subgroups in $\mathcal{H}$, since for
instance, if $H_i$ is a subgroup of $H_j$, then if $f(n)
\rightarrow f(n + 1)$ using $g \in H_i$, it could equally well be
taken to be $g \in H_j$. The characteristic of the single domain
global actions of form $\mathcal{\A} (G, \mathcal{H})$ is that
since $X_\A = G$, there is only one possible element linking each
$f(n)$ to the next $f(n + 1)$ namely $f(n + 1) f(n)^{- 1}$. We
thus have a strong link between
$$\Gamma (\A (G,\mathcal{ H})) \text{ and } \underset{\cap}{\sqcup}\mathcal{H},$$
the `amalgamated product' of $\mathcal{H}$ over its intersections, and
an analysis of homotopy classes will prove (later) that
$$\pi_1(\A (G,\mathcal{H}), 1) \cong \text{ Ker} (\underset{\cap}{\sqcup}
\mathcal{H} \rightarrow G),$$
since a based path $(g_1, g_2, \cdots, g_n)$ ends at 1 if and only if
the product $g_1 \cdots g_n = 1$. These identifications will be
investigated more fully (and justified) shortly.

\begin{rem}
Many aspects of these $\A(\mathcal{G, H})$ are considered in the
paper by Abels and Holz \cite{AbelsHolz}. In particular the above
identification of $\pi_1(N(\A (G, \mathcal{H})))$ in terms of the
kernel of the evaluation morphism is Corollary 2.5 part (b) (page
318). Their proof uses covering space techniques. We will explore
other aspects of their paper later.
\end{rem}

\begin{example}
The number of subgroups in $\mathcal{H}$ clearly determines the
dimension of $N(\A)$,  when $\A = \A(G, \mathcal{H})$. Here is
another 3 subgroup example.

Take $q 8 = \{1, i, j, k, - 1, -i,- j, -k\}$ to be the quaternion
group, so $i^4 = j^4 = k^4 = 1$, and $ ij = k$. Set $H_i = \{1, -
1, i, - i\}$ etc., so $H_i \cap H_j = H_i \cap H_k = H_j \cap H_k
= \{1, - 1\}$ and let $\mathcal{H} = \{H_i, H_j, H_k\}$,\\and $
\A_{q 8} = \A(q 8, \mathcal{H})$.

Then $N(\A_{q 8})$ is, as above in  Example \ref{K4}, a shell of an
octahedron with 4 faces missing. Note however that $V(\A_{q 8})$
has 8 vertices and, comparing with $V(\A_{K4})$, each edge of that
diagram has become enlarged to a 3-simplex. It is still feasible
to work with $V(\A_{q8})$ directly, but $N(\A_{q 8})$ gives a
clearer indication that
$$\pi_1 (\A_{q 8}, 1) \text{ is free of rank 3.}$$
\end{example}
\begin{example}
Consider next the symmetric group, $S_3$, given by the presentation
$$S_3 := \langle x_1, x_2 ~|~ x_1^2 = x_2^2 = 1, (x_1x_2)^3 = 1
\rangle$$
Take $H_1 =\langle x_1\rangle$, $H_2 = \langle x_2\rangle$ so both are of
index 3. Each coset intersects two cosets in the other list giving a nerve of
form (see below):
\begin{figure}[h]\begin{center}
\includegraphics[scale = 1.0]{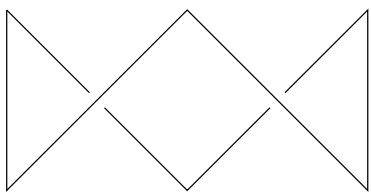}
\caption{$N(\A(S_3, \mathcal{H}))$}
\end{center}\end{figure}
so $\pi_1 (\A(S_3, \mathcal{H}))$ is infinite cyclic.
\end{example}
 \begin{example}\label{S4}
The next symmetric group, $S_4$, has presentation
$$S_4: = \langle x_1, x_2, x_3 ~|~ x_1^2 = x_2^2 = x_3^2 = 1, (x_1
x_2)^3 = (x_2 x_3)^3 = 1, (x_1 x_3)^2 = 1 \rangle.$$ Take $H_1 =
\langle x_1, x_2 \rangle$, $H_2 = \langle x_2, x_3 \rangle$, $H_3
= \langle x_1, x_3 \rangle$. $H_1$ and $H_2$ are copies of $S_3$,
but $H_3$ is isomorphic to the Klein 4 group, $K_4$. Thus there
are 4 + 4 + 6 cosets  in  all. There are 36 pairwise intersections
and each edge is in two 2-simplices. Each vertex is either at the
centre of a hexagon or a square, depending on whether it
corresponds to a coset of $H_1, H_2$ or of $H_3$. There are 24
triangles, and $N(\A(S_4, \mathcal{H}))$ is a surface. Calculation
of the Euler characteristic gives 2, so this is a triangulation of
$S^2$, the two sphere. It is almost certainly  the dual of the `permutahedron'.(Thanks to Chris Wensley for help with the calculation using GAP.)

The fundamental group of $\A(S_4,\mathcal{H})$ is thus trivial and
using the result mentioned above
$$S_4 \cong \underset{\cap}{\sqcup}  H_i,$$
the coproduct of the subgroups amalgamated over the intersection..
 \end{example}

\begin{questions}
\hspace{0.5em}

\begin{enumerate}[1)]
\item Taking $S_n := \langle x_1, \cdots, x_{n - 1} ~|~ x_i^2 =
1, i = 1, \cdots, n- 1, \text{ all }(x_i x_{i + 1})^3 = 1$, and
$(x_i x_j)^2 = 1$
  if $|i - j| \geq 2 \rangle, \mathcal{H} = \{H_i\mid i = 1, \cdots, n -
  1\}$ where $H_i$ is generated by all the $x_j$ \emph{except}
  $x_i$, investigate if $N(\A(S_n, \mathcal{H}))$ is an $(n-2)$-
  sphere. (It is known to be of the homotopy type of a bouquet of
  $S^{n - 2}s$ cf. proof in Abels and Holz's paper \cite{AbelsHolz} for Tits systems,
  $S_n$ being a Coxeter group.)
\item There may be a link between $N(\A( G, \mathcal{H}))$ and
  various lattices of subgroups, classifying spaces for families of
  subgroups (L\"uck et al, Dwyer, ...). This needs investigation.
\item Examine other classes of groups e.g. generalised Coxeter
  groups relative to parabolics (cf. Abels and Holz, \cite{AbelsHolz}), the triangle and
  von Dyck groups, $\Delta (\ell, m, n)$ and $D(\ell, m, n)$ and graph
  products of groups including graph groups.
\end{enumerate}
\end{questions}

\section{Single domain global actions I}

In this section we will continue the study of single domain global
actions including a discussion of their path spaces and fundamental
groups. Certain facets of this study must wait until we have a theory
of covering spaces for global actions/groupoid atlases.

We have so far examined in detail only those single domain global
actions of the form $\A = \A(G, \mathcal{H})$ where $G$  is a
group and $\mathcal{H}$ is a family of subgroups of $G$. There are
two obvious variants:\label{AGKH}

\begin{enumerate}[\hspace{-0.5em} (i)]
\item Usually so far we have let $\mathcal{H} = \{ H_i \mid i \in\Phi\}$ be
  a family of subgroups without any order structure on $\Phi$. More
  generally one can take $\Phi$ to have a reflexive relation on it
  mirroring the intersection and subgroup relations between subgroups
  in $\mathcal{H}$. This seems to make little difference to invariants
  such as $N (\A)$. Of course, this is related to the discussion of
  `subdivision' above.
\item If $K$ is a subgroup of $G$ then we can form a variant of
  $\A(G, \mathcal{H})$ mod $K$. We will write this as $\A = \A ((G, K),
  \mathcal{H})$. It is specified by
\begin{align*}|X_\A|& = G/K, \text{ the  set of right cosets of } K \text{ in } G,\\
(X_\A)_{\alpha} &= X_\A \text{ for all } \alpha \in \Phi\\
\intertext{(where $\Phi$ may be as in (i) above)} H_i
&\curvearrowright X_\A \text{ by left multiplication}
\end{align*}
so the local orbits are of the form $H_i x K$.

We will see later that \emph{all} connected single domain global
actions have this general form, up to isomorphism.
\end{enumerate}

We now turn to the investigation of the fundamental groups of single
domain global actions (of this simple form). First we look at some results on
paths under fixed end point homotopies.
\subsection{The based paths on $\A(G,\mathcal{H})$}
Suppose $\A$ is a general global action with $\Phi_\A$ as
coordinate system, $(X_\A)_{\alpha}$ as local sets, $\alpha \in
\Phi_\A$ and $(G_\A)_{\alpha}$ as local groups. We assume $\A$ is
connected and that a base point $a_0 \in \A$ has been chosen. In
our earlier discussion we saw that a path in $\A$ based at $a_0$
is given by a curve
$$f: L \rightarrow \A$$
and thus by a sequence of points $f(n) \in X_\A$ such that

\begin{enumerate}
\item[(i)] there is a stabilisation pair $(N^-_f, N^+_f)$, i.e., for $n \leq
  N^-_f$, $f(n) = a_0$, and for $n \geq N^+_f$, $f(n) = f(N^+_f)$;
\item[(ii)] there is a sequence of $\beta_n \in \Phi_\A$ with
arrows
  $g_{\beta_n} :f(n) \rightarrow f(n + 1)$ in the groupoid
  $(\mathcal{G}_\A)_{\beta_n}$ (and as $\A$ is a global action, we assume
  $g_{\beta_n}\in (G_\A)_{\beta_n}$), and $g_{\beta_n} f(n) = f(n + 1)$.
\end{enumerate}
(Furthermore the $g_{\beta_n}$ are assumed to stabilise to the
identity arrows for $n \geq n^+_f$, or $n \leqq n^-_f$.)

As we are considering weak curves, the particular $g_{\beta}$ used
are not really in question. If $\A = \A(G, \mathcal{H})$, then
this does not matter as there will be a unique $g_{\beta}$
satisfying $g_{\beta_n} f(n) = f(n + 1)$, namely $g_{\beta_n} =
f(n + 1) f(n)^{- 1}$. In a relative case, the $g_{\beta}$ will be
determined up to multiplication by elements of $K$, and in the
case of interrelationships between the $H_i$, the same element may
be considered as an element of different $H_i$. In complete
generality, we can thus say little about the elements $g_{\beta}$.
Because of this, we will initially assume $\A = \A(G,
\mathcal{H})$, and may therefore base it at 1. We also assume $\A$
is connected.

We can thus consider $f$ to be represented by a word in $\sqcup
H_i$. As $f$ starts at 1, $f$ looks like
\[
\xymatrix{
1 \ar@{>}[r]_{g_1} & f(1) \ar@{>}[r]_{g_2} &f(2) \ar@{>}[r]_{g_3} &
\cdots \ar@{>}[r] &f(N^+_f)
}
\]
(We have reindexed to set $N^-_f$ to 0.) The partial word $(g_m,
\cdots, g_1)$ determines $f(m)$ since $g_m \cdots g_1 = f(m)$. We
can thus think of $f$ as being this list of elements $(g_m, \cdots,
g_1)$ and hence as an element of $\sqcup H_i$. We will need to examine
the loops (in which $f(N^+_f) = 1$), but can examine homotopy of based
paths independently of what value is taken by $f(N^+_f)$.

A fixed end point homotopy from
\begin{align*}f &\leftrightarrow (g_N, \cdots, g_1)\\
\intertext{to }
f'& \leftrightarrow (g'_{N'}, \cdots, g'_1)\\
\intertext{is given by a map of $L \times L$ that stabilises in both
directions. It can therefore be thought of as a sequence of ``moves''
on $f$s each corresponding to an elementary homotopy,}
h :L \times L &\rightarrow \A\\
h (m, 0) &= f_0(m)  = f(m) \\
\intertext{with}
h (m, 1) & = f_1 (m) \\
h (m, n) & = h(m, 0) \quad n \leq 0\\
\intertext{and}
h (m, n) & = h(m, 1) \quad n \geq 1.
\end{align*}
We then can visualise $h$ as being given by a ``ladder''
\[
\CD
\xymatrix{
1 \ar@{=}[d]\ar@{>}[r]^{g_1}& f(1)
\ar@{>}[d]_{h_1}\ar@{>}[r]^{g_2}&f(2) \ar@{>}[r] \ar@{>}[d]_{h_2}&\\
1 \ar@{>}[r]_{g'_1}& f'(1) \ar@{>}[r]_{g'_2} &f'(2) \ar@{>}[r] &
}
\endCD \quad \quad \text{etc.}
\]
where $h_i f(i) = f'(i)$ and the square
\[
\xymatrix{
f(i - 1) \ar@{>}[d]_{h_{i - 1}}\ar@{>}[r]^{g_i} & f(i)
\ar@{>}[d]^{h_i}\\
f'(i - 1) \ar@{>}[r]_{g_i} & f'(i)
}
\]
``commutes'', so for some $H_{\alpha_i}, \alpha_i \in\Phi$,
\begin{enumerate}
\item[a)]$g'_i, g_i, h_{i - 1}, h_i \in H_{\alpha_i}$;
\item[b)]$g'_i h_{i - 1} = h_i g_i.$
\end{enumerate}
This of course means that $h_i \in H_{\alpha_i} \cap H_{\alpha_{i +
    1}}$

\begin{rem}
For convenience we assumed $f(0) = 1$, i.e., $N^-_f$ can be taken to be
0, but the discussion of ripple homotopies earlier shows that shifting
$f$ to left or right keeps within the fixed end point homotopy class
and, of course, does not change the representing word in $\sqcup H_i$.
\end{rem}

Returning to elementary homotopies, we can read off a new based path
from the diagram above namely
\[
\xymatrix{
1 \ar@{>}[r]^{g_1} & f(1) \ar@{>}[r]^{h_1} &f'(1) \ar@{>}[r]^{g'_2} &
f'(2) \ar@{>}[r]^{g'_3} & \cdots
}
\]
and so a representing word $(g'_{N'}, \cdots, g'_2, h_1, g_1)$. This
process can be repeated to get
\[
\xymatrix{
1 \ar@{>}[r]^{g_1} & f(1) \ar@{>}[r]^{g_2} &f(2) \ar@{>}[r]^{h_2} &
f'(2) \ar@{>}[r]^{g'_3} & \cdots
}
\]
so we can track the homotopy from $f$ to $f'$ \emph{within} the
representing words by moves that transfer
$$(\cdots, g'_i h_{i - 1}, g_{i - 1}, \cdots )$$
to
$$(\cdots, g'_i, h_{i - 1}, g_{i - 1}, \cdots )$$
and thus to
$$(\cdots, g'_i, h_{i - 1} g_{i - 1}, \cdots )$$
which of course equals
$$(\cdots, g'_i, g'_{i - 1} h_{i - 2}, \cdots ).$$

In other words homotopy between based paths corresponds exactly to
passing between representing words in $\sqcup H_i$ by the usual moves
that give the amalgamation over the (pairwise) intersections. We thus
have proved
\begin{prop}
If $\A = \A(G,{\mathcal{H}})$, then
$$|\Gamma (\A, 1)|/_{\sim} \cong \underset{\cap}{\sqcup} \mathcal{H},$$
the amalgamated coproduct of the groups in
$\mathcal{H}$.\hfill$\Box$ \end{prop}

\begin{remarks}
\hspace{0.5em}

\begin{enumerate}
\item[(i)] As the usual construction of universal covering spaces
in   topology is the analogue, there, of the left hand side of
this   isomorphism, it is natural to expect the right hand side,
the   amalgamated coproduct, to play that role here. We will look
at   coverings separately, and in some detail, shortly so here it
suffices to note that the end point map
$$\Gamma (\A, 1) \rightarrow \A$$
induces a map
$$|\Gamma (\A, 1)|/_{\sim} \rightarrow G$$
which interprets as the natural evaluation of a word $(g_m, \cdots,
g_1)$ to the product, $g_m \cdots g_1$, i.e., to the natural homomorphism
$$\underset{\cap}{\sqcup} \mathcal{H} \rightarrow G,$$
induced by the universal colimit-property of the amalgamated
coproduct and the inclusions of the subgroups $H_i$ into $G$.
\item[(ii)] We note for future examination that $ \Gamma (\A, 1)$
\emph{has a global action/groupoid atlas structure} and it is
natural to expect that the quotient by fixed end point homotopies
will inherit a similar structure, but that we have not yet
described the construction of colimits, and in particular,
quotients, in this setting. In the particular case above $\A = \A
(G, \mathcal{H})$, it is easily seen that the amalgamated
coproduct carries a global action structure:

There are inclusions
$$i_{\alpha} :H_{\alpha} \rightarrow \underset{\cap}{\sqcup}
\mathcal{H}$$
and writing
$$\tilde{ G} = \underset{\cap}{\sqcup} \mathcal{H},\hspace{1cm}
\tilde{ \mathcal{H}} = \{i_{\alpha} (H_{\alpha}) :\alpha
\in\Phi\},$$ we can construct, $\tilde{ \A} = \A (\tilde{G},
\tilde{ \mathcal{H}})$.

The map
$$\tilde \A \rightarrow \A$$
is a regular morphism of global actions. (Left as an exercise!)

\end{enumerate}
\end{remarks}

Given the analogy between the above and the topological case, it
is no surprise that restricting attention to the loops at 1 in
$\A$, the defining equation
$$\pi_1 (\A, 1) = \pi_0 (\Omega \A) = |\Omega (\A, 1)|/_{\sim}$$
gives:

\begin{cor}
If $\A = \A(G, \mathcal{H})$,\\
$$\pi_1 (\A, 1) \cong \text{ Ker}( \underset{\cap}{\sqcup} \mathcal{H}
\rightarrow G).$$ \hfill $\Box$\end{cor}

 This result in this context was found by A. Bak. Given the
 identification of $\pi_1(\A, 1)$ with $\pi_1 (V(\A), 1)$ and the Dowker
 theorem identifying this with $\pi_1 (N(\A), 1)$, it can be seen to be
 a version of a result of Abels and Holz, \cite{AbelsHolz}. They, in turn,
 relate it to earlier results of Behr, \cite{Behr}, and Soul\'e, \cite{Soule},
 and mention applications of a related result given by Tits, \cite{Tits}. The
 proof given above has the advantage of being very elementary and ``constructive''!

\begin{question}
Compare the use of $N(\A)$ as a simplicial complex
with
  $N^{\text{simp}}(\A)$, as simplicial set. The action of $G$ (which we
  will look at next) gives $N(\A) /G$ is a simplex, but Abels and Holz
  identify  $\pi_1 (N^{\text{simp}} (\A)/G)$ as being
    $\underset{\cap}{\sqcup} \mathcal{H}$. The comparison of $V(\A)$
     with the bar resolutions of $H_i$ and $G$ studied in Abels and Holz, \cite{AbelsHolz},
     also needs examining in detail (and greater generality).
\end{question}

\medskip

\subsection{Group actions on $N(\A)$}

Further information on $N(\A)$ and $V(\A)$ may be obtained by
exploiting the natural action of $G$ on these simplicial
complexes. This leads to a connection of these single domain
global actions not only with the work of Abels and Holz, but with
related work on complexes of groups by  Corson, Haefliger and
others, \cite{brid&haef,corson1,corson2,corson3,haef1,haef2}.

Again $G$ will be a group, $\mathcal{H} = \{H_i  \mid i \in \Phi\}$ a
family of subgroups and $\A = \A (G, \mathcal{H})$ the
corresponding single domain global action. We will assume that
$\A$ is connected so $G$ is generated by the union of the $H_i$s.
Recall that $N(\A)$ is the simplicial complex given by the nerve
of the covering, $\mathfrak{H}$, of $G$  by left cosets of the
$H_i$, $H_i \in \mathcal{H}$.

The group $G$ acts on $N(\A)$ by right translation. A typical
$n$-simplex of $N(\A)$ is of the form
$$\sigma = \{ H_{\alpha_0} x_0, \cdots, H_{\alpha_n} x_n\}$$
where
$$\bigcap \sigma = \cap_{i = 0}^n ~H_{\alpha_i} x_i \neq
\emptyset.$$
If $g \in G$, we can consider $\sigma. g = \{H_{\alpha_0} x_0g,
\cdots, H_{\alpha_n} x_ng\}$.

If $y \in \cap_{i = 0}^n H_{\alpha_i} x_i = \bigcap \sigma$, then $y
.g \in \bigcap \sigma. g$ so
$$\sigma.g \in N(\A).$$
This is clearly a group action. It is ``without inversion''
(Haefliger) or ``regular'' (Abels and Holz,\cite{AbelsHolz}) in as
much as if $\sigma.g = \sigma$, then $H_{\alpha_i} x_i g =
H_{\alpha_i} x_i$, since the $x_i$ used in a given $\sigma$ are
all distinct.  This implies that the orbit space of $N(\A)$ is
also a complex.

\begin{prop}
If $\sigma = \{H_{\alpha_0} x_0, \cdots, H_{\alpha_n} x_n\}$ is an
$n$-simplex of $N(\A)$ then for any $a \in \bigcap \sigma$
$$\sigma a^{- 1} = \{ H_{\alpha_0}, \cdots, H_{\alpha_n}\}$$
Moreover any finite subset $J$ of $\Phi$ corresponds to a unique
$G$-orbit of $N(\A)$ and vice versa.
\end{prop}
\begin{proof}
As $a \in \bigcap \sigma$, there are elements $h_{\alpha_i} \in
H_{\alpha_i}$ with $a = h_{\alpha_i} x_i$ for $i = 0, 1, \cdots,
n$. Thus
$$H_{\alpha_i} x_i a^{- 1} = H_{\alpha_i}$$
and $\sigma a^{- 1} = \{ H_{\alpha_0}, \cdots, H_{\alpha_n}\}$.
The orbit of $\sigma$ is thus determined by the indices of the
subgroups, $H_i \in \mathcal{H}$, used in it. The orbit of
$\sigma$ then corresponds to the finite subset $J_{\sigma} =
\{\alpha_0, \cdots, \alpha_n\}$ of $\Phi$ and conversely.
\end{proof}

\begin{cor}
For $\sigma = \{H_{\alpha_j} x_j :j = 0, \cdots, n\}$ as above, and
$a \in \bigcap \sigma$, $$Stab_G(\sigma) = a^{- 1} \left( \bigcap \{H_j :j \in
  J_{\sigma}\} \right) a.$$ \end{cor}
\begin{proof}
Write $\sigma_0 = \{H_j :j \in J_\sigma\} \in N(\A)$, then
$\sigma_0 a = \sigma$ so $g \in Stab_G \sigma$ if and only if
$\sigma_0 a g = \sigma_0 a$ i.e. if $\sigma_0 aga^{- 1} =
\sigma_0$ which just says $aga^{- 1} \in Stab_G \sigma_0$. However
$Stab_G \sigma_0$ is clearly equal to $\bigcap \sigma_0$, which
completes the proof.\end{proof}

 \begin{cor}
The space of orbits $N(\A)/G$ is a simplex of dimension
card$(\Phi) - 1$.\hfill $\Box$
\end{cor}

\begin{examples}
\hspace{0.5em}

\begin{enumerate}[\hspace{-0.5em} 1)]
\item $G = S_3, H_1 = \{ 1, (1 ~2 ~3), (1 ~ 3 ~2)\}, H_2 = \{1,
  (1~2)\}$. The nerve $N(\A)$ in this case is the graph given in example \ref{S3} with vertices
\[
\xymatrix@R=1pc{
& H_1 & & H_1 b& \\
H_2 & & H_2 a & &H_2 a^2
}
\]
(where, as there, $a = (1~2~3)$, $b = (1~2))$. The action is given by: $a$
fixes $H_1$ and $H_1b$ and permutes the cosets of $H_2$ in the obvious
way; $b$ permutes $H_1$ and $H_1 b$ \emph{and} $H_2 a$ and $H_2
a^2$, but fixes $H_2$ (of course). On 1-simplices
$$a \in H_1 \cap H_2 a \text{~~ so ~~} H_1 a^{- 1} \cap H_2 aa^{- 1} = H_1
\cap H_2 \neq \emptyset$$ and so on. It is thus easy to see that
$N(\A)/S_3 \cong \Delta [1]$.

Of more interest are the examples: \item $G = K_4 = \{1, a , b,
c\}$, $N(\A_{K4})$ is the octahedral shell
  with 4 faces removed. Using the same notation as before: $a$ fixes 1
  and 2, permutes 3 and 4, and also 5 and 6, so in the diagram in example \ref{K4}, $a$
  corresponds to a rotation through $180^{\circ}$ about the vertical
  axis. Similarly for $b$ and $c$, but about the two horizontal axes.
 The orbit space is $\Delta [2]$ as this example has 3 subgroups.
\item $G = q8$, $N(\A_{q8})$ has the action of $q8$ via the
quotient
  homomorphism to $K_4$ and the action outlined before in 2. Of course,
  $N(\A_{q8})/q8$ is again a  2-simplex.

Our final two examples are

\item $S_3$ with $H_1 = \langle (1~2)\rangle$, $H_2 = \langle
  (2~3)\rangle$, so $N(\A)$ is a hexagon (empty) and, of course, the
  $S_3$-action collapses this down to a 1-simplex.\\
\hspace*{-1cm}  and
\item $S_4$ with three subgroups,

$H_1 = \langle (1, 2), (2, 3)
  \rangle$,

 $H_2 = \langle (2, 3), (3, 4)\rangle$\\ and

$H_3 = \langle (1,2), (3, 4)\rangle$.\\ The nerve was found earlier to be a
  triangulation of $S^2$. The action can be specified, but will not be
  given here and, of course, the quotient is $\Delta [2]$.
\end{enumerate}
\end{examples}

This situation is a simple form of a general one considered by
Haefliger (cf. \cite{brid&haef,haef1,haef2}) and Corson (cf. \cite{corson1,corson2,corson3}).
They consider a simplicial complex (or
more generally a simplicial cell complex, cf. Haefliger, \cite{haef1} or a
scwol (small category without loops) cf. Bridson and Haefliger, \cite{brid&haef} ) on which a
group  $G$ acts without inversion. Then $\tilde X/G$ is also a
simplicial (cell)  complex. Their work uses complexes of groups, a notion
generalising that of graphs of groups as in Bass-Serre theory. We will give
definitions shortly, but first need to introduce some more detailed
notation and terminology relating to barycentric subdivisions.

If $K$ is a simplicial complex, we can encode the information in $K$
in a simply way by considering $K$ as a partially ordered set. The
elements of this partially ordered set are the elements of $S_K$, the
set of simplices of $K$ ordered by inclusion. The barycentric
subdivision of $K$ is then just the (categorical) nerve of the poset
$(S_K, \subseteq)$ as noted earlier. We will follow Haefliger \cite{haef1} in
orienting the edges of $K'$ in the following way:

The vertices of $K' (= Sd (K))$ are the simplices of $K$. An (unoriented) edge
of $K'$ consists of a pair $(\sigma, \tau)$ with either $\sigma
\subset \tau$ or $\tau \subset \sigma$. If $a$ is an edge of $K'$
contained in a simplex $\sigma$ of $K$, then the \emph{initial point}
$i(a)$ of $a$ is the barycentre of $\sigma$ (i.e. $\sigma$ as a
vertex of $K'$) and its \emph{terminal point}, $t(a)$, is the barycentre
of some smaller simplex, $\tau$. We write $i(a) = \sigma$, $t(a) =
\tau$ and so have $a = (\tau, \sigma)$, with $\tau \subset
\sigma$. (This is perhaps the opposite order from that which seems
natural, but it avoids considering dual posets later.)

\begin{example}
For the 2-simplex, considered as the simplicial complex of non-empty
subsets of $\{1, 2, 3\}$, this gives
\[
\xymatrix{
& & 2 & & \\
~\\
& 12 \ar@{>}[ruu] \ar@{>}[ldd]& & 23\ar@{>}[luu] \ar@{>}[rdd]\\
& & 123 \ar@{>}[lu] \ar@{>}[ru] \ar@{>}[lld] \ar@{>}[rrd]\ar@{>}[uuu]\ar@{>}[d]& & \\
1& & {13} \ar@{>}[ll] \ar@{>}[rr]& & ~3
}
\]
Although it is usual to consider partially ordered sets as
categories, because his complexes \underline{are} more general
than mere simplicial complexes, Haefliger introduces a specific
construction of a small category associated to $K$ (cf.
\cite{haef1}).

Define a category $C(K)$ with set of objects $S_K$, the set of
vertices of the barycentric subdivision $K'$ of $K$ and with
arrows $\Arr(C(K)) = E_k \sqcup S_k$, the set of edges of $K'$
together with $S_K$. (Of course, the vertices are considered as
identity arrows at themselves.) Two edges $a$ and $b$ are
considered composable if $i(a) = t(b)$ and the composite is $c =
ba$ such that $a, b, c$ form the boundary of a 2-simplex in $K'$.
\[
\xymatrix{
& ~\\
~ \ar@{>}[ru]^c \ar@{>}[r]_b & \cdot \ar@{>}[u]_{a}
}
\]

\end{example}

This category $C(K)$ is an example of a \emph{small category without loops}
as introduced by Haefliger \cite{brid&haef,haef1}. We shall consider a small
category, $\chi$, to consist of a set, $V(\chi)$, of vertices or objects (denoted
here by Greek letters, $\tau$, $\sigma$, etc.) and a set $E(\chi)$ of edges
(denoted by Latin letters, $a$, $b$, \ldots ), together with maps \\
\hspace*{2cm}$i :E(\chi) \rightarrow V(\chi)$, \quad the initial vertex or source map,\\
\hspace*{2cm}$t :E(\chi) \rightarrow V(\chi)$, \quad the terminal vertex or target map\\
and a composition \\
\hspace*{2cm}$E^{(2)}(\chi) \rightarrow E(\chi)$,\\
where $E^{(2)}(\chi) = \{(a,b) \in E(\chi) \times E(\chi) :i(a) =
t(b)\}$, together with  associativity of composition and the rules $i(ba) = i(b)$, $t(ba) =
t(a)$ for $ba$, the composite of $a$ and $b$.

The small category $\chi$ is a small category without loops, or scwol, if for
all $a$ in $E(\chi)$, $i(a) \neq t(b)$.
\begin{rem} Haefliger's definition of a small category without
loops in \cite{brid&haef} (p.521) is optimised for the statement
of the no loops condition, but omits to define composition of an
arbitrary arrow with a vertex.  This is handled correctly (p.573)
in an appendix.  This does not influence the later development.
\end{rem}

For the moment we will move attention back to $K$ and the definition of a
complex of groups.

\subsection{Complexes of groups}

A complex of groups $G(K)$ on $K$ is $(K, G_0, \psi_a, g_{a, b})$
given by
\begin{enumerate}[\hspace{-0.5em} 1)]
\item a group $G_{\sigma}$ for each simplex $\sigma$ of $K$; \item
an injective homomorphism
$$\psi_a :G_{i (a)} \rightarrow G_{t(a)}$$
for each edge $a \in E_K$ of the barycentric subdivision of $K$;
\item for two composable edges $a$ and $b$ in $E_K$, an element
  $g_{a, b} \in G_{t(a)}$ is given such that
$$g^{- 1}_{a, b} \psi_{ba} (_-) g_{a, b} = \psi_a \psi_b$$
and such that the ``cocycle condition''
$$g_{a, cb} \psi_a (g_{b, c}) = g_{ab, c} g_{a, b}$$
holds.

(If the dimension of $K$ is less than 3, this condition is void.)
\end{enumerate}

\textbf{Almost generic example}: developable complexes of groups.

Suppose we have a simplicial  complex $\tilde X$ with a right $G$
action which is ``without inversion'', i.e., if $\sigma.g = \sigma$ then
$xg = x$ for all vertices $x$ of $\sigma$. Write $X = \tilde X/G$ for
the quotient complex. We will specify a complex of groups $G(X)$ on
$X$:\\
Set $p :\tilde X \rightarrow X$ to be the quotient mapping.

For a simplex $\sigma$ of $X$, pick a $\tilde{\sigma} \in X$ with
$p(\tilde \sigma) = \sigma$, we say $\tilde{\sigma}$ is the \emph{chosen lift} of $\sigma$, and set\\
\begin{align*}
G_{\sigma} = & ~G_{\tilde \sigma},\text{ the stability subgroup of }\tilde \sigma,\\
= &~ \{g :\tilde \sigma g = \tilde \sigma\}.
 \end{align*}
For each $a \in E_X$ with $i(a) = \sigma$, let $\tilde a$ be the edge
in $\tilde \sigma$ whose projection is $a$, i.e. $p (\tilde a) = a$ and
$i(\tilde a) = \tilde \sigma$. Then there is some $h_a \in G$ with
$t(\tilde a. h_a) = \tilde \tau$ where $\tilde \tau$ is the chosen lift
of $\tau = t(a)$. (If $t(\tilde a) = \tilde \tau$ already, we agree to
take $H_a$ to be the identity of $G$.)

Define
$$\psi_a :G_{i(a)} \rightarrow G_{t(a)}$$
by
$$\psi_a(g) = h_a^{- 1} g h_a \quad \text{ for } g \in G_{i(a)}.$$
Given two composable edges $a$ and $b$ define
$$g_{a, b} = h_{ba}^{- 1} h_bh_a.$$

{\it  Verification of conditions}

(Although easy to do, this helps the intuition:)

\begin{enumerate}[(i)]
\item Suppose $g \in G_{i (a)}$, then $\tilde a = (\tilde \tau
  h_a^{- 1}, \tilde \sigma)$ or $\tilde a. h_a = (\tilde \tau, \tilde
  \sigma. h_a)$. As $\tilde \sigma g = \tilde \sigma$, and $\tilde \tau h_a^{- 1}
\subset \tilde \sigma$, we have
$$\tilde \tau h_a^{- 1} g = \tilde \tau h_a^{- 1}$$
and $h_a^{- 1} g h_a \in G_{t(a)}$, i.e. $\psi_a (g) \in
G_{t(a)}$. \item Suppose $a, b$ are composable: $i (ba) = i(b)$,

$t(ba) = t(a)$,
  then
$$\psi_b : G_{i (b)} \rightarrow G_{t(b)} = G_{i (a)}$$
$$\psi_b (g) = h_b^{- 1} g h_b.$$
Similarly $\psi_a \psi_b (g) = h_a^{-1} h_b^{- 1}gh_b h_a$,\\
whilst
$$\psi_{ba} (g) = h_{ba}^{- 1} g h_{ba}.$$
It is clear that $g_{a, b}$ as defined above does the job. \item
cocycle condition:
$$g_{a, cb} \psi_a (g_{b, c}) = h_{cba}^{- 1} h_{cb} h_a. h_a^{- 1}
h_{cb}^{- 1} h_c h_b h_a$$
$$g_{a, cb} . g_{a, b} = h_{cba}^{- 1} h_c h_{ba}. h_{ba}^{- 1} h_b
h_a,$$
so it does check out correctly.
\end{enumerate}


In the case of a single domain global action $\A = \A(G,
\mathcal{H})$ where $\mathcal{H} = \{H_1, \cdots, H_n\}$ with $H_i
< G$, then $N(\A)/G \cong \Delta^{n - 1}$. Suppose $\sigma \in
S_{\Delta^{n - 2}}$ then if $\sigma = \{\alpha_1, \cdots,
\alpha_r\}$, we can always choose $\tilde \sigma = \{H_{\alpha_1},
\cdots, H_{\alpha_r} \}$. If $a$ is an edge of $Sd (\Delta^{n -
1})$ then for $i(a) = \sigma$ and $t(a) = \tau, \tilde \tau
\subset \tilde \sigma$ and hence
\begin{align*}
G_{\tau} = G_{\tilde \tau} & = \bigcap \{H_i  \mid i \in \tilde \tau\},\\
G_{\sigma} = G_{\tilde \sigma} & = \bigcap \{ H_i  \mid i \in \tilde
\sigma\},
\end{align*}
so there is no need to have $h_a \neq 1$. Because of this, $\psi_a$ is
simply an inclusion of a subgroup and $g_{a, b}$ can be chosen to be
1. Thus single domain global actions yield simplices of groups of a
particularly simple kind. This does \emph{not} imply that the
more general case is irrelevant to global actions, merely that single
domain global action are ``untwisted''.

\medskip

\begin{question}
Are there `twisted' variants that also arise from
global actions?
\end{question}

Given a complex of groups, both Corson and Haefliger show how to
construct a universal covering complex and a fundamental group which
yields the given complex of groups, provided certain fairly mild
restrictions are satisfied.

\subsection{Fundamental group(oid) of a complex of groups.}

Let $G(K) = (K, G_\sigma, \psi_a, g_{a, b})$ as before.\\
Let $E^{\pm}_K$ be the set of edges of $K'$ with an orientation
\begin{align*}
a^+ = a, \quad &  a^- =   ~a \text{ with the opposite orientation}\\
 \text{so } i(a^-) = & t(a^+) \text{ etc.}
\end{align*}
First define $FG(K)$ to be the group generated by
$$\bigsqcup \{G_\sigma :\sigma \in V_K\} \cup  E^{\pm}_K$$
subject to the relations
\begin{enumerate}
\item[-] the relations of each $G_\sigma ,$
\item[-] $(a^{+})^{- 1} = a^- $ and $(a^-)^{- 1} = a^+$,
\item[-] $\psi_a (g) = a^- g a^+$ for $g \in G_{i(a)}$,
\item[-] $(ba)^+ g_{ab} = b^+ a^+$ for composable $a, b$.
\end{enumerate}
The image of $G_{\sigma}$ in $FG(K)$ will be denoted $\bar
G_{\sigma}$.

Haefliger defines $\pi_1 (G(K), \sigma_0)$ in two equivalent ways:
\begin{Def}

{\it  Version 1.}
If $\sigma_0, \sigma_1 \in V_K$, the vertices of $K$, a $G(K)$-path $c$
from $\sigma_0$ to $\sigma_1$ is a sequence $(g_0, e_1, g_1, \cdots,
e_n, g_n)$, where $(e_1, \cdots, e_n)$ is an edge path in $K'$ from
$i(e_1) = \sigma_0$ to $t(e_n) = e_1$ and $e_i \in E^{\pm}_K$, $ i = 1,
\cdots, n$ and where $g_K \in G_{t(e_k)} = G_{i (e_{k + 1})}$.

Such a $G(K)$-path, $c$, represents $g_0 e_1 \cdots e_ng_n \in
FG(K)$. Two such paths from $\sigma_0$ to $\sigma_1$ are said to
be homotopic if they represent the same element of $FG(K)$. We set
$\pi_1(G(K), \sigma_0, \sigma_1)$ equal to the subset of $FG(K)$
represented by $G(K)$-paths from $\sigma_0$ to $\sigma_1$. When
$\sigma_0 = \sigma_1$, we write
$$\pi_1 (G(K), \sigma_0) = \pi_1 (G(K), \sigma_0, \sigma_0).$$
This is a subgroup of $F(G)$ and is called the {\it fundamental
group} of $G(K)$.

{\it  Version 2.}
Assume $K$ is connected and pick a maximal tree $T$ in the
1-skeleton of $Sd (K) = K'$. Let $N(T)$ be the normal subgroup of
$FG(K)$ generated by $\{a^+ :a \in T\}$, then
$$\pi_1 (G(K), T) \cong FG(K)/ N(T),$$
and hence has a presentation:

\begin{tabular}{lll}
- generators & $\sqcup G_{\sigma} \sqcup E_K$\\
- relations :
& - $g_1 \cdot g_2 = g_1 g_2$ & within any particular $G_{\sigma}$\\
& - $\psi_{\alpha} (g) = \alpha^{- 1} g \alpha$ &  $g \in
  G_{i(\alpha)}$,\\
& - $(\beta \alpha) g_{\alpha, \beta} = \beta .\alpha$ & if
  $\alpha, \beta \in E_K$ are composable\\
& - $\alpha  = 1$ & if $\alpha \in T$.
\end{tabular}
\end{Def}
\begin{example}

Suppose $\A = \A(G, \mathcal{H})$, $\tilde{K} = N(\A)$,
$\mathcal{H} = \{H_1, \cdots, H_n\}$ so $K = \Delta^{n - 1}$.
Pick the maximal tree with edges radiating out from the vertex
$\{H_1\}$, e.g. if $n = 3$, we get figure \ref{Del2}
\begin{figure}[h]
\begin{center}
\includegraphics[scale = 1.5]{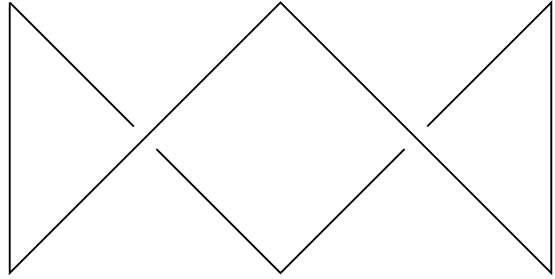}
\caption{Barycentric Subdivision of $\Delta^ 2$ with the chosen maximal tree shown.\label{Del2}}
\end{center}\end{figure}

There is an obvious collapse of $\Delta^{n - 1}$ to $T$. We have
already noted that all the $g_{a, b}$ are trivial in these examples so
we can prove (inductively via the collapsing order) that if $a$ is any
edge in $Sd (\Delta^{n - 1})$, the fact that $\alpha = 1$ for $\alpha
\in T$ implies that $a = 1$ in $\pi_1 (G(K), 1)$. Thus $\pi_1 (G(K),
1)$ has a presentation with

\begin{tabular}{lll}
- generators & $\sqcup G_{\sigma}$&\\
- relations :
& - $g_1 \cdot g_2 = g_1 g_2$ &within any particular $G_{\sigma}$\\
& - $\psi_{\alpha}(g) = g $ &for $g \in G_{i(\alpha)}$\\
\end{tabular}

As $G_{\sigma} = \bigcap \{H_i \mid i \in \sigma\}$, we have
$$\pi_1 (G(K), 1) \cong \underset{\cap}{\sqcup} H_i,$$
the coproduct of the $H_i$ amalgamated over the intersection.

It is noticeable that there is, as before, a homomorphism
$$\pi_1 (G(K), 1) \rightarrow G$$
with kernel $\pi_1 N(\A)$. \endbox
\end{example} 

Thus for single domain global actions, the fundamental group is the same as the fundamental group of the corresponding complex of groups.
\section{Coverings of global actions.}

To complete the study of the various interpretations of $\pi_1
(\A, a_0)$ we need to consider covering maps of global actions and
the analogues of the universal covering.
\subsection{Covering maps}
\begin{Def}
Let $p :\B  \rightarrow \A$ be a morphism of global actions. We
say that $p$ is a {\it covering map} if it satisfies the \emph{unique local frame lifting property}.
Explicitly, for any given local frame $x_0,\ldots,x_n$ in $\A$ and any $y_0\in p^{-1}(x_0)$ in $\B$, there exists a 
unique local frame $y_0,\ldots,y_n$ in $\B$ such that $p(y_i)=x_i$ for all $i$.
\end{Def}

In the context of groupoids, the definition of covering is based on the Stars of the objects (cf. Brown \cite{rb}).
In the global action case, we have notions of local and global Stars:

Let $\C$ be an arbitrary global action. Let $\alpha \in \Phi_\C$
and $x \in X_{\alpha}$, then Star$_\alpha (x) = \{\sigma x~|~
\sigma \in G_{\alpha}\}$. If $x \in X$, then set Star$_C (x) =
\cup \{\text{Star}_{\alpha} (x) ~|~ \alpha \in \Phi_\C \text{ and
} x \in X_{\alpha}\}$. 

\begin{prop}
Let $p :\B  \rightarrow \A$ be a covering of global actions. Then the following two conditions are
satisfied.

\begin{enumerate}
\item[1)] if $x_1, x_2 \in X_\B $ are such that $p (x_1) =
p(x_2)$, $x_1
  \neq x_2$, and there is some $\beta \in \Phi_\B $ with $x_1, x_2 \in
  X_{\beta}$ then Star$_{\beta} (x_1) \cap \text{ Star}_{\beta} (x_2)
  = \emptyset$;
\item[2)] if $x \in X_\B $ then $p
\downharpoonright_{\text{Star}_\B  (x)} :\text{Star}_\B
  (x) \rightarrow \text{Star}_\A (p(x))$ is a bijection.
\end{enumerate}
(In fact 2) implies 1))
\end{prop}

The proof is omitted. 

Note that the converse of this result is false. 
Consider for example the global action $\B$ with 3 points, say $\{a,b,c\}$ and 3 local
sets $\{a,b\},\{a,c\},\{b,c\}$. The cyclic group $C_2$ acts on each local set in the obvious way. Let 
$\A$ be a global action with the same underlying set but with only one local set on which $C_3$ acts. Then the
identity map between the underlying sets induces a morphism $\B\to \A$ which satisfies the conditions of the last
proposition but this map is not a covering map of global actions.

\smallskip

We say $p$ is {\it regular} if $\Phi_\B  = \Phi_\A$ and $p$ is a
regular morphism with $p = (p_{\Phi}, p_G, p_X)$ such that
$p_{\Phi} :\Phi_\B  \rightarrow \Phi_\A$ is the identity map.
\begin{prop}
Any covering is weakly isomorphic to a regular covering. \end{prop}

The proof is omitted, as we will nor be using the result here.

The construction and theory of coverings is very similar to that
in the classical topological theory. We will assume $\A$ is
connected and when necessary that a base point $a_0 \in X_\A$ has
been chosen.

\subsection{Path  lifting and homotopy lifting}

Suppose $p :\B  \rightarrow \A$ is a covering map of global
actions, suppose $f :\LL \rightarrow \A$ is a path and $(N^-_f,
N^+_f)$ is a stabilisation pair for $f$. Set $x_0 = f(N^-_f)$.

\begin{lem}
Given any $\tilde x_0 \in p^{- 1} (x_0)$, there is a unique path
$\tilde f :\LL \rightarrow \B $ with the same stabilisation pair
$(N^-_f, N^+_f)$ such that $p \tilde f = f$ and $\tilde f (N^-_f)
= \tilde x_0$.\end{lem}
\begin{proof}
(For simplicity of notation we assume $N^-_f = 0$, and write
$N^+_f = N$.) Suppose $f = (f(0), f(1), \cdots, f(N))$. There is
an $\alpha_1 \in \Phi_\A$ with $f(0), f(1) \in (X_\A)_{\alpha_1}$,
and a $g_1 \in (G_\A)_{\alpha_1}$ such that $g_1 f(0) = f(1)$.
Thus $x_1 = f(1) \in \text{Star}_\A (x_0)$.  As $
\downharpoonright_{\text{Star}_\B (\tilde x_0)}$ is a bijection
between Star$_\B (\tilde x_0)$ and Star$_\A(x_0)$, there is a {\it
unique} $\tilde x_1$ with $p(\tilde x_1) = x_1$ and a $\tilde
g_{\alpha_1}$ with $\tilde x_1 = \tilde g_{\alpha_1} \tilde x_0$.
A simple use of induction completes the proof. \end{proof}

Again suppose $f_0, f_1 :L \rightarrow \A$ are paths, but in
addition suppose $h :\LL \times \LL \rightarrow \A$ is a (fixed
end point) homotopy from $f_0$ to $f_1$. (This, of course, implies
that $f_0, f_1$ stabilise to points $x_0$ and $x_N$ say.) Again
let $\tilde x_0 \in p^{- 1} (x_0)$ and then we have

\begin{lem}
The two lifts $\tilde f_0, \tilde f_1$ of $f_0, f_1$ are homotopic
by a homotopy $\tilde h$ lifting $h$. In particular $\tilde f_0
(N) = f_1 (N)$.\end{lem}
\begin{proof}
We can assume $h$ stabilises outside a square. We may assume
this square is actually a 1 by 1 square as the general case
follows by induction. Initially we will need
\[
\xymatrix{
x_0 \ar@{>}[r] ^{g'_1} & x'_1\\
x_0 \ar@{>}[r]_{g_1} \ar@{=}[u] & x_1 \ar@{>}[u]_{h_1}
}
\]
but in general,
\[
\xymatrix{
x'_i \ar@{>}[r] & x'_{i + 1}\\
x_i \ar@{>}[r] \ar@{>}[u] & x_{i + 1} \ar@{>}[u]
}
\]

Since $h$ is a morphism of global actions, any such square must end up within
a single local patch of $\A$ and so can be lifted. ``Uniqueness''
ensures that it glues to any lifts constructed earlier in the
process in the obvious way. The only statement left unproved is
the last.

We end up with $\tilde x'_N$ and $\tilde x_N$ lying over $x_N$ and
the right hand side of the homotopy giving us a path from $\tilde
x_N$ to $\tilde x'_N$ which maps down (via $p$) to the identity
path from $x_N$ to itself. ``Uniqueness'' of path lifting then
shows this must be the identity path at $\tilde x_N$ and so
$\tilde x_N = \tilde x'_N$ as required. \end{proof}

\begin{cor}
If $p :\B  \rightarrow \A$ is a covering and $\A$ is connected
then all the fibres $p^{- 1} (x)$, $  x \in \A$ have the same
cardinality.
 \end{cor}
\begin{proof}
If $x_0, x_1 \in \A$, let $f :\LL \rightarrow \A$ be a path from
$x_0$ to $x_1$. Now pick $\tilde x_0 \in p^{- 1} (x_0)$, and lift
$f$ to $\tilde f$ joining $\tilde x_0$ to some uniquely determined
$\tilde x_1 \in p^{- 1} (x_1)$. This assignment is a bijection
since the reverse path is also uniquely determined.\end{proof}
\subsection{The $\pi_1(\A,a_0)$-action}
This gives a way of associating to each covering $p :\B
\rightarrow \A$, a set $F_{a_0} ( = p^{- 1} (a_0))$ together with
an action of $\pi_1 (\A, a_0)$. Alternatively the covering may be
thought of as a collection of fibres indexed by the elements of
$X_\A$ and then we get an action of $\Pi_1 \A$ on $\B $ by ``deck
transformations'' over $\A$.
\begin{thm}
Suppose $p :\B  \rightarrow \A$ is a covering map, and $b_0 \in
p^{- 1} (a_0)$ then the induced maps
$$p_{\ast} ; \Pi_1 \B  \rightarrow \Pi_1 \A$$
and
$$p_{\ast} :\pi_1 (\B , b_0) \rightarrow \pi_1 (\A, a_0)$$
are  monomorphisms.
\end{thm}
\begin{proof}
The induced maps take $\omega \in\Pi_1 \B $ with $w = [f]$ to
$p_{\ast} (\omega) = [pf] \in\Pi_1 \A$. The result is just an
immediate consequence of unique path lifting together with the
lifting of homotopies. (The proof is easy and follows that given
in many elementary homotopy texts.)\end{proof}

Now suppose $a_0 \in \A$ is the chosen base point, that $\B $ is
connected  and we have $b_0, b_1 \in p^{- 1} (a_0)$. Choose a path
class $\gamma$ from $b_0$ to $b_1$ in $\Pi_1 \B $ which thus gives
an isomorphism
$$u :\pi_1 (\B , b_0) \rightarrow \pi_1 (\B , b_1)$$
by conjugation $u (\omega) = \gamma^{- 1} \omega \gamma$ within
$\Pi_1 \B $. As in the topological case we get an inner
automorphism $\nu$ of $\pi_1 (\A, a_0), \nu (\omega) = p_{\ast}
(\gamma)^{- 1} \omega p_{\ast} (\gamma)$, and of course $p_{\ast}
(\gamma)$ is a loop since
  $p(b_0) = p(b_1)$. We thus have
\begin{prop}
The images of $\pi_1 (\B , b_0)$ and $\pi_1 (\B , b_1)$ are
conjugate subgroups of $\pi_1 (\A, a_0)$.\hfill
$\Box$\end{prop}

Path lifting then shows that any conjugate subgroup of $\pi_1 (\B
, b_0)$ in $\pi_1 (\A, a_0)$ can arise in this way. Just lift any
conjugating element to a path in $\B $.

 The theory  of coverings of global actions
  follows the same general development as the classical topological
  one (cf. Massey \cite{massey}) or the groupoid one (cf. Brown, \cite{rb}). For instance
one can easily prove the following results.

\begin{prop}
Let $p:\C\to\B$ and $q:\B\to\A$ be morphisms of global actions.
\begin{enumerate}
\item
If $p$ and $q$ are covering maps, so is $qp$.

\item
If $p$ and $qp$ are covering maps and $p$ is epi, then $q$ is a covering.
\end{enumerate}
\end{prop}

\begin{prop}
 If $p :\B \rightarrow \A$ and $q :\C \rightarrow \A$ are coverings and
  $f :\B  \rightarrow \C$ is a morphism over $\A$, then $f$ is also a
  covering.
\end{prop}

The proof of the following result is in fact easier than in the classical case. It uses the unique path lifting property of coverings.

\begin{prop}\label{clasico}
Let $p:\B\to\A$ be a covering and $a_0\in X_\A$, $b_0\in X_\B$ such that $p(b_0)=a_0$, let $f:\C\to\A$ be a morphism with $f(c_0)=a_0$ and suppose that $\C$ is connected. Then $f$ lifts
to a morphism $\tilde f:\C\to\B$, with $\tilde f(c_0)=b_0$, if and only if $f_*(\pi_1(\C,c_0))\subseteq p_*(\pi_1(\B,b_0))$.\endbox
\end{prop}

\begin{rem} Given the similarity of the development to the classical
  and groupoid cases, it should be clear that {\it all} of the above
  goes across to the context of groupoid atlases. There are also
  strong lifting properties for strong coverings. 
\end{rem}

We fix a base point $a_0$ in $\A$ and denote for simplicity $\pi_1\A=\pi_1 (\A, a_0)$.  

The action of $\pi_1\A$ on $p^{- 1} (a_0)$ extends to a
functor from the category of coverings over $\A$ to that of $\pi_1
\A$-sets, i.e., sets with $\pi_1 \A$-actions.
Explicitly if $b_0 \in p^{- 1} (a_0)$ and $\gamma \in\pi_1 (\A,
a_0), \gamma = [f]$ say; lift $\gamma$ to a path in $\B $ starting
at $b_0$. Its other end will be $b_0^{\gamma} \in p^{- 1} (a_0)$.

If $g :(\B , p) \rightarrow (C, q)$ is a morphism in the category
of coverings over $\A$, then $g$ restricts to a map $p^{- 1}(a_0)
\rightarrow q^{- 1} (a_0)$. Uniqueness of path lifting then shows
that
$$g(b_0^{\gamma}) = g(b_0)^{\gamma}$$
as hoped for.

Writing $\text{Cov}/\A$ for the category of coverings of $\A$, we
will get
$$ \text{Cov}/\A \rightarrow \pi_1 \A\text{-Sets}.$$
If $(\B , p)$ is a covering global action of $\A$, we will write
$\Aut_\A(\B , p)$ for its automorphism group (group of covering or
deck transformations) within $\text{Cov}/\A$.

The functor above gives a homomorphism
$$\Aut_\A (\B, p) \rightarrow \Aut_{\pi_1 \A\text{-Sets}}
(p^{- 1} (a_0))$$
as is easily checked.

If $\varphi :p^{- 1} (a_0) \rightarrow p^{- 1} (a_0)$ is an
automorphism of $\pi_1 \A$-sets then the isotropy subgroup of any
point $b_0 \in p^{ - 1} (a_0)$ is the same as that of $\varphi
(b_0)$, i.e., from our earlier discussion, it is easily seen to be
$p_{\ast} (\pi_1 (\B , b_0))$. Thus as a $\pi_1 \A$-set, $p^{- 1}
(a_0)$ is isomorphic to the ``coset space'' $\pi_1 \A/ p_{\ast}
\pi_1 \B $. By \ref{clasico} one has that the automorphism
$\varphi$ can be realised by a deck transformation and so $\Aut_\A
(\B , p)$ and $\Aut_{\pi_1
  \A\text{-Sets}} (p^{- 1} (a_0))$ are isomorphic.
 \begin{cor}
If $\pi_1 \B $ is trivial (i.e. the covering global action is
simply connected) then $\Aut_\A (\B , p) \cong \pi_1 (\A, a_0)$,
i.e. $(\B , p)$ is a universal covering.\hfill $\Box$\end{cor}

\subsection{The Galois-Poincar\'{e} theorem for global actions}
To complete the triple description of $\pi_1  \A$, we need to
show that a simply connected covering exists. In fact we will show
more, namely that given any conjugacy class of subgroups of $\pi_1
\A$, we can find a covering $(\B , p)$ corresponding to that
conjugacy class (i.e. $\{p_{\ast} \pi_1 (\B , p) :b \in p^{-
1} (a_0)\}$ gives exactly the given conjugacy class). This and
generalities on $\pi_1 \A$-sets will then establish that
$$\text{Cov}/_{\displaystyle \A} \rightarrow \pi_1 \A\text{-Sets}$$
is an equivalence of categories, which is the Galois-Poincar\'e
correspondence theorem in this context.

As before let $\A$ be connected and pick a base point $a_0 \in
X_\A$. Let $H$ be a subgroup of $\pi_1 (\A, a_0)$.

We have a set $\Gamma \A$ of based paths and a projection
$$ p :\Gamma \A \rightarrow \A$$
given by $p(\omega) = e^1 (\omega)$. Define an equivalence
relation $\sim_H$ on $\Gamma \A$ by $f \sim f'$ if $p(f) = p(f')$
and $[f] [f']^{- 1} \in H$, where the composition is viewed as
taking place within the fundamental groupoid $\Pi_1 \A$ and $H$ as
a subgroup of the vertex group at $a_0$.

Let $X_{\A_H}$ denote the set of equivalence classes, $\<f\>$, of
based paths under $\sim_H$.

The function  $p: X_{\Gamma \A} \rightarrow X_\A$ clearly induces
one, $p_H :X_{\A_H} \rightarrow X_\A$ given by $p_H \<f\> = e^1
(f)$.  We will give $\A_H$ a global action structure. Take
$$\Phi_{\A_H} = \Phi_\A$$

For $\alpha \in \Phi_{\A_H},$  $(X_{\A_H})_{\alpha} = \{ \omega
\in X_{\A_H} ~|~ p_H (\omega) \in (X_\A)_{\alpha}\}$ and
$(G_{\A_H})_{\alpha} = (G_\A)_{\alpha}$. The action of
$(G_{\A_H})_{\alpha}$ on $(X_{\A_H})_{\alpha}$ is as follows:

Let $f \in \Gamma \A$ ( and we as usual assume $N^-_f = 0$ and
$N^+_f = n$, say), then $p(f) = f(n)$. Suppose $f \in
(X_{\A_H})_{\alpha}$, so $f(n) \in (X_{\A})_{\alpha}$ and let
$\sigma \in (G_\A)_{\alpha}$. Define a path $\sigma^f_n$ by
\[
\sigma^f_n (m) =
\begin{cases}
f(n) & \text{if } m \leq n,\\
\sigma. f(n) & \text{if } m \geq n + 1
\end{cases}
\]
and set
$$\sigma. f = f \ast \sigma^f_n. $$
Thus
\[
\sigma. f(m) =
\begin{cases}
f(m) & \text{if } m \leq n\\
\sigma. f(n) & \text{if } m \geq n+1
\end{cases}
\]
It is easily checked that this gives an action on equivalence
classes by
$$\sigma. \<f\> = \<\sigma. f\>$$
since it just adds one extra ``link'' to the path. Clearly $p_H
\sigma. \<f\> = \sigma. p_H \<f\> = \sigma. f(n)$, $\A_H$ is a
global action and $p_H$ a regular morphism of global actions.

We can now prove that $(\A_H,p_H)$ is a covering of $\A$.

Suppose that $x_0,\ldots,x_n$ is a local frame in $\A$ and that $\omega_0=\<f_0\>$ is an element in 
$\A_H$ such that $p_H(\omega_0)=x_0$. Since $x_0,\ldots,x_n$ is a local frame, there exists $\alpha\in\Phi_{\A}=\Phi_{\A_H}$ and
$g_i\in G_{\alpha}$ such that $x_i=g_ix_0$. For each $i$ take $\omega_i=g_i\omega_0$. It is clear by definition, that $\omega_0,\ldots,\omega_n$ is
a local frame in $\A_H$ and that $p_H(\omega_i)=x_i$. Moreover, it is the unique local frame  with this property. 
This proves that $(\A_H,p_H)$ is a covering.

We note that $p^{- 1}_H (a_0)$ is the set of $\sim_{H}$ equivalence
classes of loops at $a_0$. If we look at the equivalence relation
$\sim_H$, it is clearly made up of two parts:

\begin{enumerate}[(i)]
\item if $f \sim f'$ (i.e. fixed end point homotopic) then
  clearly $[f] [f']^{- 1} = 1_{a_0} \in\Pi_1 \A$, and so $f \sim_H
  f'$ whatever $H$ is chosen;
\item if $f$ is a loop at $a_0$ then $f \sim_H  1_{a_0}$ if $[f]
  \in H$.
\end{enumerate}

Together these imply that $p^{- 1} (a_0) \cong G/H$ where we have
written $G$ for $\pi_1 (\A, a_0)$.

Since  $(\A_H,p_H)$ is a covering of $\A$, there is a $G$-action on $p^{- 1} (a_0)$ making this a $G$-set
isomorphism. 

Any path $\varphi$ at $a_0$ in $\A$ will lift to a path given a
choice of initial point. Fix $\tilde a_0$ to be the class of the
constant path at $a_0$, so $\tilde a_0 \in \A_H$ and $p_H (\tilde
a_0) = a_0$. If $\tilde \varphi$ is the lift of $\varphi$ starting
at $\tilde a_0$ then $\tilde \varphi (n)$ is the element of $\A_H$
represented by the partial path from $\varphi (0)$ to $\varphi
(n)$. If $n > N^+_{\varphi}$ and $\varphi$ is a loop at $a_0$
representing an element of $H \leq \pi_1 (\A, a_0)$ then the end
point of $\tilde \varphi$ is the point of $\A_H$ represented by
the path from $\varphi (0)$ to $\varphi (n)$, i.e., $\<\varphi\>$,
but $[\varphi] \in H$ so $\<\varphi\> = \<\tilde a_0\>$ and
$\varphi$ lifts to a loop in $\A_H$. Conversely any loop in $\A_H$
at $a_0$ is the lift of a loop at $a_0$ which represents an
element of $H$.

A similar argument implies that if $\bar \varphi$ and $\bar
\varphi'$ are homotopic loops at $\tilde a_0$ in $\A_H$ then they
are lifts of homotopic loops at $a_0$ in $\A$ which then of course
represent the same element of $H$. We thus have
\begin{prop}
The induced homomorphism
$$p_{H \ast} :\pi_1 (\A_H, \tilde a_0) \rightarrow \pi_1 (\A, a_0)$$
is a monomorphism with image, $H$.
\end{prop}
 \begin{proof}
The proof is by direct calculation using the explicitly defined
lifts of paths and homotopies.\end{proof}

\begin{rem}
Much of the above would work for strong paths, but the proof that
the strong version of $(\A_H, p_H)$ is a covering would seem to
depend on a local condition which is in some way analogous to
``semi locally simply connected''. This would say that small
strong loops were strongly null-homotopic. Here by small we mean
\[
\xymatrix{
a \ar@{>}[r]^{g} &  b \ar@{>}[r]^{g'} & a.
}
\]
This is clearly satisfied for many examples.
\end{rem}

\begin{question}
Adapt the above discussion to handle strong
coverings and /or groupoid atlases.
\end{question}

To summarise:
\begin{thm}
Given $(\A, a_0), \A$ connected, and $H \leq \pi_1 (\A, a_0)$ then
there is a connected covering space $(\A_H, p_H)$ with
$$p_{H_{\ast}} (\pi_1 (\A_H, \tilde a_0)) = H.$$
In particular corresponding to $H = 1$,  the trivial subgroup of
$\pi_1 (\A, a_0)$, one has a simply connected covering space
$(\tilde \A,p).$ \hfill $ \Box$\end{thm}

 Of course by an earlier result $\Aut_{\A} (\tilde \A, p) \cong \pi_1
(\A, a_0)$.

Now set $G = \pi_1 (\A, a_0)$. Any $G$-set $X$ can be decomposed
as a disjoint union of ``connected'' $G$-sets. Here ``connected''
merely means single orbit or transitive $G$-sets. These all have
form $G/H$ and  we can note that $p_H^{- 1} (a_0) \cong G/H$ as
$G$-sets. Using disjoint unions of $(\A_H, p_H)$s for various
subgroups $H$ will yield a covering space $(\B , p)$ with $p^{- 1}
(a_0) \cong X$ as $G$-sets. It is then more or less routine to
check that
$$\text{Cov}/\A  \leftrightarrows G\text{-Sets}$$
is an equivalence of categories.

Thus we have three descriptions of $\pi_1 (\A, a_0)$ for a
connected global action $\A$:
\begin{enumerate}
\item[(i)] equivalence (homotopy) classes of loops at $a_0$,
\item[(ii)] $\pi_0 (\Omega \A)$, \item[(iii)] $\Aut_{\A}(\tilde \A, p)$
and  thus the group $G$ in the above   equivalence.
\end{enumerate}

 \section{Single domain global actions II.}

In this section we will examine general single domain global
actions and their coverings.
\subsection{General single domain global actions and $\A(G,\mathcal{H})$s}
Suppose $\A = (X_\A, \Phi_\A, G_\A)$ is a single domain global
action. We thus have that $\Phi_\A$ is a set with a reflexive
relation $\leq$ defined on it, then
$$G_\A :\Phi_\A \rightarrow \Groups$$
can be considered as a (generalised) functor and we can form its
colimit $G = \text{colim} ~G_\A$. Each $G_{\alpha}, \alpha \in
\Phi_\A$ acts on $G$ by left multiplication via its image in $G$.
(Note: $G_{\alpha}$ need not be isomorphic to a subgroup of $G$,
but other than that one has virtually the situation of $\A(G,
\mathcal{H})$.)

Define a global action $\mathsf{G}$ with $|\mathsf{G}| = G$
\begin{align*}
\Phi_{\mathsf{G}} & = \Phi_\A\\
(G_{\mathsf{G}})_{\alpha} & = G_{\alpha}\\
(X_\mathsf{G})_{\alpha} & = G,
\end{align*}
so $\mathsf{G}$ is a single domain global action.

If $H$ is any subgroup of $G$, we can form a quotient global
action $\mathsf{G}/H$ with $|\mathsf{G}/H| = |\mathsf{G}|/H$, the set of
right cosets $\Phi_{\mathsf{G}/H} = \Phi_{\mathsf{G}}$, etc, so
$\mathsf{G}/H$ is again a single domain global action.

It is clear that
$$\pi_0 (\mathsf{G}) = \pi_0 (\mathsf{G}/H) = 1.$$

\begin{thm}
Any connected single domain global action is regularly isomorphic
to some $\mathsf{G}/H$.\end{thm}
\begin{proof}
Let $\A$ be a connected single domain global action and let $a_0
\in \A$ be a chosen basepoint. Let $\mathsf{G}_\A$ be the global
action constructed above from colim $(G_\A :\Phi_\A \rightarrow
\Groups)$.

The group $G$ acts on $\A$ and also, of course, on
$|\mathsf{G}_\A|$. Define a function
$$p :|\mathsf{G}_\A| \rightarrow |\A|$$
(using the base point) sending $\omega \in G_\A$ to $\omega. a_0$,
i.e., read the word $\omega$ off from the right acting on $a_0$,
inductively. The possible ambiguities in the word are due to cases
of $\alpha \leq \beta$ and the compatibility condition ensures
this does not matter.

Since $g_{\alpha} \cdot \omega$ gets sent to $(g_{\alpha} \cdot \omega
a_0) = g_{\alpha} (\omega. a_0)$, this defines a regular morphism of
global actions.

Let $H_\A = p^{- 1} (a_0)$, which is the stabiliser of $a_0$ in
$\mathsf{G}_\A$.

Clearly
$$\mathsf{G}_\A/{H_\A} \cong \A.$$\end{proof}

The only difference therefore between single domain global actions
of the form $\A((G, K),\mathcal{H})$ as introduced in section \ref{AGKH}
and the general case is that the $H_i$ may not be subgroups and
may have interrelations between them.

\begin{examples}
\hspace{0.5em}
\begin{enumerate}[1.]

\item As before take $S_3 = \langle a,b  ~| ~a^3 = b^2 = (ab)^2 = 1
\rangle$, $H_1 = \langle a \rangle$, $H_2 = \langle b \rangle$ to
get $\A(S_3, \{\langle a \rangle,\langle b \rangle \})$.  Then the
colimit group is $C_3 \ast C_2$ and, of course, the stabiliser of
1 in $\mathsf{G}_\A$ is merely $Ker( C_3 \ast C_2 \rightarrow
S_3)$, that is $\pi_1(\A(S_3, \{\langle a \rangle,\langle b
\rangle \})$.   The quotient map $$p: C_3 \ast C_2 \rightarrow
S_3$$ is that with kernel the normal closure of $(ab)^2 $ and
writing $K$ for that kernel (and thus for $\pi_1$), we have
$$\A(S_3, \{\langle a\rangle,\langle b \rangle \})\cong \A((C_3 \ast C_2, K) \{\langle a
\rangle,\langle b \rangle \}).$$

\medskip

A similar picture emerges with the other examples.

\item For $\A_{K4} = \A(K_4, \{\langle a\rangle,\langle b
\rangle,\langle c \rangle\})$, the colimit group is $C_2\ast
C_2\ast C_2 = C_2^{(3)}$ and the stabiliser is the normal closure
of $abc$. This normal subgroup has rank 3. Thus $\A_{K4}$ has a
second description as $\mathcal{\A}((C_2^{(3)},K),\mathcal{H})$,
where $\mathcal{H} = \{\langle a \rangle, \langle b \rangle,
\langle c \rangle\}$, these subgroups being here subgroups of
$C_2^{(3)}$, not of $K_4$.  Of course, $K \cong \pi_1(\A_{K4})$.

\item  The only change for $\A_{q8}$ is that the colimit group is
$(C_4\ast C_4\ast C_4) /\{1,-1\}$, the free product with
amalgamation.

\item Taking $S_3$ again, but with presentation $\langle x_1, x_2 ~|~
x_1^2 = x_2^2= 1, (x_1x_2)^3 = 1\rangle,$ $H_1 = \langle x_1
\rangle$, $H_2 = \langle x_2 \rangle$,  gives colimit group
$C_2\ast C_2$.  The stabiliser of 1 / fundamental group is free on
$(x_1x_2)^3$.  We again get a second description as a `relative'
$\mathcal {\A}((G,K), \mathcal{H})$.

\item For our final example, $S_4$ with presentation $\langle
x_1,x_2, x_3 ~|~ x_1^2, i = 1,2,3, (x_1x_2)^3 = 1 =(x_2x_3)^3,
(x_1x_3)^2 = 1\rangle $ is $S_4$ itself and the given description
is the one we have found earlier through the general process.

\end{enumerate}
\end{examples}

\subsection{Coverings of single domain global actions}

If $\A \cong \mathsf{G}_\A/{H_\A}$ as above then for $\tilde \A$,
its universal or simply connected covering, $\tilde \A$ is also a
single domain global action and as the diagram $G_{\tilde \A} =
G_\A$,
$$\tilde \A \cong \mathsf{G}_\A/{H_{\tilde \A}}.$$
Clearly there is a diagram
\[
\xymatrix{ \mathsf{G}_{\tilde \A}/{H_{\tilde \A}} \ar@{>}[d]
\ar@{>}[r]^{\quad\cong}&
\tilde \A \ar@{>}[d]\\
\mathsf{G}_\A/{H_\A} \ar@{>}[r]^{\quad\cong} & \A }
\]
so $H_{\tilde \A} \subseteq H_\A$.

It remains to relate $H_{\tilde \A}$ more closely to $H_\A$.

If $\alpha \in \Phi_{\tilde \A} = \Phi_\A$, then set
$(G'_\A)_{\alpha} = \text{ image} \left( (G_\A)_{\alpha}
\rightarrow G \right)$.

Let
\begin{align*}
\mathcal{H}_A = \{ H \leq |{G}_\A| ~|  ~\text{ for all }
\sigma  \in |{G}_\A|, ~\sigma H \sigma^{- 1} \cap
(G'_\A)_{\alpha} = ~& \sigma H_\A \sigma^{- 1} \cap
(G'_\A)_{\alpha}
 \\ &\text{ \hspace{1cm}  for all } \alpha \in \Phi \}
\end{align*}
and let $H_{\tilde \A} = \bigcap\{H \in \mathcal{H}\}$. Then $
H_{\tilde \A} \in \mathcal{H}_\A $ and is minimal.
 Of course $H_{\tilde \A} \triangleleft H_\A$ and $H_\A/{H_{\tilde \A}} \cong \pi_1
 (\A)$.  Thus as a corollary of the main classification theorem for coverings we get:
\begin{prop} General connected coverings of $\A$ correspond  bijectively to intermediate
 groups  between $H_{\tilde \A}$ and $H_\A$.\endbox\end{prop}
\begin{rem}
If $\pi_0(\A)$ is not trivial,  i.e. $\A$ is not connected, then
it is important to remember that $\tilde{\A}$ only covers the
connected component of the basepoint. 
\end{rem}
Again we turn to our examples to see what these results give there.

\medskip

In each case we have a description of $\A$ as
$\A((G_\A,K_\A),\mathcal{H})$ and
 as our initial situation had
$$K_\A = \Ker(p: G \to X_\A),$$
where $G$ is the colimit group of the original system, we have
$$K_\A \cong \pi_1(\A,1)$$in each case.  This implies that $K_{\tilde \A}$ is
trivial.  We could also deduce this from the description of
$\mathcal{H}_\A$ with $K_{\tilde{\A}} = \bigcap \mathcal{H}_\A$.
($K_\A$ is normal in $G_\A$ so $ \sigma K_\A\sigma^{-1} \cap
(G^\prime_\A)_\alpha = K_\A \cap (G^\prime_\A)_\alpha $.
 In each example this is trivial, so $\mathcal{H}_\A$ contains the trivial
 group and hence has that group as its intersection.)

We can describe the simply connected covering of $\A$ in each
case:

1. $\A = \mathcal{\A}(S_3,\{\langle a\rangle, \langle
b\rangle\})$, $\tilde{\A} = \mathcal{\A}(C_3\ast C_2, \{\langle
a\rangle, \langle b\rangle\})$, where as before $\langle a \rangle
$ is to be interpreted \emph{in context} as a subgroup of the
corresponding group.

2.  $\A_{K4} = \mathcal{\A}(K_4, \{\langle a\rangle, \langle
b\rangle,\langle c\rangle\})$, $\tilde{\A_{K4}} =
\mathcal{\A}(C_3^{(3)}, \{\langle a\rangle, \langle
b\rangle,\langle c\rangle\})$.

3. $\A_{q8} = \mathcal{\A}(q8, \{\langle i\rangle, \langle
j\rangle,\langle k\rangle\})$, $\tilde{\A_{q8}} =
\mathcal{\A}(G_\A, \{\langle i\rangle, \langle j\rangle,\langle
k\rangle\})$, where $G_\A$ has  presentation
$$\langle i,j,k ~|~ i^4 = 1, i^2 = j^2 = k^2 \rangle.$$

4.  $\A = \mathcal{\A}(S_3, \mathcal{H})$, with $\mathcal{ H} = \{
H_1,H_2\}$ with each $H_i$ generated by a transposition,
$\tilde{\A} = \mathcal{\A}(C_2^{(2)}, \{\langle x_1\rangle,
\{\langle x_2\rangle\}).$

5.  The case $\mathcal{\A}(S_4, \mathcal{H})$ with $\mathcal{ H} =
\{H_1,H_2,H_3\}$, as above (example \ref{S4}), is already simply connected, so is its
own simply connected cover.

\medskip

In the next section we examine a more complex example, namely the
elementary matrix group, $\E_n(R)$ of a ring $R$, which forms the
connected component of the identity in the global action
$\Gla_n(R)$.

\section{The Steinberg Group and Coverings of  $\Gla_n (R)$ }

A particularly important example of a single domain global action
is the General Linear Global Action $\Gla_n (R)$.  We saw,
example \ref{GLn}, that $\pi_0(\Gla_n (R))$ was the set of right
cosets of the group $\Gl_n(R)$ modulo the subgroup $\E_n(R)$ of
elementary matrices and hence was identifiable as being
$K_1(n,R)$.  We asked ``\emph{is $K_2(n,R)
\cong \pi_0(\Gla_n (R))$?}''  It is to this question we now turn.
\subsection{The Steinberg group $\St_n(R)$}
The usual approach to $K_2(n,R)$ is via the Steinberg group $\St_n(R)$.  We
earlier, example \ref{GLn}, introduced  the notation $\epsilon_{ij}(r)$ for the elementary
matrix with
\[
 \epsilon_{ij} (r)_{k, l} = \begin{cases}
1 & \text{if } k = l,\\
r & \text{if } (k, l) = (i, j),\\
0 & \text{otherwise .}
 \end{cases}
\]
Here, of course, $(i,j) \in \Delta$, the set of non-diagonal
positions in an $n\times n$ array. Elementary matrices satisfy
certain standard relations and the Steinberg group is obtained by
considering the group having generators $x_{ij}(r)$, abstracting
the elementary matrices, and having as relations just these
standard, almost universal, relations.  More precisely (and a
standard reference is Milnor's notes, \cite{Milnor}), $\St_n(R)$
is given by generators $x_{ij}(r)$, $r\in R$, $i,j = 1,2, \ldots,
n$, $i \neq j$, which are subject to the relations:

St1 \quad $x_{i,j}(a)x_{i,j}(b) =  x_{i,j}(a + b)$;

St2 \quad $[ x_{i,j}(a),x_{k,\ell }(b)] = \left\{ \begin{array}{ll}
      1 & \textrm{if } i \neq \ell, j\neq k\\
      x_{i,\ell}(ab) & i \neq \ell, j = k
\end{array}\right.$

These are called the \emph{Steinberg relations}.

There is an epimorphism
$$\varphi :\St_n(R) \rightarrow \E_n(R)$$given by mapping $x_{i,j}(a)$ to
$\epsilon_{i,j}(a)$.  The second (unstable) $K$-group, $K_2(n, R)$ is then
defined to be $Ker~\varphi$.
\subsection{A more detailed look at $\Gla_n(R)$}We will construct a global action analoguous to $St_n(R)$ but for this we need to understand $\Gla_n(R)$ better.
Early in this paper, \S\ref{GLn}, p.\pageref{GLn}, we
introduced the elementary matrix groups, $\Gl_n (R)_\alpha$.
Recall that we let $\Delta$ be the set of off-diagonal positions in an
$n\times n$ array and called a subset $\alpha \subseteq \Delta$
\emph{closed} if it corresponded to a transitive relation, i.e. if $(i,j)
\in \alpha$ and $(j,k) \in \alpha$, then $(i,k)\in \alpha$.
The general linear global action $\Gla_n(R)$ then had coordinate
system $\Phi$, the set of closed subsets of $\Delta$ ordered by
inclusion.  The underlying set of $\Gla_n(R)$ was the general
linear group ${ \Gl}_n(R)$ and for $\alpha \in \Phi$,
$\Gl_n(R)_\alpha$ was the group of elementary matrices generated
by the $\epsilon_{i,j}(r)$ with $(i,j)\in \alpha$. 

Clearly this
single domain global action is of the form $\A(G,\cal{H})$. To
form its connected covering (which will cover the connected
component of 1, that is, will cover the sub-global action, $\E_n(R)$, determined by the elementary matrices), we need to take the
colimit of the $\Gl_n(R)_\alpha$. Clearly to examine this colimit
we need to see what the maximal elements of $\Phi$ are, and to
examine the corresponding $\Gl_n(R)_\alpha$.

\begin{lem}
If $\alpha \in \Phi$ is maximal, then it is a total order on $\{1,
\ldots, n\}$.\end{lem}
 \begin{proof}
Transitivity follows from closedness.  If $(i,j)\in \alpha$, then
$(j,i)\notin \alpha$, since no diagonal elements are in $\alpha$,
but as $\alpha$ is maximal, one or other of  $(i,j)$ and $(j,i)$
must be in $\alpha$ - otherwise we could add it in! \end{proof}

\begin{lem}
Let $T_n(R)$ be the group of upper triangular $n \times n$ matrices over $R$.
If $\alpha \in \Phi $ is maximal, then
$$\Gl_n(R)_\alpha \cong T_n(R).$$
\end{lem}
 \begin{proof}
 Pick an order isomorphism, $f$ between $\alpha$ and the total order $1 < 2 <
\ldots < n$.   Map the generator $\epsilon_{i,j}(r)$ to
$\epsilon_{f(i),f(j)}(r)$.  This extends to the required
isomorphism from $\Gl_n(R)_\alpha$  to $T_n(R)$.\end{proof}
\subsection{The Steinberg global action $\mathsf{St}_n(R)$}
It is known (cf. Milnor's notes \cite{Milnor}) that $T_n(R)$ has a
presentation given by the $x_{k,\ell}(r)$, with $1 \leq k < \ell
\leq n$, and with the Steinberg relation (restricted to those
indices $(k,\ell)$ with $k < \ell$) between them.  Let
$\St_n(R)_\alpha$ be the group given by generators $x_{k,\ell}(r)$
with $(i,j)\in \alpha$ and with the corresponding Steinberg
relations, then
$$colim_{\alpha \in \Phi} \St_n(R)_\alpha = \St_n(R),$$ so if we define a global action as below, it will be connected.
\begin{Def}
Let $\mathsf{St}_n(R)$ be the global action having $ \St_n(R)$ as its underlying set, $\Delta$, above, as its coordinate system and, for $\alpha\in \Delta$, $\St_n(R)_\alpha$ as the corresponding local group.
\end{Def} 
The isomorphism $\Gl_n(R)_\alpha \cong T_n(R)$ for maximal $\alpha$
together with the fact that $T_n(R) = \St_n(R)_{\alpha_0}$ for
$\alpha_0 = \{ (i,j) :i < j\}$  gives that there is an
isomorphism $\varphi_\alpha :\St_n(R)_\alpha
\stackrel{\cong}{\rightarrow} \Gl_n(R)_\alpha$ for maximal
$\alpha$ compatible with the inclusions into $\St_n(R)$ and
$\Gl_n(R)$ and the homomorphism $\varphi$ introduced earlier.

Two maximal $\alpha$ can be linked with each other by a zig-zag where
intermediate maximal elements (total orders) differ by the transposition of
two elements only.  The corresponding  isomorphisms $\Phi_\alpha$ agree on
intersections of the corresponding groups thus giving an isomorphism
$$\St_n(R) = colim_{\alpha \in \Phi} \St_n(R)_\alpha\stackrel{\cong}{\rightarrow}  colim_{\alpha \in \Phi}
\Gl_n(R)_\alpha = \widetilde{\Gl_n(R)}.$$
The resulting map is then
easily shown to be $\varphi$. It remains to analyse this map a little more.

The kernel,
$$H_\A = \Ker(\St_n(R)\simeq colim_{\alpha \in \Phi}
  \Gl_n(R)_\alpha\rightarrow \E_n(R)),$$
is central (again see Milnor's notes), hence
$$\sigma H_\A\sigma^{-1} \cap Im \Gl_n(R)_\alpha = H_\A \cap 
  \Gl_n(R)_\alpha = \{1\}$$as no element of $
  \Gl_n(R)_\alpha $ vanishes when mapped into $
  \Gl_n(R)$ as the mapping is an inclusion.  Thus the family, whose minimal
element we need, contains the trivial subgroup!  Hence that must be the
minimal element.
\end{proof}

We have proved
\begin{thm}
The simply connected universal covering of $\Gla_n(R)$ is isomorphic to  $\mathsf{St}_n(R)$.  The covering map is given on elements by  the evaluation mapping
$$\varphi :\St_n(R) \rightarrow \Gl_n(R).$$\hfill $ \Box$
\end{thm}
\begin{cor}
The second $K$-group $K_2(n,R)$ is isomorphic to $\pi_1(\Gla_n(R))$.\hfill $ \Box$
\end{cor}


\begin{thebibliography}{10}
\newcommand{\enquote}[1]{`#1'}

\bibitem{AbelsHolz}
\textsc{H.~Abels} and \textsc{S.Holz}, \enquote{Higher generation by
  subgroups.}
\newblock \textit{J. Alg} 160 (1993) 311-- 341.

\bibitem{bak1}
\textsc{A.~Bak}, \enquote{Global actions: The algebraic counterpart of a
  topological space.}
\newblock \textit{Uspeki Mat. Nauk., English Translation: Russian Math.
  Surveys} 525 (1997) 955 --996.

\bibitem{bak2}
\textsc{A.~Bak}, \enquote{Topological methods in algebra.}
\newblock \enquote{Rings, Hopf Algebras and Brauer Groups,}  (eds
  \textsc{S.~Caenepeel} and \textsc{A.~Verschoren}) (M. Dekker, New York,
  1998), no. 197 in Lect. Notes in Pure and Applied Math .

\bibitem{bakbmp99}
\textsc{A.~Bak}, \textsc{R.~Brown}, \textsc{G.~Minian} and \textsc{T.~Porter},
  \enquote{Global actions, groupoid atlases and related topics.}
\newblock \textit{University of Wales, Bangor, Mathematics Preprint} 99.27.

\bibitem{Behr}
\textsc{H.~Behr}, \enquote{Explizite pr\"asentationen von chevalleygruppen
  \"uber $\mathbb{Z}$.}
\newblock \textit{Math. Z.} 141 (1975) 235--241.

\bibitem{brid&haef}
\textsc{M.~Bridson} and \textsc{A.~Haefliger}, \textit{Metric Spaces of
  Non-Positive Curvature}.
\newblock No.~31 in Grundlehren der Math. Wiss (Springer, 1999).

\bibitem{rb}
\textsc{R.~Brown}, \textit{Topology and Groupoids} (Booksurge LLC, S. Carolina,
  2006).

\bibitem{corson1}
\textsc{J.~M. Corson}, \enquote{Complexes of groups.}
\newblock \textit{Proc. London Math. Soc.} 65 (1992) 199--224.

\bibitem{corson2}
\textsc{J.~M. Corson}, \enquote{Groups acting on complexes and complexes of
  groups.}
\newblock \enquote{Geometric Group Theory,}  (eds \textsc{Davis Charney} and
  \textsc{Shapiro}) (de Gruyter, Berlin, 1995) pp. 79 --97, pp. 79 --97.

\bibitem{corson3}
\textsc{J.~M. Corson}, \enquote{Howie diagrams and complexes of groups.}
\newblock \textit{Comm. Algebra} 23 (1995) 5221 -- 5242.

\bibitem{Dowker}
\textsc{C.~H. Dowker}, \enquote{Homology groups of relations.}
\newblock \textit{Annals of Maths} 56 (1952) 84--95.

\bibitem{haef1}
\textsc{A.~Haefliger}, \enquote{Complexes of groups and orbihedra.}
\newblock \enquote{Group Theory from a Geometric viewpoint,} ICTP, Trieste, 26
  March- 6 April 1990 (World Scientific, 1991) pp. 504 -- 540, pp. 504 -- 540.

\bibitem{haef2}
\textsc{A.~Haefliger}, \enquote{Extensions of complexes of groups.}
\newblock \textit{Annales Inst.Fourier} 42 (1992) 275 -- 311.

\bibitem{massey}
\textsc{W.~M. Massey}, \textit{Algebraic Topology, an introduction} (Harcourt,
  Brace \& World, 1967).

\bibitem{Milnor}
\textsc{J.~Milnor}, \textit{Introduction to Algebraic K-theory}, vol.~72 of
  \textit{Annals of Mathematics Sudies} (Princeton, 1971).

\bibitem{Soule}
\textsc{C.~Soul\'e}, \enquote{Groupes op\'erants sur un complexe simplicial
  avec domain fondamental.}
\newblock \textit{C.R. Acad. Sci. Paris} S\'er A 276 (1973) 607--609.

\bibitem{Tits}
\textsc{J.~Tits}, \enquote{Ensembles ordonn\'es, immeubles et sommes
  amalgam\'es.}
\newblock \textit{Bull. Soc. Math. Bel.} S\'er A 38 (1986) 367--387.

\bibitem{volodin}
\textsc{I.A. Volodin}, \enquote{Algebraic k-theory as extraordinary homology
  theory on the category of associiative rings with unity.}
\newblock \textit{Izv. Akad. Nauk. SSSR} 35.
\newblock (Translation: Math. USSR Izvestija Vol. 5 (1971) No. 4, 859-887).

\end{thebibliography}

 
 \end{document}